\documentclass[10pt]{article}
\usepackage{geometry}
\usepackage[T1]{fontenc}
\usepackage{charter}
\usepackage{amsmath,amssymb,amsfonts,amsopn,amsthm}
\usepackage{mathtools}
\usepackage{dsfont}
\usepackage{bm}
\usepackage{bbm}
\usepackage{nicefrac}
\usepackage{mathrsfs}
\usepackage{graphicx}
\usepackage{epstopdf}
\usepackage{caption}
\usepackage{subcaption}         
\usepackage{overpic}
\usepackage{float}
\usepackage{chngcntr}
\counterwithout{figure}{section}
\usepackage[section]{placeins}
\usepackage{esint}
\usepackage{tikz}
\usepackage{pgfplots}
\pgfplotsset{compat=1.12}
\usepackage{tikz-3dplot}
\usetikzlibrary{%
	patterns, shapes, positioning, shadings, calc,
	decorations.pathreplacing, decorations.pathmorphing,
	angles, quotes, backgrounds, intersections, plotmarks}
\tikzset{snake it/.style={decorate, decoration=snake}}
\usepgfplotslibrary{fillbetween}
\usepackage{etoolbox}
\apptocmd{\thebibliography}{%
	\setlength{\itemsep}{-1pt}%
	\setlength{\parsep}{0pt}%
	\setlength{\parskip}{0pt}%
}{}{}
\usepackage{titlesec}
\titleformat{\section}
{\large\bfseries}{\thesection.}{0.5em}{}
\titleformat{\subsection}
{\normalsize\bfseries}{\thesubsection.}{0.4em}{}
\titleformat{\subsubsection}
{\normalsize\itshape}{\thesubsubsection.}{0.3em}{}
\titlespacing*{\section}{0pt}{1.1ex plus .3ex}{0.6ex}
\titlespacing*{\subsection}{0pt}{0.8ex plus .2ex}{0.4ex}
\titlespacing*{\subsubsection}{0pt}{0.6ex plus .2ex}{0.3ex}
\usepackage[numbers,sort&compress]{natbib}
\usepackage{hyperref}
\hypersetup{
	colorlinks = true,
	linkcolor  = blue!70,
	citecolor  = red,
	urlcolor   = red!20,
	filecolor  = red!20,
}
\usepackage{enumitem}
\setlist{
	topsep=2pt,
	itemsep=2pt,
	parsep=0pt,
	rightmargin=0pt
}
\usepackage{longtable}
\usepackage{supertabular}
\usepackage{siunitx}
\usepackage{textcomp}
\usepackage{url}
\usepackage{ifthen}
\usepackage[normalem]{ulem}     
\usepackage{listings}
\usepackage{authblk}

\definecolor{navyblue}{rgb}{0.0,0.0,0.5}
\definecolor{skyblue}{RGB}{135,206,235}
\definecolor{deepskyblue}{RGB}{0,191,255}

\numberwithin{equation}{section}
\newtheorem{theorem}{Theorem}[section]
\newtheorem{lemma}[theorem]{Lemma}

\newtheorem{proposition}[theorem]{Proposition}

\newtheorem{definition}[theorem]{Definition}

\newtheorem{remark}[theorem]{Remark}

\def\md{{\mathrm{d}}}

\def\O{\Omega}

\def\e{\varepsilon}

\def\wo{\overline}

\def\fU{{\mathfrak{U}}}

\def\Te{{\cal T}_\e}

\def\fC{\mathfrak C}

\def\fW{\mathfrak{W}}

\def\Dc{\mathcal{D}}

\def\N{{\mathbb{N}}}

\def\R{{\mathbb{R}}}

\def\A{{\mathbb{A}}}
\def\D{{\mathbb{D}}}

\def\G{{\mathbb{G}}}

\def\Ec{{\mathcal{E}}}

\def\Kc{{\mathcal{K}}}

\def\Rc{{\mathcal{R}}}

\def\Vc{{\mathcal{V}}}

\def\wt{\widetilde}
\def\X{\times}
\def\div{{\mathrm{div}}}

\def\Gg{{\bf g}}

\def\Gn{{\bf n}}
\def\Gu{{\bf u}}
\def\Gv{{\bf v}}\def\Gw{{\bf w}}

\def\GH{{\bf H}}
\def\GL{{\bf L}}\def\GO{{\bf O}}

\def\GT{{\bf T}}

\def\bv{\mathbf{v}}

\def\bw{\mathbf{w}}
\def\per{\mathrm{per}}

\def\a{\alpha}
\def\bn{\mathbf{n}}

\def\bpsi{{\boldsymbol{\psi}}}
\begin{document}
	
	%
	
	\title{Navier-Stokes-Cahn-Hilliard system in a $3$D perforated domain with free slip and source term: Existence and homogenization}

	\author{%
			Amartya Chakrabortty\thanks{
			Department of Mathematics, IIT Kharagpur,  721302, Kharagpur
			India; \texttt{amartyat26@kgpian.iitkgp.ac.in}},~~ 
		Haradhan Dutta\thanks{Department of Mathematics, IIT Kharagpur, 721302, Kharagpur India; \texttt{hdutta1412@kgpian.iitkgp.ac.in}},~~
		Hari Shankar Mahato\thanks{Department of Mathematics, IIT Kharagpur, 721302, Kharagpur India; \texttt{hsmahato@maths.iitkgp.ac.in}}%
	}

	
	\date{} 
	
	\newcommand{\keywords}[1]{\par\smallskip\noindent\textbf{Keywords: }#1}
	\newcommand{\subjclass}[1]{\par\noindent\textbf{MSC 2020: }#1}
	
	\maketitle
	
	\vspace{-1.0cm}
\begin{abstract}
	We study a time-dependent Navier--Stokes--Cahn--Hilliard system for a
	binary incompressible mixture in a periodically perforated domain
	$\Omega_p^\varepsilon\subset\mathbb{R}^3$. The obstacles have diameter
	of order $\varepsilon^\alpha$, with $\alpha>3$, and mutual distance of
	order $\varepsilon$. Hence
	$\varepsilon^\alpha/\varepsilon^3\to0$, which corresponds to the
	subcritical dilute regime. The model includes a periodically oscillating
	viscosity tensor, a nonconservative source term in the Cahn--Hilliard
	equation, no-slip conditions on the outer boundary, and free-slip
	conditions on the obstacle surfaces. The capillary coefficient
	$\lambda^\varepsilon>0$ depends on $\varepsilon$.
	
	For every fixed $\varepsilon>0$, we prove existence of a weak solution
	and derive estimates uniform in $\varepsilon$ with explicit
	$\lambda^\varepsilon$-scaling. Assuming
	$\lambda^\varepsilon\to\lambda\in[0,\infty)$, we derive the homogenized
	system on the whole domain. The subcritical obstacles leave no additional
	resistance term, and the cell problems are posed on the full periodic cell.
	The scalar correctors vanish, so the scalar diffusion operators remain
	unchanged, while the oscillating viscosity gives a time-dependent
	effective viscosity tensor. If $\lambda=0$, the limit decouples into an
	effective unsteady Stokes system and a Cahn--Hilliard system with source.
	If $\lambda>0$, the limit retains the Navier--Stokes--Cahn--Hilliard
	coupling, with convection, phase transport, and capillary forcing weighted
	by $\sqrt{\lambda}$. We also prove convergence of the time-integrated
	normalized microscopic energy to the corresponding macroscopic energy.
\end{abstract}

	
	
	{\small \tableofcontents}
	
\section{Introduction}

Multiphase transport in porous and microstructured media arises in groundwater
remediation, enhanced oil recovery, composite materials, hydrogels, and tissue
scaffolds. At the pore scale, the flow is governed by the interaction of viscous
transport, capillarity, interfacial energy, and phase separation. A corresponding
mathematical model must couple the hydrodynamics of an incompressible mixture
with the evolution of its diffuse interface.

Diffuse-interface models describe the fluid phases through a smooth order
parameter instead of a sharp interface. They accommodate changes of topology,
including coalescence and breakup, without tracking the fluid-fluid interface
explicitly. The order parameter evolves according to a Cahn--Hilliard equation
derived from a free-energy functional. Coupling this equation with the
incompressible Navier--Stokes equations gives the
Navier--Stokes--Cahn--Hilliard system, a standard model for binary
incompressible mixtures; see
\cite{anderson1998diffuse,feng2006fully,kim2012phase} and the references
therein.

The Navier--Stokes--Cahn--Hilliard system on fixed bounded domains has been
studied under several boundary conditions and constitutive assumptions
\cite{CH01,CH02,CH03,CH04,CH05,CH06}. A fluid-mechanical account of
diffuse-interface models is given in \cite{anderson1998diffuse}. Local existence
of weak solutions in three dimensions for unmatched viscosities was proved in
\cite{boyer1999mathematical}, while a frame-indifferent model for fluids with
different densities was introduced in \cite{abels_garcke_grun_2012}. For
matched densities, weak solutions and an associated energy law were obtained in
\cite{abels2009diffuse}. The case of degenerate mobility and different densities
was treated in \cite{abels_depner_garcke_2013}. Related models include systems
with moving contact lines \cite{CHB01}, linear source terms of Oono type
\cite{CHS01}, and Cahn--Hilliard--Brinkman systems with dynamic boundary
conditions \cite{twophase_CH_Brinkman}. For the stationary Navier--Stokes
equations with Navier boundary conditions, we refer to
\cite{tapia2021stokes}.

Phase-field models have also been used to describe two-phase flow in porous
media. A Darcy-scale model based on the Cahn--Hilliard equation was proposed in
\cite{ganesan2000_porous}. Several effective models have been derived by formal
asymptotic expansions. Schmuck \cite{schmuck2012_stokesch} derived an effective
Stokes--Cahn--Hilliard system for periodic microstructures at low Reynolds
numbers, and an asymptotic derivation of an effective
Navier--Stokes--Cahn--Hilliard system was given in
\cite{homog_two_phase_NSCH}. Further upscaled Cahn--Hilliard models were obtained
in \cite{schmuck2012upscaled}. Rigorous results are less numerous. Daly
et al.~\cite{daly2015_twofluid} derived Darcy-scale equations from a
pore-scale Stokes--Cahn--Hilliard model, while two-scale derivations of
macroscopic Stokes--Cahn--Hilliard systems were established in
\cite{Hari01,hari02}. Homogenization of Cahn--Hilliard-type equations by
evolutionary convergence methods was studied in
\cite{liero2018homogenization}. In a two-dimensional nonperforated setting, an
effective Navier--Stokes--Cahn--Hilliard system was obtained in
\cite{bunoiu2020homogenization}.

The limiting effect of periodically distributed obstacles depends on their size
relative to their separation. Periodic obstacles of size comparable to the
period lead to the classical porous-medium regime. For dilute inclusions, the
relevant three-dimensional capacity parameter is
\[
\rho_\e:=\frac{a_\e}{\e^3},
\]
where \(a_\e\) is the obstacle diameter and \(\e\) is the distance between
neighbouring obstacle centres. A nontrivial resistance may remain in the
capacity-critical regime, whereas the obstacle contribution disappears when
\(\rho_\e\to0\). For single-phase incompressible flow, the stationary and
evolutionary Stokes and Navier--Stokes equations in perforated domains have
been studied in
\cite{allaire1992_unsteady_stokes,Masmoudi02,feireisl_lu_homogenization_ns,
	lu2023homogenization}. Stokes equations with nonzero boundary data on the
obstacles were considered in \cite{hillairet2018homogenization}. Homogenization under slip conditions on the internal boundaries has also been
studied for single-phase flows. Allaire
\cite{slip1991allairehomogenization} considered the Stokes and
Navier--Stokes equations in domains perforated by periodically distributed
small obstacles with slip boundary conditions and identified the corresponding
Stokes, Brinkman, and Darcy regimes. More recently, Fabricius and Gahn
\cite{fabricius2023homogenization} studied the Stokes problem with Navier-slip
conditions in thin periodically perforated layers, combining homogenization
with dimension reduction. These results concern single-phase flow, whereas the
present work treats a time-dependent Navier--Stokes--Cahn--Hilliard system in a
three-dimensional domain with dilute subcritical obstacles and free-slip
conditions.

In this paper, the obstacles have diameter $a_\e=\e^\alpha$,  $\alpha>3$
and their centres are separated by distances of order \(\e\). Consequently,
\[
\rho_\e=\frac{a_\e}{\e^3}=\e^{\alpha-3}\longrightarrow0.
\]
The obstacles therefore lie in the subcritical dilute regime. Their total
volume also vanishes, but it is the preceding capacity scaling that excludes an
additional Darcy or Brinkman resistance in the limit.

Let \(S=(0,T)\), and let \(\O_p^\e\subset\R^3\) be the fluid part of a
periodically perforated domain. We study the system
\begin{equation}\label{MainP01}
	\left\{\begin{aligned}
		& \partial_t \mathbf u^\e + (\mathbf u^\e\!\cdot\nabla)\mathbf u^\e
		- \div\!\big(\mathbb A^\e\,D( \mathbf u^\e)\big)
		+ \nabla p^\e 
		+ \lambda^\e \phi^\e \nabla \mu^\e
		= \mathbf g^\e
		&& \text{in } S\times \Omega_p^\e,\\
		& \nabla\cdot \mathbf u^\e = 0 
		&& \text{in } S\times \Omega_p^\e,\\
		& \mathbf u^\e = \mathbf 0\quad \text{on }  S\times \partial\Omega,\qquad
		\mathbf u^\e\cdot\Gn^\e = 0,\quad 
		\big(\mathbb A^\e D(\mathbf u^\e)\mathbf n^\e\big)\!\cdot\boldsymbol\tau^\e = 0
		&& \text{on }  S\times \Gamma^\e_s,\\
		& \partial_t \phi^\e + \mathbf u^\e\!\cdot\nabla \phi^\e
		- \Delta \mu^\e + G(\phi^\e) = 0
		&& \text{in } S\times \Omega_p^\e,\\
		& \mu^\e = -\Delta \phi^\e + F^{'}(\phi^\e)
		&& \text{in } S\times \Omega_p^\e,\\
		& \mathbf n^\e\!\cdot\nabla \phi^\e = 0,\qquad 
		\mathbf n^\e\!\cdot\nabla \mu^\e = 0
		&& \text{on } S\times \partial\Omega_p^\e,\\
		& \mathbf u^\e(0,x) = \mathbf u^\e_0(x),
		\qquad
		\phi^\e(0,x) = \phi_0(x)
		&& \text{in } \Omega_p^\e,
	\end{aligned}\right.
\end{equation}
where $\O=\O_p^\e\cup\Gamma_s^\e\cup\O_s^\e$,
where \(\O_s^\e\) is the union of the solid obstacles and
\(\Gamma_s^\e=\partial\O_s^\e\) is the fluid-solid interface. The normal
\(\Gn^\e\) is directed outward from \(\O_p^\e\), and
\(\boldsymbol\tau^\e\) denotes a tangential direction on \(\Gamma_s^\e\).
The geometry is specified in Section~2.

The unknowns \(\Gu^\e\), \(p^\e\), \(\phi^\e\), and \(\mu^\e\) are the
velocity, hydrodynamic pressure, phase field, and chemical potential,
respectively. The tensor \(\A^\e\) describes a periodically heterogeneous
viscous response and acts on the symmetric gradient $D(\Gu^\e)
:=
\frac12\bigl(\nabla\Gu^\e+\nabla(\Gu^\e)^\top\bigr)$.
It is uniformly bounded and elliptic and is periodic in the fast variable; see
Subsection~\ref{SSec32}. The smooth double-well potential \(F\) determines the
interfacial free energy. The function \(G\) is a strictly monotone,
globally Lipschitz source term and makes the phase equation nonconservative.
The coefficient \(\lambda^\e>0\) controls the strength of the capillary force
\(\phi^\e\nabla\mu^\e\).

The outer boundary is subject to the no-slip condition. On each obstacle,
the velocity has zero normal component and zero tangential viscous traction.
These free-slip conditions allow a nonzero tangential velocity on
\(\Gamma_s^\e\). Consequently, the zero extension of \(\Gu^\e\) through the
obstacles is divergence-free but need not belong to \(H^1(\O)^3\). This
distinction between the zero extension and a primal \(H^1\)-extension is one
of the technical features of the compactness argument.

We assume
\begin{equation}\label{Lam01}
	\lambda^\e>0,
	\qquad
	\lambda^\e\longrightarrow\lambda\in[0,\infty).
\end{equation}
With bars denoting zero extensions to \(\O\), the initial velocities and body
forces are scaled so that
\begin{equation}\label{Assgg01}
	\frac{\overline{\Gu_0^\e}}{\sqrt{\lambda^\e}}
	\longrightarrow
	\Gu_0
	\quad\text{strongly in }L^2(\O)^3,\quad\text{and}\quad\frac{\overline{\Gg^\e}}{\sqrt{\lambda^\e}}
	\rightharpoonup
	\Gg
	\quad\text{weakly in }L^2(\O_T)^3.
\end{equation}
The energy estimates then show that \(\Gu^\e\) is of order
\(\sqrt{\lambda^\e}\). We therefore introduce the normalized velocity $\frac{\Gu^\e}{\sqrt{\lambda^\e}}$
and divide the momentum equation by \(\sqrt{\lambda^\e}\). Under this
normalization, the time derivative and viscous terms remain of order one,
whereas convection, phase transport, and capillary coupling carry the factor
\(\sqrt{\lambda^\e}\).

For every fixed \(\e>0\), we prove the existence of a weak solution of
\eqref{MainP01}. The estimates have constants independent of \(\e\) and retain
the explicit dependence on \(\lambda^\e\). The source term prevents conservation
of the spatial mean of \(\phi^\e\). Its strict monotonicity provides the
mean-value estimate required to control \(\phi^\e\) and \(\mu^\e\) uniformly.

The homogenization argument combines primal and zero extensions with periodic
unfolding. For the phase field, the zero extension carries the time derivative,
while the primal extension provides uniform spatial regularity. For the
velocity, a solenoidal primal extension gives an \(H^1\)-bounded field on the
fixed domain, whereas the divergence-free zero extension retains the uniform
\(L^\infty(S;L^2(\O)^3)\) estimate. Their difference converges strongly to zero
after normalization. Strong velocity compactness is then obtained from
time-compact scalar pairings against solenoidal restrictions of fixed-domain
test functions.

Passing to the limit in the momentum equation also requires oscillating
divergence-free test functions satisfying the impermeability condition on the
shrinking obstacles. These functions are constructed by cutting off a
two-scale solenoidal field near each obstacle and applying local Bogovski\u{i}
corrections in surrounding annuli. Since
\(a_\e/\e^3\to0\), the corrections vanish in the relevant norms and preserve the
required two-scale limits. The construction is given in appendix.

At the periodic scale, the obstacles collapse to a point, and the limiting cell
is the whole unit cell \(Y\). The scalar cell problems therefore have only
constant solutions. With the zero-mean normalization, $\widehat\phi=0=\widehat\mu$.
Thus the obstacles do not modify the scalar diffusion operators. The velocity
corrector remains nontrivial because \(\A^\e\) oscillates in the fast variable.
It determines a time-dependent effective viscosity tensor
\(\A^{\hom}(t)\). No additional obstacle resistance appears in the homogenized
momentum equation.

For completeness, we state the resulting equations in formal strong form. The
analysis itself is carried out in the divergence-free weak formulation. In
particular, the pressure appearing below is recovered on the fixed domain by a
de~Rham argument and is not identified as a limit of the microscopic pressures.

If \(\lambda=0\), the homogenized system is
\begin{equation}\label{HomStr01}
	\left\{
	\begin{aligned}
		\partial_t\Gu
		-\div\bigl(\A^{\hom}D(\Gu)\bigr)
		+\nabla p
		&=\Gg,
		\\
		\div\Gu&=0,
		\\
		\partial_t\phi-\Delta\mu+G(\phi)&=0,
		\\
		\mu&=-\Delta\phi+F'(\phi),
	\end{aligned}
	\right.
	\qquad
	\text{in }\O_T.
\end{equation}
It is supplemented by
\begin{equation}\label{HomBc}
	\begin{aligned}
		\Gu&=0
		\quad\text{on }S\times\partial\O,\quad
		\nabla\phi\cdot\Gn
		=\nabla\mu\cdot\Gn=0
		\text{on }S\times\partial\O,\\
		\Gu(0)&=\Gu_0,
		\qquad
		\phi(0)=\phi_0
		\quad\text{in }\O.
	\end{aligned}
\end{equation}
The velocity satisfies an effective unsteady Stokes system, while
\((\phi,\mu)\) satisfies a Cahn--Hilliard system with source \(G(\phi)\).
The two systems are decoupled.

If \(\lambda>0\), the limit retains the Navier--Stokes--Cahn--Hilliard
coupling:
\begin{equation}\label{LHomStr01}
	\left\{
	\begin{aligned}
		\partial_t\Gu
		+\sqrt{\lambda}\,(\Gu\cdot\nabla)\Gu
		-\div\bigl(\A^{\hom}D(\Gu)\bigr)
		+\nabla p
		+\sqrt{\lambda}\,\phi\nabla\mu
		&=\Gg,
		\\
		\div\Gu&=0,
		\\
		\partial_t\phi
		+\sqrt{\lambda}\,\Gu\cdot\nabla\phi
		-\Delta\mu+G(\phi)
		&=0,
		\\
		\mu&=-\Delta\phi+F'(\phi),
	\end{aligned}
	\right.
	\qquad
	\text{in }\O_T,
\end{equation}
with the boundary and initial conditions \eqref{HomBc}. The factor
\(\sqrt{\lambda}\) determines the strength of convection, phase transport, and
capillary coupling. In both regimes, only the effective viscosity is
nontrivial; the scalar diffusion coefficients remain unchanged.

We also prove convergence of the time-integrated scaled microscopic energy:
\[
\lim_{\e\to0}
\int_S
\frac{\GT^\e(\Gu^\e(t),\phi^\e(t))}{\lambda^\e}
\,\md t
=
\int_S
\GT(\Gu(t),\phi(t))
\,\md t.
\]
This result uses strong convergence of the normalized velocity and a
norm-identification argument for the phase gradients.

Section~\ref{Sec03} describes the perforated geometry, assumptions, and weak formulation.
Section~\ref{Sec04} collects the uniform extension, divergence, Korn, Poincar\'e, and
Sobolev estimates used throughout the paper. Section~\ref{Sec05} establishes the
a priori estimates and existence of microscopic weak solutions. Section~\ref{SSec61}
introduces periodic unfolding and proves the compactness results. Section~\ref{Sec06}
derives the two-scale limit system, identifies the cell problems, and obtains
the homogenized equations and energy convergence. The technical constructions
of the oscillating solenoidal test functions and solenoidal restrictions are
given in the Appendix~\ref{AppSec}.
	
\section*{Notation and function spaces}

We introduce the functional setup and general notation on $X\subset\R^3$, a bounded
domain with Lipschitz boundary. We denote by $L^r(X)$ and $W^{k,r}(X)$ the usual
Lebesgue and Sobolev spaces with their standard norms. In particular,
$H^k(X)=W^{k,2}(X)$, and $H^1_0(X)$ denotes the subspace of $H^1(X)$ with zero
trace on $\partial X$. The mean-zero subspace of $L^2(X)$ is denoted by $L^2_0(X)$.
We write $H^1_{0,\div}(X)^3$ for the divergence-free subspace of $H^1_0(X)^3$, with
dual $H^{-1}_{0,\div}(X)^3:=(H^1_{0,\div}(X)^3)^\ast$; similarly $H^{-1}(X)$
denotes the dual of $H^1(X)$. We denote by $\langle\cdot,\cdot\rangle$ and
$(\cdot,\cdot)$ the duality pairing $\langle\cdot,\cdot\rangle_{H^{-1},H^1}$ and
the inner product $(\cdot,\cdot)_{L^2}$, respectively.

The space of $m$-times continuously differentiable functions on $X$ is $C^m(X)$.
We denote by $C^\infty_c(X)$ the space of smooth compactly supported functions in
$X$, and by $C^\infty_{c,\div}(X)^3$ its divergence-free subspace. On the
reference cell $Y$, $C^\infty_{\per}(Y)$ denotes the space of $Y$-periodic smooth
functions, with divergence-free subspace $C^\infty_{\#,\div}(Y)^3$, and
$C^\infty_\#(Y)$ denotes the space of $Y$-periodic smooth functions with zero mean
over $Y$. Similarly, $H^1_\per(Y)$ denotes the $Y$-periodic $H^1$-functions on $Y$,
and $H^1_\#(Y)$ the $Y$-periodic $H^1$-functions on $Y$ with zero mean.

Moreover, we set
\begin{equation}\label{eq:set}
\begin{aligned}
	H^\sigma
	&:=
	\overline{C^\infty_{c,\div}(\O)^3}^{\,L^2(\O)^3},\quad
	H_\e^\sigma
	:=
	\overline{V_\e^\sigma}^{\,L^2(\O_p^\e)^3},\\ 		V_\e^\sigma&:=\bigl\{\bv\in V_\e:\div ~\bv=0\text{ in }\O_p^\e,\ \bv=0\text{ on }\partial\O,~\bv\cdot\Gn^\e=0\text{ on }\Gamma_s^\e\bigr\}.
\end{aligned}
\end{equation}
In estimates we write $L^2(X)$ for $L^2(X)^3$ or $L^2(X)^{3\times3}$, indicating
the full product space only when stating weak or strong convergence. We use the
Einstein summation convention. The symbol $C$ denotes a generic positive constant
independent of $\e$, whose value may change from line to line.
	
	\section{Problem setting}\label{Sec03}

\subsection{Domain description}\label{SSec31}

Let $\O\subset\mathbb{R}^3$ be a bounded domain with Lipschitz boundary, and let
$Y=(0,1)^3$ be the unit reference cell. Fix a non-empty open set
$Y_s\Subset Y$ with Lipschitz boundary such that $Y_p:=Y\setminus\overline{Y_s}$
is open and connected, and fix a point $y_0\in Y_s$. The centered reference
obstacle is $T:=Y_s-y_0$, so that $0\in T$ and $y_0+T=Y_s$. Since $Y_s\Subset Y$,
we have $\operatorname{dist}(\overline{Y_s},\partial Y)\ge d_0$ for some $d_0>0$.

Throughout, we fix $\alpha>3$ and set
\begin{equation}\label{EqHoleSize}
	a_\e:=\e^\alpha,
	\qquad
	\delta_\e:=\frac{a_\e}{\e}=\e^{\alpha-1},
\end{equation}
so that $\delta_\e\to0$ as $\e\to0$. For $0<\e<1$, let
\[
\Kc_\e:=\bigl\{\kappa\in\mathbb Z^3:\e(\kappa+\overline Y)\subset\O\bigr\}
\]
index the complete periodic cells, with union
\[
\widehat\O^\e:=\operatorname{int}\Bigl(\,\bigcup_{\kappa\in\Kc_\e}\e(\kappa+\overline Y)\Bigr),
\qquad
\Lambda^\e:=\O\setminus\overline{\widehat\O^\e}.
\]
Since $\partial\O$ is Lipschitz, there is a constant $c>0$, independent of $\e$,
with $\Lambda^\e\subset\{x\in\O:\operatorname{dist}(x,\partial\O)<c\e\}$; hence
$|\Lambda^\e|\to0$ as $\e\to0$.

For $\kappa\in\Kc_\e$, the obstacle in the cell $\e(\kappa+Y)$ is centered at
$x_{\e,\kappa}:=\e(\kappa+y_0)$ and defined by
\begin{equation}\label{EqSingleHole}
	T_{\e,\kappa}
	:=x_{\e,\kappa}+a_\e T
	=\e\kappa+\e y_0+\e^\alpha(Y_s-y_0).
\end{equation}
Thus each obstacle has diameter of order $\e^\alpha$, while neighbouring centres
are at distance of order $\e$. At the reference-cell scale we set
\[
Y_{s,\e}:=y_0+\delta_\e T,
\qquad
Y_{p,\e}:=Y\setminus\overline{Y_{s,\e}},
\qquad\text{so that}\qquad
T_{\e,\kappa}=\e(\kappa+Y_{s,\e}).
\]
Since $\delta_\e\to0$, the rescaled obstacle collapses to the point $y_0$, and
\begin{equation}\label{EqSolidCellVol}
	|Y_{s,\e}|=\delta_\e^3|T|=\e^{3(\alpha-1)}|T|\longrightarrow0.
\end{equation}
The limiting periodic cell is therefore the whole unit cell $Y$, not the fixed
perforated cell $Y_p$.

The solid and fluid domains are
\[
\O_s^\e:=\bigcup_{\kappa\in\Kc_\e}T_{\e,\kappa},
\qquad
\O_p^\e:=\O\setminus\overline{\O_s^\e},
\qquad
\widehat\O_p^\e:=\widehat\O^\e\setminus\overline{\O_s^\e}.
\]
No obstacle lies in the boundary layer, so $\O_p^\e=\widehat\O_p^\e\cup\Lambda^\e$
up to their common interface. The union of the solid--fluid interfaces is
\[
\Gamma_s^\e:=\partial\O_s^\e
=\bigcup_{\kappa\in\Kc_\e}\partial T_{\e,\kappa}
=\bigcup_{\kappa\in\Kc_\e}\bigl(x_{\e,\kappa}+\e^\alpha\partial T\bigr),
\]
whence $\partial\O_p^\e=\partial\O\cup\Gamma_s^\e$ and
$\overline\O=\overline{\O_p^\e}\cup\overline{\O_s^\e}$.

Since $T$ is bounded and $y_0\in Y$, there is $\e_0>0$ such that
$\overline{Y_{s,\e}}\subset Y$, equivalently
$\overline{T_{\e,\kappa}}\subset\e(\kappa+Y)$ for all $\kappa\in\Kc_\e$, whenever
$0<\e<\e_0$. Shrinking $\e_0$ if necessary, there is a constant $c_0>0$,
independent of $\e$, with
\[
\operatorname{dist}\bigl(T_{\e,\kappa},\partial(\e(\kappa+Y))\bigr)\ge c_0\e
\qquad\text{for all }\kappa\in\Kc_\e,
\]
so the obstacles are pairwise disjoint and uniformly separated at the periodicity
scale. As $\#\Kc_\e=\GO(\e^{-3})$,
\begin{equation}\label{EqSolidVol}
	|\O_s^\e|=\#\Kc_\e\,\e^{3\alpha}|T|=\GO\bigl(\e^{3(\alpha-1)}\bigr)\longrightarrow0
	\qquad\text{as }\e\to0.
\end{equation}
Finally, $\rho_\e=\frac{a_\e}{\e^3}=\e^{\alpha-3}\longrightarrow0$,
which measures the cumulative influence of the obstacles. Since $\alpha>3$, this
is the subcritical regime.
\subsection{Model assumptions}\label{SSec32}

Throughout the paper the exponent $\alpha>3$ is fixed. Unless stated otherwise,
$C>0$ denotes a generic constant independent of $\e$, whose value may change from
line to line. We collect here the assumptions used in the existence and
homogenization analysis.

\textbf{Viscosity.}
Let $\mathbb A\in
L^\infty\bigl(S;L^\infty_{\rm per}(Y)^{3\times3\times3\times3}\bigr)$
be a $Y$-periodic fourth-order tensor, and set
\[
\mathbb A^\e(t,x):=\mathbb A\left(t,\frac{x}{\e}\right),
\qquad (t,x)\in S\times\O.
\]
For $\Xi\in\mathbb R_{\rm sym}^{3\times3}$, we define
\[
\bigl(\mathbb A(t,y)\Xi\bigr)_{ij}
:=\sum_{k,l=1}^3\mathbb A_{ijkl}(t,y)\Xi_{kl}.
\]
We assume the minor and major symmetries
\[
\mathbb A_{ijkl}
=\mathbb A_{jikl}
=\mathbb A_{ijlk}
=\mathbb A_{klij}
\]
and that there exist $\kappa_0,\kappa_1>0$ such that, for a.e.\
$(t,y)\in S\times Y$,
\begin{equation}\label{Coe01}
	\|\mathbb A\|_{L^\infty(S\times Y)^{3\times3\times3\times3}}
	\leq\kappa_0,
	\qquad
	\mathbb A(t,y)\Xi:\Xi\geq\kappa_1|\Xi|^2
	\quad
	\text{for all }\Xi\in\mathbb R_{\rm sym}^{3\times3}.
\end{equation}
Thus $\mathbb A^\e$ has the same symmetries and satisfies \eqref{Coe01}
uniformly in $\e$. Its restriction to $S\times\O_p^\e$ is denoted by the
same symbol.

For the periodic unfolding operator $\Te$ defined in
Section~\ref{SSec61},
\[
\Te(\mathbb A^\e)(t,x,y)=\mathbb A(t,y)
\quad\text{for a.e. }(t,x,y)\in S\times\widehat\O^\e\times Y,
\]
whereas $\Te(\mathbb A^\e)=0$ on $S\times\Lambda^\e\times Y$. Since
$|\Lambda^\e|\to0$, we have
\begin{equation}\label{PeriodicCoefficientConvergence}
	\Te(\mathbb A^\e)\longrightarrow\mathbb A
	\quad\text{strongly in }
	L^q(\O_T\times Y)^{3\times3\times3\times3}
\end{equation}
for every $1\leq q<\infty$.

\textbf{Capillary strength.}
The sequence $\{\lambda^\e\}_{\e>0}$ consists of positive reals satisfying
\eqref{Lam01}.

\textbf{Initial data.}
Let $\phi_0\in H^1(\O)$, we set $\phi_0^\e:=\phi_0|_{\O_p^\e}$. Moreover, let $\Gu_0^\e\in H_\e^\sigma$ (see \eqref{eq:set}), and $\overline\Gu_0^\e$ denote the zero extension of
$\Gu_0^\e$ to $\O$. We assume that there exist $\Gu_0\in H^\sigma$ such that
\begin{equation}\label{ICA01}
	\frac{\overline\Gu_0^\e}{\sqrt{\lambda^\e}}
	\longrightarrow
	\Gu_0
	\quad\text{strongly in }L^2(\O)^3,
	\qquad
	\|\Gu_0^\e\|_{L^2(\O_p^\e)}
	\leq
	\kappa_4\sqrt{\lambda^\e},
	\qquad
	\|\phi_0\|_{H^1(\O)}
	\leq
	\kappa_3,
\end{equation}
for constants $\kappa_3,\kappa_4>0$ independent of $\e$.

\textbf{External force.}
Let $\Gg\in L^2(S\times\O)^3$. The microscopic body force is
\begin{equation}\label{Ass01}
	\Gg^\e=\sqrt{\lambda^\e}\,\Gg\quad\text{in }S\times\O_p^{\e},
	\qquad
	\|\Gg\|_{L^2(S\times\O)}\leq\kappa_5,
\end{equation}
for a constant $\kappa_5>0$ independent of $\e$.

\textbf{Double-well potential.}
The smooth double-well potential $F:\mathbb{R}\to\mathbb{R}$ is
\begin{equation}\label{PotentialF}
	F(s):=\tfrac{1}{4}(s^2-1)^2,\qquad s\in\mathbb{R},
\end{equation}
with derivative $F'(s)=s^3-s$. It satisfies
\begin{equation}\label{F01}
	\tfrac{1}{2}s^2-\tfrac{3}{4}\leq F(s)\leq\tfrac{1}{4}s^4+\tfrac{1}{4},
	\qquad
	|F'(s)|\leq|s|^3+|s|\leq C(1+|s|^3),
	\qquad s\in\mathbb{R}.
\end{equation}

\textbf{Source term.}
The source term $G:\mathbb{R}\to\mathbb{R}$ satisfies
\begin{equation}\label{AssG01}
	G\in C^1(\mathbb{R}),\qquad G(0)=0,\qquad
	0<c_1\leq G'(s)\leq c_2\quad\text{for all }s\in\mathbb{R},
\end{equation}
for constants $c_1,c_2>0$. Consequently,
\begin{equation}\label{SourceGrowth}
	c_1|s|\leq|G(s)|\leq c_2|s|,
	\qquad
	c_1s^2\leq G(s)s\leq c_2s^2,
	\qquad s\in\mathbb{R}.
\end{equation}
%

\begin{remark}
	The particular choice \eqref{Ass01} may be replaced by the
	assumption
	\[
	\frac{\overline{\Gg^\e}}{\sqrt{\lambda^\e}}
	\rightharpoonup
	\Gg
	\quad\text{weakly in }L^2(\O_T)^3,
	\]
	where $\overline{\Gg^\e}$ denotes the zero extension of $\Gg^\e$ from
	$\O_p^\e$ to $\O$. This convergence implies the uniform estimate $\|\Gg^\e\|_{L^2(S\times\O_p^\e)}
	\leq C\sqrt{\lambda^\e}$.
	The arguments below remain unchanged under this assumption.
\end{remark}

\subsection{Weak formulation}\label{SSec33}

We first introduce the velocity space adapted to the mixed no-slip/free-slip
boundary conditions:
\[
\begin{aligned}
	\GH^1(\O_p^\e)&:=\bigl\{\mathbf v\in H^1(\O_p^\e)^3:
	\mathbf v=0\ \text{on }\partial\O,\ \mathbf v\cdot\mathbf n^\e=0\ \text{on }\Gamma_s^\e\bigr\},\\
	\GH^1_{\div}(\O_p^\e)&:=\bigl\{\mathbf v\in\GH^1(\O_p^\e):
	\nabla\!\cdot\mathbf v=0\ \text{in }\O_p^\e\bigr\},
\end{aligned}
\]
with dual $\GH^{-1}_{\div}(\O_p^\e)$. The admissible spaces are
\[
\begin{aligned}
	\mathfrak{U}^\e&=L^\infty(S;L^2(\O_p^\e))^3\cap L^2(S;\GH^1_{\div}(\O_p^\e))\cap W^{1,4/3}(S;\GH^{-1}_{\div}(\O_p^\e)),\\
	\mathfrak{C}^\e&=L^\infty(S;H^1(\O_p^\e))\cap H^1(S;H^{-1}(\O_p^\e)),\qquad
	\mathfrak{W}^\e=L^2(S;H^1(\O_p^\e)),
\end{aligned}
\]
with macroscopic counterparts
\[
\begin{aligned}
	\mathfrak{U}&=L^\infty(S;L^2(\O))^3\cap L^2(S;H^1_{0,\div}(\O)^3),\\
	\mathfrak{C}&=L^\infty(S;H^1(\O))
	\cap
	C([0,T];L^2(\O)),\qquad
	\mathfrak{W}=L^2(S;H^1(\O)).
\end{aligned}
\]

\begin{definition}[Weak solution]\label{Def01}
	A triplet
	$(\mathbf u^\e,\phi^\e,\mu^\e)\in\mathfrak{U}^\e\times\mathfrak{C}^\e\times\mathfrak{W}^\e$
	is a weak solution of \eqref{MainP01} if $\mathbf u^\e(0)=\mathbf u_0^\e$ and
	$\phi^\e(0)=\phi_0$ in $\O_p^\e$, and the following identities hold for all
	$\varphi_1\in L^4(S;\GH^1_{\div}(\O_p^\e))$ and
	$\varphi_2,\varphi_3\in L^2(S;H^1(\O_p^\e))$:
	\begin{align}
		\label{weak-momentum1}
		&\int_S\langle\partial_t\mathbf{u}^\e,\varphi_1\rangle\,\md t
		+\int_{S\times\O_p^\e}(\mathbf u^\e\!\cdot\nabla)\mathbf u^\e\cdot\varphi_1\,\md(x,t)
		+\int_{S\times\O_p^\e}\A^\e D(\mathbf{u}^\e):D(\varphi_1)\,\md(x,t)\nonumber\\
		&\qquad\qquad
		+\int_{S\times\O_p^\e}\lambda^\e\phi^\e\nabla\mu^\e\cdot\varphi_1\,\md(x,t)
		=\int_{S\times\O_p^\e}\Gg^\e\cdot\varphi_1\,\md(x,t),\\[1ex]
		\label{weak-phi1}
		&\int_S\langle\partial_t\phi^\e,\varphi_2\rangle\,\md t
		+\int_{S\times\O_p^\e}G(\phi^\e)\varphi_2\,\md(x,t)
		+\int_{S\times\O_p^\e}\nabla\mu^\e\cdot\nabla\varphi_2\,\md(x,t)
		+\int_{S\times\O_p^\e}(\mathbf{u}^\e\!\cdot\nabla\phi^\e)\varphi_2\,\md(x,t)=0,\\[1ex]
		\label{weak-mu1}
		&\int_{S\times\O_p^\e}\mu^\e\varphi_3\,\md(x,t)
		=\int_{S\times\O_p^\e}\nabla\phi^\e\cdot\nabla\varphi_3\,\md(x,t)
		+\int_{S\times\O_p^\e}F'(\phi^\e)\varphi_3\,\md(x,t).
	\end{align}
\end{definition}

	\section{Preliminary results}\label{Sec04}

\subsection{Bogovski\u{\i} operators}\label{SSec41}
We introduce right inverses of the divergence, the Bogovski\u{\i} operators on the solid obstacle domain $\O_s^\e$;
see \cite{BogovM01,novotny2004,galdi2011,DieningFeireislLu17,Allaire89,Masmoudi02,LuUniformStokes21,NecasovaPan22}. On a fixed bounded Lipschitz domain the
existence of such a right inverse is classical; see
\cite{AcostaDuranMuschietti06,CostabelMcIntosh10} and the references therein.

Recall that
\[
T_{\e,\kappa}=x_{\e,\kappa}+a_\e T,
\qquad
x_{\e,\kappa}=\e(\kappa+y_0),
\qquad
a_\e=\e^\alpha,\quad\alpha>3.
\]

\begin{lemma}\label{Bogov01+}
	Let
	\[
	L^2_{0,\mathrm{comp}}(\O_s^\e)
	:=\Bigl\{g\in L^2(\O_s^\e):
	\int_{T_{\e,\kappa}}g\,\md x=0\ \text{for every }\kappa\in\Kc_\e\Bigr\}
	\]
	be the space of functions with zero mean on each connected component of
	$\O_s^\e$. Then $L^2_{0,\mathrm{comp}}(\O_s^\e)\subset L^2_0(\O_s^\e)$, and
	there exists a linear operator
	\[
	\mathcal B_s^\e:L^2_{0,\mathrm{comp}}(\O_s^\e)\longrightarrow H^1_0(\O_s^\e)^3
	\]
	such that, for every $g\in L^2_{0,\mathrm{comp}}(\O_s^\e)$, the field
	$\phi:=\mathcal B_s^\e(g)$ satisfies
	\[
	\div\phi=g\quad\text{in }\O_s^\e,
	\qquad
	\phi=0\quad\text{on }\partial\O_s^\e
	\]
	in the sense of distributions and traces, respectively. Moreover, there exists
	a constant $C_s>0$, independent of $\e$, $\kappa$, and $g$, such that
	\begin{equation}\label{Bogov02+-gradient}
		\|\nabla\mathcal B_s^\e(g)\|_{L^2(\O_s^\e)}\leq C_s\|g\|_{L^2(\O_s^\e)},
		\qquad
		\|\mathcal B_s^\e(g)\|_{L^2(\O_s^\e)}\leq C_s a_\e\|g\|_{L^2(\O_s^\e)}.
	\end{equation}
\end{lemma}
The proof is given in Appendix \ref{App01}.

\subsection{Further results in the perforated domain}\label{SSec42}

We construct, for each field, a primal extension and a zero extension. The
primal extension is uniformly bounded in the natural energy space and carries
the weak limits together with, under a time-derivative bound, the strong limits.
The zero extension carries the distributional time derivative and, for
solenoidal fields, the divergence-free constraint. The two extensions differ in
$L^2$ by a term of order $a_\e/\e$, so they share the same $L^2$ limit whenever
either one converges.

\begin{lemma}[Extension on the primal level]\label{Emain}
	There exist $\e_0\in(0,1)$ and $C>0$, independent of
	$\e\in(0,\e_0)$, such that, for every $\e\in(0,\e_0)$, there is a
	bounded linear operator $E^\e:H^1(\O_p^\e)\longrightarrow H^1(\O)$
	such that
	$$
	E^\e u=u
	\quad\text{a.e. in }\O_p^\e.
	$$
	Moreover,
	\begin{equation}\label{eq:extension-primal-estimates}
		\|E^\e u\|_{L^2(\O)}
		\leq
		C\|u\|_{L^2(\O_p^\e)},
		\qquad
		\|\nabla E^\e u\|_{L^2(\O)}
		\leq
		C\|\nabla u\|_{L^2(\O_p^\e)}
	\end{equation}
	for every $u\in H^1(\O_p^\e)$.
	Furthermore, if $u=0$ on $\partial\O$ in the sense of traces, then  $E^\e u\in H^1_0(\O)$.
\end{lemma}
	The result follows from the standard extension theorem for periodically
perforated domains with isolated Lipschitz holes; see
\cite{Hopker16,oleinik01,AdamsFournier03}.

\begin{lemma}[Primal and zero extensions of scalar fields]
	\label{lem:scalar-two-extensions}
	Let $v^\e\in \frak{C}^\e$
	satisfy
	\[
	\|v^\e\|_{L^\infty(S;H^1(\O_p^\e))}
	+
	\|\partial_tv^\e\|_
	{L^2(S;H^{-1}(\O_p^\e))}
	\leq C.
	\]
	Define the primal extension $\widetilde v^\e:=E^\e v^\e$
	and the zero extension $\overline v^\e$ to $S\X\O$.
	Then
	\begin{equation}
		\label{ScalarPrimalBound}
		\|\widetilde v^\e\|_
		{L^\infty(S;H^1(\O))}
		\leq C,\quad 	
		\overline v^\e
		\in
		H^1(S;H^{-1}(\O)),
		\qquad
		\|\partial_t\overline v^\e\|_
		{L^2(S;H^{-1}(\O))}
		\leq
		C.
	\end{equation}
	Moreover,
	\begin{equation}
		\label{ScalarExtensionsDifference}
		\|\widetilde v^\e-\overline v^\e\|_
		{L^\infty(S;L^2(\O))}
		\leq
		C\frac{a_\e}{\e}.
	\end{equation}
	Consequently, there exist $v\in\frak{C}$
	and a subsequence such that
	\begin{equation}
		\label{ScalarTwoExtensionStrongCompactness}
		\widetilde v^\e
		\longrightarrow v,
		\qquad
		\overline v^\e
		\longrightarrow v
		\quad\text{strongly in }
		C([0,T];L^2(\O)).
	\end{equation}
\end{lemma}
The proof of the above lemma is presented in Appendix \ref{App01}.

\begin{lemma}[Primal and zero extensions of solenoidal fields]\label{LemD02}
	Let
	\[
	\begin{aligned}
		V_\e&:=H^1(\O_p^\e)^3,\quad  V:=H^1(\O)^3,\\
		V^\sigma&:=\bigl\{w\in V:\div w=0\text{ in }\O,\ w=0\text{ on }\partial\O\bigr\}.
	\end{aligned}
	\]
	There exist $\e_0\in(0,1)$ and $C>0$, independent of $\e\in(0,\e_0)$, such that
	for every $\e\in(0,\e_0)$ there is a linear operator
	$\mathcal E^\e:V_\e^\sigma\to V^\sigma$ with
	\begin{equation}\label{SolenoidalPrimalExtensionIdentity}
		(\mathcal E^\e v)|_{\O_p^\e}=v\qquad\forall\,v\in V_\e^\sigma,
	\end{equation}
	\begin{equation}\label{SolenoidalPrimalExtensionBounds}
		\|\mathcal E^\e v\|_{L^2(\O)}\leq C\bigl(\|v\|_{L^2(\O_p^\e)}+a_\e\|\nabla v\|_{L^2(\O_p^\e)}\bigr),
		\qquad
		\|\nabla\mathcal E^\e v\|_{L^2(\O)}\leq C\|\nabla v\|_{L^2(\O_p^\e)},
	\end{equation}
	and, in particular,
	\begin{equation}\label{SolenoidalPrimalExtensionH1Bound}
		\|\mathcal E^\e v\|_{H^1(\O)}\leq C\|v\|_{H^1(\O_p^\e)}.
	\end{equation}
	For $v\in V_\e^\sigma$, let $\overline v^\e$ be its zero extension through the
	obstacles, $\overline v^\e:=v$ in $\O_p^\e$ and $\overline v^\e:=0$ in
	$\O_s^\e$. Then
	\begin{equation}\label{SolenoidalZeroExtension}
		\|\overline v^\e\|_{L^2(\O)}=\|v\|_{L^2(\O_p^\e)},
		\qquad
		\div\overline v^\e=0\ \text{in }\mathcal D'(\O).
	\end{equation}
	Moreover, for every $v^\e\in L^2(S;V_\e^\sigma)$,
	\begin{equation}\label{PrimalZeroExtensionDifference}
		\|\mathcal E^\e v^\e-\overline v^\e\|_{L^2(S\times\O)}
		\leq C\,\frac{a_\e}{\e}\,\|v^\e\|_{L^2(S;H^1(\O_p^\e))}.
	\end{equation}
	Since $a_\e/\e=\e^{\alpha-1}\to0$, if
	$\sup_{\e\in(0,\e_0)}\|v^\e\|_{L^2(S;H^1(\O_p^\e))}<\infty$, then
	$\mathcal E^\e v^\e-\overline v^\e\to0$ strongly in $L^2(S\times\O)^3$; in
	particular the two extensions share the same weak or strong
	$L^2(S\times\O)^3$ limit whenever either exists.
\end{lemma}
The proof of the above lemma is given in Appendix \ref{App01}.

\begin{remark}
	Under the free-slip condition, the zero extension
	$\overline v^\e$ is divergence-free in $\mathcal D'(\O)$, but it
	need not belong to $H^1(\O)^3$ because the tangential trace of $v$
	may be nonzero on $\Gamma_s^\e$. The primal extension
	$\mathcal E^\e v$ provides global $H^1$ regularity, while the
	Bogovski\u\i\ correction restores the divergence-free constraint
	inside the obstacles.
	
	If the free-slip condition is replaced by the no-slip condition
	$v=0$ on $\Gamma_s^\e$, the zero extension belongs to
	$H_0^1(\O)^3$ and is divergence-free. It can then be used as the
	solenoidal $H^1$ extension.
\end{remark}
\begin{lemma}[Poincar\'e inequality on the free-slip space]\label{lem:Poincare-GH1}
	There exists a constant $C>0$, independent of $\e$, such that
	\begin{equation}\label{eq:Poincare-GH1}
		\|\bv\|_{L^2(\O_p^\e)}\leq C\|\nabla\bv\|_{L^2(\O_p^\e)}
		\qquad\forall\,\bv\in\GH^1(\O_p^\e).
	\end{equation}
\end{lemma}

\begin{proof}
	Let $\bv\in\GH^1(\O_p^\e)$ and let $E^\e$ be the extension operator of
	Lemma~\ref{Emain}. Since $\bv=0$ on $\partial\O$ in the sense of traces,
	$E^\e\bv\in H^1_0(\O)^3$ and $E^\e\bv=\bv$ a.e.\ in $\O_p^\e$. By the classical
	Poincar\'e inequality on $\O$ and the uniform estimates
	\eqref{eq:extension-primal-estimates},
	\[
	\|\bv\|_{L^2(\O_p^\e)}\leq\|E^\e\bv\|_{L^2(\O)}
	\leq C\|\nabla E^\e\bv\|_{L^2(\O)}\leq C\|\nabla\bv\|_{L^2(\O_p^\e)}.
	\]
	The constant depends only on $\O$ and on the extension bound, both independent
	of $\e$.
\end{proof}

\begin{lemma}[Korn's inequality]
	\label{Korn}
	There exists a constant $C>0$, independent of $\e$, such that
	\begin{equation}\label{eq:Korn-uniform}
		\|\Gu\|_{L^2(\O_p^\e)}
		+
		\|\nabla\Gu\|_{L^2(\O_p^\e)}
		\leq
		C\|D(\Gu)\|_{L^2(\O_p^\e)}
	\end{equation}
	for every $\Gu\in\GH^1(\O_p^\e)$.
\end{lemma}

\begin{proof}
	By the uniform elasticity extension theorem for periodically
	perforated domains, there exists a linear operator $\mathscr E^\e:
	H^1(\O_p^\e)^3
	\longrightarrow
	H^1(\O)^3$
	such that
	\[
	\mathscr E^\e\Gu=\Gu
	\quad\text{a.e. in }\O_p^\e,
	\qquad
	\|D(\mathscr E^\e\Gu)\|_{L^2(\O)}
	\leq
	C\|D(\Gu)\|_{L^2(\O_p^\e)},
	\]
	with $C$ independent of $\e$; see
	\cite{Ext01,oleinik01}. Since the extension preserves the outer
	zero trace, $\mathscr E^\e\Gu\in H^1_0(\O)^3$. Korn's inequality
	on the fixed domain $\O$ therefore gives
	\[
	\|\mathscr E^\e\Gu\|_{H^1(\O)}
	\leq
	C\|D(\mathscr E^\e\Gu)\|_{L^2(\O)}
	\leq
	C\|D(\Gu)\|_{L^2(\O_p^\e)}.
	\]
	Restricting the left-hand side to $\O_p^\e$ proves
	\eqref{eq:Korn-uniform}.
\end{proof}

We next recall the Poincar\'e--Wirtinger inequality on $\O^\e_p$; see
\cite{slip1991allairehomogenization,LuUniformStokes21,oleinik01,Ext01,Ext02} for
related inequalities on perforated domains.

\begin{lemma}[Poincar\'e--Wirtinger on periodic perforations]\label{LemPW}
	For every $u\in H^1(\O^\e_p)$ with zero mean on $\O^\e_p$, i.e.\
	$u\in L^2_0(\O^\e_p)$,
	\[
	\|u\|_{L^2(\O^\e_p)}\leq C\|\nabla u\|_{L^2(\O^\e_p)},
	\]
	with $C>0$ independent of $\e$.
\end{lemma}

\begin{lemma}\label{Embed01}
	For every $u\in H^1(\O^\e_p)$ and $r\in[2,6]$, $u\in L^r(\O^\e_p)$ with
	\[
	\|u\|_{L^r(\O^\e_p)}\leq C\|u\|_{H^1(\O^\e_p)},
	\]
	where $C$ is independent of $\e$ and $u$.
\end{lemma}

The proof follows from the uniform extension of Lemma~\ref{Emain} and the standard
Sobolev embedding on $\O$. Thus $H^1(\O^\e_p)$ embeds continuously into
$L^r(\O^\e_p)$ for $r\in[2,6]$, with embedding constant independent of $\e$.

\subsection{Convective trilinear form}\label{SSec43}

Consider the convective trilinear form
\[
b^\e(\Gu,\Gv,\Gw):=\int_{\O^\e_p}(\Gu\cdot\nabla)\Gv\cdot\Gw\,\md x,
\qquad\Gu,\Gv,\Gw\in L^2(S;\GH^1_{\div}(\O^\e_p)).
\]

\begin{lemma}\label{LemTrilinear}
	For all $\Gu,\Gv,\Gw\in L^2(S;\GH^1_{\div}(\O^\e_p))$, the form $b^\e$ satisfies:
	\begin{itemize}
		\item[\textit{(i)}] \textit{Skew-symmetry:}
		$b^\e(\Gu,\Gv,\Gw)=-b^\e(\Gu,\Gw,\Gv)$;
		\item[\textit{(ii)}] \textit{Orthogonality:}
		\begin{equation}\label{skew01}
			b^\e(\Gu,\Gv,\Gv)=\int_{\O^\e_p}(\Gu\cdot\nabla)\Gv\cdot\Gv\,\md x=0;
		\end{equation}
		\item[\textit{(iii)}] \textit{Identity:}
		\begin{equation}\label{Id01}
			b^\e(\Gu,\Gu,\Gv)=-\int_{\O^\e_p}(\Gu\otimes\Gu):\nabla\Gv\,\md x,
		\end{equation}
		where $\otimes$ is the outer product and $:$ the Frobenius product;
		\item[\textit{(iv)}] \textit{Continuity:}
		\begin{equation}\label{117}
			|b^\e(\Gu,\Gv,\Gw)|\leq C_0\|\Gu\|_{H^1(\O^\e_p)}\|\Gv\|_{H^1(\O^\e_p)}\|\Gw\|_{H^1(\O^\e_p)},
		\end{equation}
		with $C_0$ independent of $\e$.
	\end{itemize}
\end{lemma}

\begin{proof}
	Skew-symmetry follows by integration by parts and the divergence theorem, using
	$\int_{\partial\O^\e_p}(\Gu\cdot\Gn^\e)(\Gv\cdot\Gw)\,\md S=0$, a consequence of
	the boundary conditions on $\Gu$. Setting $\Gv=\Gw$ gives
	$b^\e(\Gu,\Gv,\Gv)=-b^\e(\Gu,\Gv,\Gv)$, hence \eqref{skew01}.
	
	Since $\nabla\cdot\Gu=0$ in $S\times\O^\e_p$ and $\Gu=0$ on
	$S\times\partial\O$,
	\[
	b^\e(\Gu,\Gu,\Gv)
	=-\int_{\O^\e_p}(\Gu\otimes\Gu):\nabla\Gv\,\md x
	+\int_{\Gamma^\e_s}(\Gu\cdot\Gn^\e)\,\Gu\cdot\Gv\,\md S,
	\]
	and $\Gu\cdot\Gn^\e=0$ on $S\times\Gamma^\e_s$ yields \eqref{Id01}. Finally,
	H\"older's inequality gives
	\[
	|b^\e(\Gu,\Gv,\Gw)|\leq\|\Gu\|_{L^4(\O^\e_p)}\|\nabla\Gv\|_{L^2(\O^\e_p)}\|\Gw\|_{L^4(\O^\e_p)},
	\]
	and \eqref{117} follows from the uniform embedding $H^1(\O^\e_p)\hookrightarrow
	L^4(\O^\e_p)$ of Lemma~\ref{Embed01}.
\end{proof}

	\section{Existence and a priori estimates}\label{Sec05}
	This section presents the main existence result and a priori estimates for weak
	solutions of the NSCH system \eqref{MainP01}. First, we derive a priori
	estimates and the associated energy dissipation law for weak solutions. Then, we
	prove existence by means of a Galerkin approximation.
	
	\subsection{A priori estimates and energy dissipation law}\label{SSec51}
	\begin{lemma}\label{Est01}
		Let $(\Gu^\e,\phi^\e,\mu^\e)\in \mathfrak{U}^{\varepsilon } \times \mathfrak{C}^{\varepsilon} \times \mathfrak{W}^{\varepsilon}$ be a weak solution of  \eqref{MainP01} in the sense of Definition \ref{Def01}. Then, we have the following a priori estimates
		\begin{equation}\label{MainE01}
			\begin{aligned}
				\|\Gu^\e\|_{L^\infty(S;L^2(\O_p^\e))}
				+ \|\Gu^\e\|_{L^2(S; H^1(\O_p^\e))}+\|\partial_t\Gu^\e\|_{L^{4/3}(S;\GH^{-1}_{\div}(\O^\e_p))}&\leq C\sqrt{\lambda^\e},\\
				\|\phi^\e\|_{L^\infty(S;H^1(\O^\e_p))}+\|\partial_t\phi^\e\|_{L^2(S;H^{-1}(\O^\e_p))}+\|\mu^\e\|_{L^2(S;H^1(\O^\e_p))}&\leq C.
			\end{aligned}
		\end{equation}
The constants are independent of $\e$.
	\end{lemma}
	\begin{proof}
		The proof is divided into $8$ steps given below.
		
		{\bf Step 1.} Since the momentum equation is tested against functions in
		\(L^4(S;\GH^1_{\div}(\O^\e_p))\), the choice
		\(\varphi_1=\Gu^\e\) is not made directly in the weak formulation. Instead,
		the momentum part of the energy estimate is obtained at the Galerkin level,
		where testing by the approximate velocity is admissible, and then by passing
		to the limit. Equivalently, one may use a standard time-regularisation
		argument. In this way, using the skew-symmetry \eqref{skew01} of the
		convective term, we obtain, for a.e. \(t\in S\),
		\begin{multline}
			\label{est-momentum}
			\frac12\|\Gu^\e(t)\|^2_{L^2(\O^\e_p)}
			-\frac12\|\Gu^\e_0\|^2_{L^2(\O^\e_p)}
			+\int_0^t\!\!\int_{\O^\e_p}
			\A^\e D(\Gu^\e):D(\Gu^\e)\,\md x\,\md s \\
			+\int_0^t\!\!\int_{\O^\e_p}
			\lambda^\e \phi^\e\nabla\mu^\e\cdot\Gu^\e\,\md x\,\md s
			\leq
			\int_0^t\!\!\int_{\O^\e_p}\Gg^\e\cdot\Gu^\e\,\md x\,\md s .
		\end{multline}
		On the other hand, testing \eqref{weak-phi1} with
		\(\lambda^\e\mu^\e\in L^2(S;H^1(\O^\e_p))\), using
		\[
		\mu^\e=-\Delta\phi^\e+F'(\phi^\e)
		\quad\text{in }H^{-1}(\O^\e_p),
		\]
		and applying the standard chain-rule identities, gives, for a.e. \(t\in S\),
		\begin{multline}
			\label{est-phi-mu}
			\lambda^\e\int_0^t\!\!\int_{\O^\e_p}
			|\nabla\mu^\e|^2\,\md x\,\md s
			+\lambda^\e\int_0^t\langle G(\phi^\e),\mu^\e\rangle\,\md s
			-\lambda^\e\int_0^t\!\!\int_{\O^\e_p}
			\phi^\e\Gu^\e\cdot\nabla\mu^\e\,\md x\,\md s \\
			+\frac{\lambda^\e}{2}\|\nabla\phi^\e(t)\|^2_{L^2(\O^\e_p)}
			-\frac{\lambda^\e}{2}\|\nabla\phi_0\|^2_{L^2(\O^\e_p)}
			+\lambda^\e\int_{\O^\e_p}F(\phi^\e(t))\,\md x
			-\lambda^\e\int_{\O^\e_p}F(\phi_0)\,\md x
			=0 .
		\end{multline}
		
		{\bf Step 2.} Adding \eqref{est-momentum} and \eqref{est-phi-mu}, the two
		capillary transport terms cancel. Using the initial-data bounds
		\eqref{ICA01}--\eqref{F01}, the coercivity \eqref{Coe01} of \(\A^\e\), and
		the Sobolev embedding Lemma~\ref{Embed01}, we obtain, for a.e. \(t\in S\),
		\begin{multline}
			\label{4.6}
			\frac12\|\Gu^\e(t)\|^2_{L^2(\O^\e_p)}
			+\kappa_1\|D(\Gu^\e)\|^2_{L^2((0,t)\times\O^\e_p)}
			+\lambda^\e\|\nabla\mu^\e\|^2_{L^2((0,t)\times\O^\e_p)}
			+\frac{\lambda^\e}{2}\|\nabla\phi^\e(t)\|^2_{L^2(\O^\e_p)} \\
			+\lambda^\e\int_{\O^\e_p}F(\phi^\e(t))\,\md x
			+\lambda^\e\int_0^t\langle G(\phi^\e),\mu^\e\rangle\,\md s
			\leq
			C\lambda^\e
			+\int_0^t\!\!\int_{\O^\e_p}\Gg^\e\cdot\Gu^\e\,\md x\,\md s .
		\end{multline}
		By H\"older's and Young's inequalities, together with Korn's inequality
		(Lemma~\ref{Korn}), for every \(\delta>0\) we have
		\[
		\left|
		\int_0^t\!\!\int_{\O^\e_p}\Gg^\e\cdot\Gu^\e\,\md x\,\md s
		\right|
		\leq
		\delta\|D(\Gu^\e)\|^2_{L^2((0,t)\times\O^\e_p)}
		+
		C_\delta\|\Gg^\e\|^2_{L^2((0,t)\times\O^\e_p)} .
		\]
		Using \(\|\Gg^\e\|_{L^2(S\times\O^\e_p)}\leq C\sqrt{\lambda^\e}\), choosing
		\(\delta=\kappa_1/2\), and absorbing the corresponding term into the
		left-hand side of \eqref{4.6}, we arrive at
		\begin{multline}
			\label{4.7}
			\frac12\|\Gu^\e(t)\|^2_{L^2(\O^\e_p)}
			+\frac{\kappa_1}{2}\|D(\Gu^\e)\|^2_{L^2((0,t)\times\O^\e_p)}
			+\lambda^\e\|\nabla\mu^\e\|^2_{L^2((0,t)\times\O^\e_p)}
			+\frac{\lambda^\e}{2}\|\nabla\phi^\e(t)\|^2_{L^2(\O^\e_p)} \\
			+\lambda^\e\int_{\O^\e_p}F(\phi^\e(t))\,\md x
			+\lambda^\e\int_0^t\langle G(\phi^\e),\mu^\e\rangle\,\md s
			\leq C\lambda^\e ,
		\end{multline}
		where \(C>0\) is independent of \(\e\) and $T$. Since \(t\in S\) is arbitrary, this
		gives the corresponding \(L^\infty(S;L^2)\) and \(L^2(S;H^1)\) bounds.
		
		{\bf Step 3.} We estimate $\lambda^\e\int_{S\times\O_p^\e}G(\phi^\e)\mu^\e\,\md(x,t)$
		from below. Recalling $\mu^\varepsilon=-\Delta\phi^\varepsilon+F'(\phi^\varepsilon)$
		and using integration by parts with the Neumann condition \eqref{MainP01}$_6$,
		we decompose for a.e.\ $t\in(0,T)$:
		\begin{equation}\label{E01}
			\int_{\O_p^\e}G(\phi^\e)\mu^\e\,\md x
			= \underbrace{\int_{\O_p^\e}G'(\phi^\e)|\nabla\phi^\e|^2\,\md x}_{=:I_1}
			+ \underbrace{\int_{\O_p^\e}G(\phi^\e)F'(\phi^\e)\,\md x}_{=:I_2}.
		\end{equation}
		From \eqref{AssG01}, $I_1$ satisfies
		\begin{equation}\label{E02}
			c_1\|\nabla\phi^\e\|^2_{L^2(\O^\e_p)}\leq I_1\leq c_2\|\nabla\phi^\e\|^2_{L^2(\O^\e_p)}.
		\end{equation}
		The assumption \eqref{AssG01} and the fundamental theorem of calculus give
		\begin{equation}\label{ass:S-structure}
			c_1|s|\leq|G(s)|\leq c_2|s|, \qquad
			c_1 s^2\leq G(s)s\leq c_2 s^2 \qquad\forall s\in\R.
		\end{equation}
		A pointwise case analysis ($|s|\geq 1$ and $|s|<1$) using \eqref{ass:S-structure}$_1$
		gives
		\begin{equation}\label{E03}
			c_1\int_{\O^\e_p}\!\left((\phi^\e)^2-\tfrac{1}{2}\right)^2\md x
			-\frac{c_2}{4}|\O^\e_p|\leq I_2.
		\end{equation}
		Combining \eqref{E01}--\eqref{E03} with \eqref{4.7} and dropping non-negative
		left-hand-side terms, we conclude
		\begin{equation}\label{U131}
			\begin{aligned}
				\|\Gu^\e\|_{L^\infty(S;L^2(\O^\e_p))}
				+\|\Gu^\e\|_{L^2(S;H^1(\O^\e_p))}
				&\leq C\sqrt{\lambda^\e},\\
				\|\nabla\mu^\e\|_{L^2(S\times\O^\e_p)}
				+\|\nabla\phi^\e\|_{L^\infty(S;L^2(\O^\e_p))}
				+\int_{\O_p^\varepsilon}F(\phi^\e)\,\md x
				+\int_{S\times\O^\e_p}\!\left((\phi^\e)^2-\tfrac{1}{2}\right)^2\md(x,t)
				&\leq C,
			\end{aligned}
		\end{equation}
		with constants independent of $\e$.
		
		{\bf Step 4.} We prove that there exists a constant $C$ independent of  $\e$ such that
		\begin{equation}\label{IIMb}
			\|\bar \phi^\e\|_{L^\infty((0,T))}\leq C ,\quad\text{where}\quad \wo{\phi}^\e=\frac{1}{|\O^\e_p|}\int_{\O^\e_p}\phi^\e\,\md x.
		\end{equation}
		Using the Neumann boundary conditions in \eqref{MainP01}$_6$ and the
		divergence-free condition in \eqref{MainP01}$_2$, and taking the mean over $\O^\e_p$ on the equation \eqref{MainP01}$_4$, with the above conditions, we get
		\[
		\frac{d}{dt} \wo{\phi}^\e +G({\bar\phi^\e})+\frac{1}{|\O^\e_p|}\int_{\O^\e_p}\left[G(\phi^\e)-G({\bar\phi^\e})\right]\,\md x=\frac{d}{dt} \wo{\phi}^\e +\frac{1}{|\O^\e_p|}\int_{\O^\e_p}G(\phi^\e)\,\md x=0.
		\]
		By assumption, $G$ is Lipschitz continuous with Lipschitz constant $c_2$. Therefore,
		\begin{equation}\label{Ext01}
			|r^\e(t)|=\left|\frac{1}{|\O^\e_p|}\int_{\O^\e_p}\left[G(\phi^\e)-G({\bar\phi^\e})\right]\,\md x\right|\leq \frac{c_2}{|\O^\e_p|^{1/2}}\|\phi^\e-\bar\phi^\e\|_{L^2(\O^\e_p)}\le \frac{C}{|\O^\e_p|^{1/2}}\|\nabla\phi^\e\|_{L^2(\O^\e_p)},
		\end{equation}
		due to the fact that $\phi^\e-\bar\phi^\e$ has mean zero and we applied the Poincar\'e-Wirtinger inequality Lemma \ref{LemPW}. So, the ODE for the mean satisfies 
		$$\frac{d}{dt} \wo{\phi}^\e +G({\bar\phi^\e})+r^\e(t)=0.$$
		Since $G$ satisfies \eqref{ass:S-structure}, we have
		$$G(s)\mathrm{sgn}(s)\geq c_1|s|,\quad \forall\,s\in \R.$$
		Then, multiplying by $\mathrm{sgn}(\wo{\phi}^\e)$, we get the above inequality for $|\wo{\phi}^\e|$ as
		$$ \frac{d}{dt} |\wo{\phi}^\e| \leq  -c_1|\bar\phi^\e|+|r^\e(t)|,$$
		which implies, by Gronwall's inequality,
		$$ |\bar\phi^\e(t)|\leq |\bar\phi_0|\exp(-c_1t)+\int_0^t\exp(-c_1(t-s))|r^\e(s)|\,ds\leq |\bar\phi_0|+\sup_{s\in[0,T]}|r^\e(s)|\int_0^t\exp(-c_1(t-s))\,ds,$$
		where 	$\overline{\phi^\e}(0)
		=
		\frac{1}{|\O_p^\e|}
		\int_{\O_p^\e}\phi_0\,dx$.
		Since $|\O_p^\e|$ is uniformly bounded away from zero,
		\[
		\bigl|\overline{\phi^\e}(0)\bigr|
		\leq
		|\O_p^\e|^{-1/2}\|\phi_0\|_{L^2(\O)}
		\leq C.
		\]
		Consequently, $\|\overline{\phi^\e}\|_{L^\infty(S)}\leq C$,
		with $C$ independent of $\e$.
		 By \eqref{Ext01}, we have
		$$\sup_{s\in[0,T]}|r^\e(s)|\leq \frac{C}{|\O^\e_p|^{1/2}}\|\nabla\phi^\e\|_{L^\infty(S;L^2(\O^\e_p))}\leq C\sqrt\frac{2|Y|}{|Y_p||\O|}\|\nabla\phi^\e\|_{L^\infty(S;L^2(\O^\e_p))}\leq C_1,$$
		using the estimate \eqref{U131}$_4$. This implies, together with \eqref{ICA01} 
		$$\|\bar\phi^\e\|_{L^\infty(S)}\leq |\bar\phi_0|+C_2.$$
		The upper bound is independent of $\e$. Hence we obtain \eqref{IIMb}.
		
		Then, the estimate \eqref{U131}$_4$ together with the estimate \eqref{IIMb} and Poincar\'e-Wirtinger inequality (Lemma \ref{LemPW}) give
		\begin{equation}\label{133}
			\|\phi^\e\|_{L^\infty(S;H^1(\O^\e_p))}\leq C.
		\end{equation}
The constants are independent of $\e$.
		
		{\bf Step 5.} Estimate for $\partial_t \Gu^\varepsilon$ in
		$L^{4/3}(0,T;H^{-1}_{\div}(\Omega^\varepsilon_p))$.
		
		Take arbitrary $\psi\in H^1_{\div}(\Omega^\varepsilon_p)$ with
		$\|\psi\|_{H^1(\Omega^\varepsilon_p)}\le 1$. Using \eqref{weak-momentum1},
		the embedding $H^1(\Omega^\varepsilon_p)\hookrightarrow L^4(\Omega^\varepsilon_p)$
		(Lemma~\ref{Embed01}), Korn's inequality, H\"older's inequality, and \eqref{Ass01},
		we obtain
		\[
		\|\partial_t \Gu^\varepsilon\|_{H^{-1}_{\div}(\Omega^\varepsilon_p)}
		\le C\Bigl(\|\Gu^\e\|^{1/2}_{L^2(\O^\e_p)}\|\Gu^\varepsilon\|^{3/2}_{H^1(\Omega^\varepsilon_p)}
		+ \|\Gu^\varepsilon\|_{H^1(\Omega^\varepsilon_p)}
		+ \lambda^\varepsilon\|\varphi^\varepsilon\|_{H^1(\Omega^\varepsilon_p)}
		\|\nabla\mu^\varepsilon\|_{L^2(\Omega^\varepsilon_p)}
		+ \kappa_5\sqrt{\lambda^\varepsilon}\Bigr).
		\]
		Raising to the power $4/3$, integrating in time, and applying the a priori
		estimates \eqref{U131}, \eqref{133}, the dominant convective contribution
		satisfies (see \cite{temam1977})
		\[
		\int_0^T\|\Gu^\e\|^{2/3}_{L^2(\O^\e_p)}\|\Gu^\varepsilon\|^2_{H^1(\O^\e_p)}\,dt
		\le \Bigl(\operatorname{ess\,sup}_{t}\|\Gu^\varepsilon\|_{L^2(\O^\e_p)}^{2/3}\Bigr)
		\int_0^T\|\Gu^\varepsilon\|_{H^1(\O^\e_p)}^{2}\,dt
		\le C(\lambda^\varepsilon)^{4/3},
		\]
		while the remaining terms contribute at most $C\sqrt{\lambda^\varepsilon}$
		after H\"older's inequality. Consequently,
		\[
		\|\partial_t \Gu^\varepsilon\|_{L^{4/3}(0,T;\,H^{-1}_{\div}(\Omega^\varepsilon_p))}
		\le C\sqrt{\lambda^\varepsilon},
		\]
		with $C$ independent of $\varepsilon$.
		
		{\bf Step 6.} Estimate for $\partial_t\phi^\e$ in $L^2(S;H^{-1}(\O^\e_p))$.
		
		From \eqref{weak-phi1} with arbitrary $\varphi_2\in H^1(\O^\e_p)$, using
		\eqref{ass:S-structure}, $H^1(\O^\e_p)\hookrightarrow L^4(\O^\e_p)$, and
		H\"older's inequality,
		\[
		\sup_{\|\varphi_2\|_{H^1(\O^\e_p)}\leq 1}
		\bigl|\langle\partial_t\phi^\e,\varphi_2\rangle\bigr|
		\leq C\bigl(\|\phi^\e\|_{L^2(\O^\e_p)}
		+\|\nabla\mu^\e\|_{L^2(\O^\e_p)}
		+\|\phi^\e\|_{H^1(\O^\e_p)}\|\Gu^\e\|_{H^1(\O^\e_p)}\bigr).
		\]
		Squaring, integrating in time, and applying \eqref{U131}, \eqref{133} give
		\begin{equation}\label{136}
			\|\partial_t\phi^\e\|_{L^2(S;\,H^{-1}(\O^\e_p))}\le C,
		\end{equation}
		with constant independent of $\e$.
		
		{\bf Step 7.} $L^2$-bound for $\mu^\e$.
		
		Testing the chemical potential identity $\mu^\e=-\Delta\phi^\e+F'(\phi^\e)$
		with $1$ and using the Neumann condition \eqref{MainP01}$_6$ gives
		$\bar\mu^\e = |\O^\e_p|^{-1}\int_{\O^\e_p}F'(\phi^\e)\,\md x$.
		Since $|F'(s)|\leq C(|s|^3+1)$, estimates \eqref{133} and Lemma~\ref{Embed01}
		yield $\|\bar\mu^\e\|_{L^\infty(S)}\leq C$. Combined with \eqref{U131}$_3$
		and the Poincar\'e--Wirtinger inequality,
		\begin{equation}\label{137}
			\|\mu^\e\|_{L^2(S;\,H^1(\O^\e_p))}\leq C,
		\end{equation}
		with constants independent of $\e$.
		
		Thus we obtain \eqref{MainE01} from the above estimates, for a weak solution $(\Gu^\e,\phi^\e,\mu^\e)\in \mathfrak{U}^{\varepsilon } \times \mathfrak{C}^{\varepsilon} \times \mathfrak{W}^{\varepsilon}$.

	\end{proof}
	
\begin{lemma}[Global fields]\label{Lem08+}
	Let $(\Gu^\e,\phi^\e,\mu^\e)\in\mathfrak U^\e\times\mathfrak C^\e\times\mathfrak W^\e$
	be a weak solution of \eqref{MainP01}. Let $\wt\Gu^\e:=\mathcal E^\e\Gu^\e$ be the
	primal velocity extension (Lemma~\ref{LemD02}), $\overline\Gu^\e$ its zero
	extension through the obstacles, $\wt\phi^\e:=E^\e\phi^\e$ and
	$\wt\mu^\e:=E^\e\mu^\e$ the primal scalar extensions (Lemma~\ref{Emain}), and
	$\overline\phi^\e$ the zero extension of $\phi^\e$. Then
	$\wt\Gu^\e|_{\O_p^\e}=\overline\Gu^\e|_{\O_p^\e}=\Gu^\e$,
	$\wt\phi^\e|_{\O_p^\e}=\overline\phi^\e|_{\O_p^\e}=\phi^\e$,
$\wt\mu^\e|_{\O_p^\e}=\mu^\e$, and
\[
\div\wt\Gu^\e
=
\div\overline\Gu^\e
=
0
\quad\text{in }\mathcal D'(S\times\O).
\]
Moreover
	\begin{equation}\label{EXEst01}
		\begin{aligned}
			\|\overline\Gu^\e\|_{L^\infty(S;L^2(\O))}
			+\|\wt\Gu^\e\|_{L^2(S;H^1(\O))}
			&\leq C\sqrt{\lambda^\e},\\
			\|\wt\phi^\e\|_{L^\infty(S;H^1(\O))}
			+\|\partial_t\overline\phi^\e\|_{L^2(S;H^{-1}(\O))}
			+\|\wt\mu^\e\|_{L^2(S;H^1(\O))}
			&\leq C,
		\end{aligned}
	\end{equation}
	with $C>0$ independent of $\e$. Furthermore,
	\begin{equation}\label{GlobalVelocityExtensionsDifference}
		\|\wt\Gu^\e-\overline\Gu^\e\|_{L^2(S\times\O)}\leq C\,\tfrac{a_\e}{\e}\,\sqrt{\lambda^\e},
		\qquad
		\|\wt\phi^\e-\overline\phi^\e\|_{L^\infty(S;L^2(\O))}\leq C\,\tfrac{a_\e}{\e},
	\end{equation}
	so $\tfrac{1}{\sqrt{\lambda^\e}}(\wt\Gu^\e-\overline\Gu^\e)\to0$ strongly in
	$L^2(S\times\O)^3$ and $\wt\phi^\e-\overline\phi^\e\to0$ strongly in
	$L^\infty(S;L^2(\O))$; the normalized primal and zero velocity extensions, and
	the primal and zero scalar extensions, therefore share the same limit whenever
	either exists.
\end{lemma}
\begin{proof}
	The identity and divergence properties follow from
	Lemmas~\ref{Emain} and~\ref{LemD02}. The spatial estimates in
	\eqref{EXEst01} follow from the pore-domain estimates
	\eqref{MainE01} and the uniform extension bounds of
	Lemmas~\ref{Emain} and~\ref{LemD02}.
	
	For the scalar time derivative, consider the restriction operator
	\[
	r^\e:H^1(\O)\longrightarrow H^1(\O_p^\e),
	\qquad
	r^\e\eta:=\eta|_{\O_p^\e}.
	\]
	Since $\|r^\e\|\leq1$, its dual operator
	\[
	(r^\e)^*:H^{-1}(\O_p^\e)\longrightarrow H^{-1}(\O)
	\]
	has norm at most one. Moreover, $\partial_t\overline\phi^\e
	=
	(r^\e)^*(\partial_t\phi^\e)$
	in $L^2(S;H^{-1}(\O))$. Consequently,
	\[
	\|\partial_t\overline\phi^\e\|_{L^2(S;H^{-1}(\O))}
	\leq
	\|\partial_t\phi^\e\|_
	{L^2(S;H^{-1}(\O_p^\e))}
	\leq C.
	\]
	
	The estimates in
	\eqref{GlobalVelocityExtensionsDifference} follow from
	Lemmas~\ref{lem:scalar-two-extensions} and~\ref{LemD02}, together
	with \eqref{MainE01}. Since $\frac{a_\e}{\e}=\e^{\alpha-1}\longrightarrow0$,
	they imply the stated strong convergences and the equality of the
	corresponding limits.
\end{proof}
	
	\begin{lemma}\label{Lem13}
		Let $(\Gu^\e,\phi^\e,\mu^\e)
		\in \mathfrak{U}^{\varepsilon } \times
		\mathfrak{C}^{\varepsilon} \times \mathfrak{W}^{\varepsilon}$
		be a weak solution of \eqref{MainP01} in the sense of
		Definition~\ref{Def01}. Define
		\begin{equation}\label{ExpE01}
			\GT^\e(\Gu^\e,\phi^\e)
			=\GT^\e_K(\Gu^\e)+\GT^\e_F(\phi^\e)=
			\frac{1}{2}\|\Gu^\e\|^2_{L^2(\O^\e_p)}
			+\frac{\lambda^\e}{2}\|\nabla \phi^\e\|^2_{L^2(\O^\e_p)}
			+\lambda^\e\int_{\O_p^\varepsilon} F(\phi^\e)\,\md x .
		\end{equation}
		Then, in the sense of distributions in \(S=(0,T)\), the energy satisfies
		\begin{equation}\label{EDL}
			\frac{d}{dt}\GT^\e(\Gu^\e,\phi^\e)
			+\frac{\kappa_1}{2}\|D(\Gu^\e)\|^2_{L^2(\O_p^\e)}
			+\lambda^\e\|\nabla\mu^\e\|^2_{L^2(\O_p^\e)}
			+\lambda^\e\langle G(\phi^\e),\mu^\e\rangle
			\leq
			C\|\Gg^\e\|^2_{L^2(\O_p^\e)} 
		\end{equation}
		where the constant \(C>0\) are independent of \(\e\). Moreover,
		\begin{equation}\label{EstEn01}
			\GT^\e(\Gu^\e,\phi^\e)\leq C\lambda^\e .
		\end{equation}
		The constant \(C\) is independent of \(\e\).
	\end{lemma}
	
	\begin{proof}
		The energy estimate is first obtained at the Galerkin level, where testing
		the momentum equation by the approximate velocity is admissible. Passing to
		the limit gives the corresponding inequality for weak solutions. Proceeding
		as in Lemma~\ref{Est01}, Step~1, and using the skew-symmetry of the
		convective term, we obtain
		\begin{equation*}
			\frac{1}{2}\frac{d}{dt}\|\Gu^\e\|^2_{L^2(\O^\e_p)}
			+\int_{\O_p^\e}\A^\e D(\Gu^\e):D(\Gu^\e)\,\md x
			+\lambda^\e\int_{\O_p^\e}
			\phi^\e\nabla\mu^\e\cdot\Gu^\e\,\md x
			\leq
			\int_{\O_p^\e}\Gg^\e\cdot\Gu^\e\,\md x .
		\end{equation*}
		Testing the Cahn--Hilliard equation with \(\lambda^\e\mu^\e\), using \eqref{MainP01}$_5$ holds in $H^{-1}(\O^\e_p)$, and applying the standard chain-rule identities gives
		\begin{equation*}
			\lambda^\e\|\nabla\mu^\e\|^2_{L^2(\O_p^\e)}
			+\lambda^\e\langle G(\phi^\e),\mu^\e\rangle
			-\lambda^\e\int_{\O_p^\e}
			\phi^\e\Gu^\e\cdot\nabla\mu^\e\,\md x
			+\frac{\lambda^\e}{2}\frac{d}{dt}
			\|\nabla\phi^\e\|^2_{L^2(\O_p^\e)}
			+\lambda^\e\frac{d}{dt}
			\int_{\O_p^\e}F(\phi^\e)\,\md x
			=0 .
		\end{equation*}
		Adding the last two relations, the capillary transport terms cancel.
		Using the coercivity of \(\A^\e\), we get
		\begin{equation*}
			\frac{d}{dt}\GT^\e(\Gu^\e,\phi^\e)
			+\kappa_1\|D(\Gu^\e)\|^2_{L^2(\O_p^\e)}
			+\lambda^\e\|\nabla\mu^\e\|^2_{L^2(\O_p^\e)}
			+\lambda^\e\langle G(\phi^\e),\mu^\e\rangle
			\leq
			\int_{\O_p^\e}\Gg^\e\cdot\Gu^\e\,\md x .
		\end{equation*}
		By H\"older's inequality, Korn's inequality, and Young's inequality,
		\[
		\left|\int_{\O_p^\e}\Gg^\e\cdot\Gu^\e\,\md x\right|
		\leq
		\frac{\kappa_1}{2}\|D(\Gu^\e)\|^2_{L^2(\O_p^\e)}
		+C\|\Gg^\e\|^2_{L^2(\O_p^\e)} .
		\]
		Hence, we obtain \eqref{EDL}.
		This completes the proof.
	\end{proof}
	
	\begin{remark}
		Observe that $\GT^\e_K$ is the kinetic energy and $\GT_F^\e$ is the bulk free energy. Moreover, the advection term $\Gu^\e\cdot\nabla\phi^\e$ and the surface tension term $\phi^\e\nabla\mu^\e$ do not contribute to the total energy $\GT^\e$.  That is, they satisfy the zero-energy-contribution property,
		ensuring that the dissipation of order parameter and the dissipation of momentum
		remain decoupled. 
	\end{remark}

	\subsection{Existence of a weak solution}\label{SSec52}
	
	Below, we give the existence of weak solution for the microscopic NSCH system \eqref{MainP01}.
	\begin{theorem}\label{Th1}
		There exists a weak solution $(\Gu^\e,\phi^\e,\mu^\e)\in \mathfrak{U}^{\varepsilon } \times \mathfrak{C}^{\varepsilon} \times \mathfrak{W}^{\varepsilon}$ of the problem \eqref{MainP01} in the sense of Definition~\ref{Def01}.
	\end{theorem}
	\begin{proof}
		We prove the existence by a four-step Galerkin procedure.
		
		\textbf{Step 1. Galerkin basis.}
		
		Set $V:=\GH^1_{\div}(\O^\e_p)$ and $H:=\overline{V}^{\,L^2(\O^\e_p)^3}$.
		By Korn's inequality and \eqref{Coe01}, we have
		$$\kappa_1 a_0^\e(\bw,\bw)\leq a_t^\e(\bw,\bw)\leq \kappa_0 a_0^\e(\bw,\bw),\quad\forall\, \bw\in V,$$
		where
		$a^\e_0(\bw,\bv):=\int_{\O^\e_p}D(\bw):D(\bv)\,\md x$ and $a_t^\e(\bw,\bv)=\int_{\O^\e_p}\A^\e(t,x)D(\bw):D(\bv)\,\md x$. The bilinear form $a_0^\e$ is coercive on $V$, so
		the associated compact self-adjoint operator on $H$ (via Lax--Milgram and
		the canonical injection $J:H\to V^\ast$) admits a spectral basis
		$\{\bpsi^\e_k\}_{k\in\N}\subset V$ that is $L^2(\O^\e_p)^3$-orthonormal,
		$a^\e_0$-orthogonal, divergence-free, vanishes on $\partial\O$, and satisfies
		$\bpsi^\e_k\cdot\Gn^\e=0$ on $\Gamma^\e_s$.
		We also fix an $L^2$-orthonormal basis $\{\zeta_i\}_{i\in\N}\subset H^1(\O^\e_p)$,
		orthogonal in $H^1$, and set
		$V_n:=\mathrm{span}\{\bpsi^\e_1,\dots,\bpsi^\e_n\}$,
		$W_n:=\mathrm{span}\{\zeta_1,\dots,\zeta_n\}$.
		
		\textbf{Step 2. Galerkin approximation and reduction to ODEs.}
		
		For fixed $\varepsilon>0$ (suppressing superscripts on the basis), we seek
		\[
		\mathbf u_n^\e = \textstyle\sum_{j=1}^n a_j^n(t)\bpsi_j,\qquad
		\phi_n^\e   = \textstyle\sum_{j=1}^n b_j^n(t)\zeta_j,\qquad
		\mu_n^\e    = \textstyle\sum_{j=1}^n c_j^n(t)\zeta_j,
		\]
		with initial coefficients given by the Leray $L^2$-projections of
		$\mathbf u^\e_0,\phi_0$ onto $V_n,W_n$. The Galerkin system reads: for
		a.e.\ $t\in S$ and all $\varphi_1\in V_n$, $\varphi_2,\varphi_3\in W_n$,
		\begin{align}
			\label{weak-momentum}
			\langle\partial_t\mathbf u_n^\e,\varphi_1\rangle
			+\int_{\O^\e_p}(\mathbf u_n^\e\cdot\nabla)\mathbf u_n^\e\cdot\varphi_1\,\md x
			+a^\e_t(\mathbf u_n^\e,\varphi_1)
			+\int_{\O^\e_p}\lambda^\e\phi_n^\e\nabla\mu_n^\e\cdot\varphi_1\,\md x
			&=\int_{\O^\e_p}\mathbf g^\e\cdot\varphi_1\,\md x,\\
			\label{weak-phi}
			\langle\partial_t\phi_n^\e,\varphi_2\rangle
			+\int_{\O^\e_p}G(\phi_n^\e)\varphi_2\,\md x
			+\int_{\O^\e_p}\nabla\mu_n^\e\cdot\nabla\varphi_2\,\md x
			&=\int_{\O^\e_p}\phi_n^\e\mathbf u_n^\e\cdot\nabla\varphi_2\,\md x,\\
			\label{weak-mu}
			-\int_{\O^\e_p}\mu_n^\e\varphi_3\,\md x
			+\int_{\O^\e_p}\nabla\phi_n^\e\cdot\nabla\varphi_3\,\md x
			+\int_{\O^\e_p}F'(\phi_n^\e)\varphi_3\,\md x
			&=0.
		\end{align}
		Testing with $\varphi_1=\psi_i$, $\varphi_2=\varphi_3=\zeta_i$ and using
		the orthonormality of $\{\zeta_i\}$ to eliminate $c^n$ via
		$Sc^n=Db^n+F(b^n)$ reduces the system to
		\begin{equation}\label{eq:coupled-system}
			\dot a^n=\Phi_1(a^n,b^n,t),\qquad \dot b^n=\Phi_2(a^n,b^n,t),
		\end{equation}
		where $\Phi_1,\Phi_2$ are defined by the standard matrix expressions against
		the bases $\{\bpsi_i\},\{\zeta_i\}$ (see \cite[Chapter~I]{coddington1956theory}).
		Since $F'\in\mathcal{C}^\infty(\R)$ and $G\in\mathcal{C}^1(\R)$, the maps
		$\Phi_1,\Phi_2$ are measurable in $t$ and locally Lipschitz in $(a^n,b^n)$.
		The Carathéodory theorem \cite[Chapter~I]{coddington1956theory} yields a unique
		maximal solution $(a^n,b^n)\in W^{1,1}_{\mathrm{loc}}([0,T_n^\ast);\R^{2n})$,
		which extends to $[0,T)$ once the uniform bounds of Step~3 preclude blow-up.
		
		\textbf{Step 3. Global-in-time bounds.}
		
		Since $\mathbf u_n^\e\in V_n$ and $\phi_n^\e,\mu_n^\e\in W_n$, proceeding
		exactly as in Lemma~\ref{Est01} (Steps~1--3) and Lemma~\ref{Lem13} gives
		\begin{multline}\label{316}
			\frac{1}{2}\frac{d}{dt}\|\mathbf u_n^\e\|^2_{L^2(\O^\e_p)}
			+c\|\Gu^\e_n\|^2_{H^1(\O^\e_p)}
			+\lambda^\e\|\nabla\mu_n^\e\|^2_{L^2(\O^\e_p)}
			+\frac{\lambda^\e}{2}\frac{d}{dt}\|\nabla\phi_n^\e\|^2_{L^2(\O^\e_p)}
			+\lambda^\e\frac{d}{dt}\int_{\O^\e_p}F(\phi_n^\e)\,\md x\\
			+\lambda^\e c_1\|\nabla\phi^\e_n\|^2_{L^2(\O^\e_p)}
			+\lambda^\e c_1\int_{\O^\e_p}\!\left((\phi^\e_n)^2-\tfrac{1}{2}\right)^2\md x
			\leq C(\e).
		\end{multline}
		Define the energy $E_n(t):=\frac{1}{2}\|\mathbf u_n^\e\|^2_{L^2}
		+\frac{\lambda^\e}{2}\|\nabla\phi_n^\e\|^2_{L^2}
		+\lambda^\e\int_{\O^\e_p}F(\phi_n^\e)\,\md x$.
		Dropping $\lambda^\e\|\nabla\mu_n^\e\|^2_{L^2}$ and using the algebraic identity
		$(s^2-\tfrac12)^2=4F(s)+s^2-\tfrac34$ to absorb the last term of \eqref{316}
		into $E_n$, together with the Poincaré bound
		$\kappa_p\|\Gu^\e_n\|^2_{L^2}\leq\|\Gu^\e_n\|^2_{H^1}$, yields
		\begin{equation}\label{eq:en-diff}
			\frac{d}{dt}E_n(t)+\beta E_n(t)\leq C(\e),
			\qquad\beta:=\min\{2c\kappa_p,2c_1,4c_1\}>0,
		\end{equation}
		with $\beta,C(\e)$ independent of $n$ and $T_n$. Gronwall's lemma gives
		\begin{equation}\label{eq:en-uniform-214}
			E_n(t)\leq\max\!\left\{E_n(0),\,\frac{C(\e)}{\beta}\right\}=:C(\e)
			\qquad\forall\,t\in(0,T_n),
		\end{equation}
		where $E_n(0)\leq C\lambda^\e$ uniformly in $n$ by \eqref{ICA01}.
		Hence $\|\mathbf u_n^\e\|_{L^2}+\|\nabla\phi_n^\e\|_{L^2}\leq C(\e)$, and
		the Poincaré--Wirtinger inequality (as in Step~4 of Lemma~\ref{Est01}) gives
		$\|\phi^\e_n\|_{L^2}\leq C(\e)$, both independent of $n$ and $T_n$.
		Since $\{\psi_i\},\{\zeta_i\}$ are $L^2$-orthonormal,
		these bounds imply $\|(a^n(t),b^n(t))\|_{\R^{2n}}\leq C(\e)$, precluding
		blow-up; hence $T_n=T$ for each $n\in\N$.
		
		\textbf{Step 4. Uniform estimates and passage to the limit.}
		
		Integrating \eqref{316} over $(0,T)$ and arguing as in Lemma~\ref{Est01}
		gives, for a constant $C>0$ independent of $n$,
		\begin{equation}\label{214}
			\begin{aligned}
				&\|\mathbf u_n^\e\|_{L^\infty(S;L^2(\O^\e_p))}
				+\|\mathbf u_n^\e\|_{L^2(S;H^1(\O^\e_p))}
				+\|\partial_t\mathbf u_n^\e\|_{L^{4/3}(S;\GH^{-1}_{\div}(\O^\e_p))}
				\leq C,\\
				&\|\phi_n^\e\|_{L^\infty(S;H^1(\O^\e_p))}
				+\|\partial_t\phi_n^\e\|_{L^2(S;H^{-1}(\O^\e_p))}
				+\|\mu_n^\e\|_{L^2(S;H^1(\O^\e_p))}
				\leq C.
			\end{aligned}
		\end{equation}
		
		\noindent
		\textbf{Step 4.1: Compactness.}
		From \eqref{214} and the compact embedding
		$H^1(\O^\e_p)\subset\subset L^p(\O^\e_p)\hookrightarrow H^{-1}(\O^\e_p)$
		for $2\leq p<6$, we extract a subsequence (not relabelled) with
		\begin{equation}\label{SC0}
			\begin{aligned}
				\mathbf u_n^\e &\rightharpoonup\mathbf u^\e
				&&\text{weakly in }L^2(S;H^1(\O^\e_p))^3
				\text{ and weak-* in }L^\infty(S;L^2(\O^\e_p))^3,\\
				\phi_n^\e &\rightharpoonup\phi^\e
				&&\text{weakly in }L^2(S;H^1(\O^\e_p))
				\text{ and weak-* in }L^\infty(S;H^1(\O^\e_p)),\\
				\mu_n^\e &\rightharpoonup\mu^\e
				&&\text{weakly in }L^2(S;H^1(\O^\e_p)).
			\end{aligned}
		\end{equation}
		The Aubin--Lions--Simon lemma \cite{simon1987} applied to $\phi_n^\e$
		(exponents $(2,2)$) and to $\mathbf u_n^\e$ (exponents $(2,4/3)$) gives
		\begin{align}
			\label{SC01}
			\phi_n^\e &\to\phi^\e
			\quad\text{strongly in }L^2(S;L^p(\O^\e_p))\ \forall\,2\leq p<6
			\text{ and in }C([0,T];L^2(\O^\e_p)),\\
			\label{SC01u}
			\mathbf u_n^\e &\to\mathbf u^\e
			\quad\text{strongly in }L^2(S;L^2(\O^\e_p))^3
			\text{ and in }C([0,T];\GH^{-1}_{\div}(\O^\e_p)).
		\end{align}
		Interpolating \eqref{SC01u} with the $L^\infty(S;L^2)\cap L^2(S;H^1)$ bound
		upgrades to $\mathbf u_n^\e\to\mathbf u^\e$ strongly in $L^2(S;L^p(\O^\e_p))^3$
		for all $2\leq p<6$, and $\mathbf u^\e\in
		L^2(S;\GH^1_{\div}(\O^\e_p))\cap W^{1,4/3}(S;\GH^{-1}_{\div}(\O^\e_p))$.
		
			\textbf{Step 4.2: Nonlinear terms.}
	The Lipschitz continuity of $G$ and \eqref{SC01} imply
	\[
	G(\phi^\e_n)\to G(\phi^\e)
	\quad\text{strongly in }L^2(S;L^p(\O^\e_p)),
	\qquad 2\leq p<6.
	\]
	Moreover, interpolation between the convergence in
	$C([0,T];L^2(\O^\e_p))$ and the uniform bound in
	$L^\infty(S;L^6(\O^\e_p))$ gives
	\begin{equation}\label{SC01a}
		\phi^\e_n\to\phi^\e
		\quad\text{strongly in }C([0,T];L^p(\O^\e_p)),
		\qquad 2\leq p<6.
	\end{equation}
	Since $F'(s)=s^3-s$, the uniform $L^\infty(S;L^6)$ bound and
	the pointwise convergence yield
	\begin{equation}\label{SC02}
		F'(\phi_n^\e)\rightharpoonup F'(\phi^\e)
		\quad\text{weakly in }L^2(S\times\O^\e_p).
	\end{equation}
	Using \eqref{SC01a} with $p=4$ and the strong convergence
	$\mathbf u^\e_n\to\mathbf u^\e$ in $L^2(S;L^4(\O^\e_p))^3$, we obtain
	\begin{equation}\label{SC03}
		\phi_n^\e\mathbf u_n^\e\to\phi^\e\mathbf u^\e
		\quad\text{strongly in }L^2(S\times\O^\e_p)^3.
	\end{equation}
	Similarly, \eqref{SC01a} with $p=3$ and
	$\nabla\mu^\e_n\rightharpoonup\nabla\mu^\e$ in
	$L^2(S;L^2(\O^\e_p))^3$ give
	\begin{equation}\label{SC04}
		\phi_n^\e\nabla\mu_n^\e\rightharpoonup
		\phi^\e\nabla\mu^\e
		\quad\text{weakly in }
		L^2(S;L^{6/5}(\O^\e_p))^3.
	\end{equation}
	The convergences \eqref{SC0}--\eqref{SC04} allow us to pass to the limit in
	\eqref{weak-momentum}--\eqref{weak-mu}. Hence
	\[
	(\mathbf u^\e,\phi^\e,\mu^\e)
	\in
	\mathfrak U^\e\times\mathfrak C^\e\times\mathfrak W^\e
	\]
	satisfies \eqref{weak-momentum1}--\eqref{weak-mu1}.
		
		This completes the proof.
	\end{proof}
	

	
	
	\section{Unfolding operator and compactness  for a perforated domain}\label{SSec61}
	The main tool for homogenization in the perforated domain $\Omega_p^\varepsilon$
	and for the derivation of the two-scale limit system is the periodic unfolding
	operator adapted to perforated media, denoted by $\mathcal T_\varepsilon^\ast$.
	It was first introduced in \cite{cior03} and further developed in
	\cite{cioranescu2012periodic}. For a detailed presentation we refer to
	\cite[Chapter~4]{CDG}.
	
	We introduce two unfolding operators. The first resolves oscillations at
	the periodicity scale $\e$, while the second resolves the local behavior
	around the obstacles at the scale
	$a_\e=\e^\a=\e\delta_\e,$ 
	$ \quad \delta_\e=\e^{\a-1}.$
	
	Throughout this section, $E^\e$ denotes the scalar extension operator
	and $\mathcal E^\e$ denotes the divergence-preserving vector extension
	operator from Lemmas~\ref{Emain}--\ref{LemD02}.
	
	\subsection{Periodic unfolding}
	
	For every $x\in\mathbb{R}^3$, we use the decomposition
	$$
	\frac{x}{\e}
	=
	\left[\frac{x}{\e}\right]
	+
	\left\{\frac{x}{\e}\right\},
	\qquad
	\left[\frac{x}{\e}\right]\in\mathbb{Z}^3,
	\qquad
	\left\{\frac{x}{\e}\right\}\in Y,
	$$
	where $\left[x/\e\right]$ denotes the unique integer vector such that
	$\left\{x/\e\right\}\in Y$.
	
	Let $\widehat\O^\e$ and $\Lambda^\e$ be the union of the complete cells and
	the corresponding boundary layer introduced in Subsection~\ref{SSec31}.
	The periodic unfolding operator
	$$
	\Te:
	L^2(S\times\O)
	\longrightarrow
	L^2(\O_T\times Y)
	$$
	is defined by
	$$
	\Te(w)(t,x,y)
	:=
	\begin{cases}
		\displaystyle
		w\left(
		t,
		\e\left[\frac{x}{\e}\right]
		+
		\e y
		\right),
		&
		(t,x,y)\in S\times\widehat\O^\e\times Y,
		\\[2mm]
		0,
		&
		(t,x,y)\in S\times\Lambda^\e\times Y.
	\end{cases}
	$$
	Thus, for fixed $y\in Y$, the function $\Te(w)$ is constant with respect to
	$x$ on each complete cell $\e(\kappa+Y)$, while the microscopic oscillations
	are represented by the variable $y$. The operator is set equal to zero on the
	boundary layer $\Lambda^\e$, whose measure tends to zero as $\e\to0$.
	
	\begin{definition}[Perforated unfolding operator]
		\label{DefPerforatedUnfolding}
		Let
		$
		v^\e\in L^2(S;H^1(\O_p^\e)),
		\quad
		\bv^\e\in L^2(S;\GH^1_{\div}(\O_p^\e)).
		$
		Using the extension operators introduced in section \ref{Sec04}, we define
		$$
		\Te^\ast(v^\e)
		:=
		\Te(E^\e v^\e),
		\qquad
		\Te^\ast(\bv^\e)
		:=
		\Te(\Ec^\e\bv^\e).
		$$
		For the corresponding gradients, we set
		$$
		\Te^\ast(\nabla v^\e)
		:=
		\Te(\nabla E^\e v^\e),
		\qquad
		\Te^\ast(\nabla\bv^\e)
		:=
		\Te(\nabla\Ec^\e\bv^\e).
		$$
	\end{definition}
	
	The values of the extensions inside the obstacles are auxiliary and must
	not contribute to integrals over the fluid domain. To remove these
	contributions, we introduce
	$$
	\chi_{p,\e}
	:=
	\chi_{Y_{p,\e}},
	\qquad
	Y_{p,\e}
	:=
	Y\setminus\overline{Y_{s,\e}},
	\qquad
	Y_{s,\e}
	:=
	y_0+\delta_\e T.
	$$
	Since
	$$
	|Y_{s,\e}|
	=
	\delta_\e^3|T|
	\longrightarrow0,
	$$
	we have
	\begin{equation}\label{CellCharacteristicConvergence}
		\chi_{p,\e}
		\longrightarrow
		1
		\quad\text{strongly in }L^r(Y)
		\quad
		\text{for every }1\leq r<\infty.
	\end{equation}
	
	\begin{proposition}[Properties of the perforated unfolding operator]
		\label{PropPerforatedUnfolding}
		The following properties hold.
		
		\begin{enumerate}[label=\textnormal{(\roman*)}]
			\item For every
			$v^\e\in L^2(S;H^1(\O_p^\e))$, there exists a constant
			$C>0$, independent of $\e$, such that
			\begin{equation}\label{PerforatedUnfoldingEstimate}
				\|\Te^\ast(v^\e)\|_{L^2(\O_T\times Y)}
				+
				\|\Te^\ast(\nabla v^\e)\|_{L^2(\O_T\times Y)}
				\leq
				C
				\|v^\e\|_{L^2(S;H^1(\O_p^\e))}.
			\end{equation}
			The corresponding estimate holds componentwise for vector
			fields.
			
			\item For every
			$v^\e\in L^2(S;H^1(\O_p^\e))$,
			\begin{equation}\label{PerforatedGradientIdentity}
				\nabla_y\Te^\ast(v^\e)
				=
				\e\Te^\ast(\nabla v^\e)
				\quad
				\text{a.e. in }\O_T\times Y.
			\end{equation}
			The same identity holds componentwise for vector fields.
			
			\item Let
			$v^\e\in L^2(S;H^1(\O_p^\e))$, and let
			$H:\mathbb R\to\mathbb R$ be continuous with $H(0)=0$ and
			$
			H(v^\e)\in L^1(S\times\O_p^\e).
			$
			Then
			\begin{equation}\label{PerforatedCompositionIdentity}
				\int_{\O_T\times Y}
				\chi_{p,\e}(y)
				H\bigl(\Te^\ast(v^\e)\bigr)
				\,\md(y,x,t)
				=
				\int_{S\times\widehat\O_p^\e}
				H(v^\e)
				\,\md(x,t).
			\end{equation}
			
			\item Let $v^\e$ and $w^\e$ admit the extensions appearing in
			Definition~\ref{DefPerforatedUnfolding}, and assume that
			$
			v^\e w^\e
			\in
			L^1(S\times\O_p^\e).
			$
			Then
			\begin{equation}\label{PerforatedProductIdentity}
				\int_{\O_T\times Y}
				\chi_{p,\e}(y)
				\Te^\ast(v^\e)
				\Te^\ast(w^\e)
				\,\md(y,x,t)
				=
				\int_{S\times\widehat\O_p^\e}
				v^\e w^\e
				\,\md(x,t).
			\end{equation}
			
			\item If, in addition,
			$$
			\int_{S\times\Lambda^\e}
			|F(v^\e)|
			\,\md(x,t)
			\longrightarrow0,
			$$
			then
			\begin{equation}\label{PerforatedCompositionCriterion}
				\int_{\O_T\times Y}
				\chi_{p,\e}(y)
				F\bigl(\Te^\ast(v^\e)\bigr)
				\,\md(y,x,t)
				-
				\int_{S\times\O_p^\e}
				F(v^\e)
				\,\md(x,t)
				\longrightarrow0.
			\end{equation}
			Similarly, if
			$$
			\int_{S\times\Lambda^\e}
			|v^\e w^\e|
			\,\md(x,t)
			\longrightarrow0,
			$$
			then
			$$
			\int_{\O_T\times Y}
			\chi_{p,\e}(y)
			\Te^\ast(v^\e)
			\Te^\ast(w^\e)
			\,\md(y,x,t)
			-
			\int_{S\times\O_p^\e}
			v^\e w^\e
			\,\md(x,t)
			\longrightarrow0.
			$$
		\end{enumerate}
	\end{proposition}
	
	\begin{proof}
		The estimate \eqref{PerforatedUnfoldingEstimate} follows from the
		continuity of the classical unfolding operator and the uniform
		bounds for $E^\e$ and $\Ec^\e$. Moreover,
		$
		\nabla_y\Te(z)
		=
		\e\Te(\nabla z)
		$
		for every $z\in L^2(S;H^1(\O))$, which gives
		\eqref{PerforatedGradientIdentity}.
		
		For $x\in\widehat\O^\e$ and $y\in Y_{p,\e}$, the point
		$
		\e\left[\frac{x}{\e}\right]+\e y
		$
		belongs to $\widehat\O_p^\e$. Hence, the extensions coincide with
		the original functions at the corresponding unfolded points. Using
		the composition and product properties of $\Te$, together with its
		integration formula, we obtain
		\eqref{PerforatedCompositionIdentity} and
		\eqref{PerforatedProductIdentity}. The final assertions follow from
		$
		\O_p^\e
		=
		\widehat\O_p^\e
		\cup
		\Lambda^\e
		$
		up to a set of measure zero.
	\end{proof}

	\subsection{Compactness properties}\label{SSec62}

\begin{lemma}[Periodic-scale compactness]
	\label{lem:unfolding-compactness}
	 The following assertions hold, up to the
	extraction of subsequences which are not relabeled.
	
	\begin{enumerate}[label=\textnormal{(\roman*)}]
		\item
		Let $v^\e\in L^2(S;H^1(\O_p^\e))$
		satisfy
		\[
		\|v^\e\|_{L^2(S;H^1(\O_p^\e))}
		\leq C.
		\]
		Then there exist $v\in L^2(S;H^1(\O))$
		and $\widehat V
		\in
		L^2(\O_T;H^1_\#(Y))$
		such that
		\begin{equation}
			\label{PeriodicScalarCompactness}
			\begin{aligned}
				\Te^\ast(v^\e)
				&\rightharpoonup
				v
				&&\text{weakly in }L^2(\O_T\times Y),
				\\
				\Te^\ast(\nabla v^\e)
				&\rightharpoonup
				\nabla_xv+\nabla_y\widehat V
				&&\text{weakly in }
				L^2(\O_T\times Y)^3.
			\end{aligned}
		\end{equation}
		
		If, in addition,
		\[
		\|v^\e\|_{L^\infty(S;H^1(\O_p^\e))}
		+
		\|\partial_tv^\e\|_
		{L^2(S;H^{-1}(\O_p^\e))}
		\leq C,
		\]
		then
		\begin{equation}
			\label{PeriodicScalarStrongCompactness}
			\Te^\ast(v^\e)
			\longrightarrow
			v
			\quad\text{strongly in }
			L^2(S;L^s(\O\times Y))
			\quad
			\text{for every }2\leq s<6.
		\end{equation}
		
		\item
		Let $\bv^\e
		\in
		L^2(S;\GH^1_{\div}(\O_p^\e))$
		satisfy
		\begin{equation}
			\label{VelocityCompactnessAssumptions}
			\|\bv^\e\|_
			{L^\infty(S;L^2(\O_p^\e)^3)}
			+
			\|\bv^\e\|_
			{L^2(S;H^1(\O_p^\e)^3)}
			+
			\|\partial_t\bv^\e\|_
			{L^{4/3}(S;\GH^{-1}_{\div}(\O_p^\e))}
			\leq C.
		\end{equation}
		Then there exist $\bv
		\in
		L^\infty(S;L^2(\O)^3)
		\cap
		L^2(S;H^1_{0,\div}(\O)^3)$
		and $\widehat\bv
		\in
		L^2(\O_T;H^1_\#(Y)^3)$,
		$\nabla_y\cdot\widehat\bv=0$
		such that
		\begin{equation}
			\label{PeriodicVelocityCompactness}
			\begin{aligned}
				\Te^\ast(\bv^\e)
				&\longrightarrow
				\bv
				&&\text{strongly in }
				L^2(\O_T\times Y)^3,
				\\
				\Te^\ast(\nabla\bv^\e)
				&\rightharpoonup
				\nabla_x\bv+\nabla_y\widehat\bv
				&&\text{weakly in }
				L^2(\O_T\times Y)^{3\times3}.
			\end{aligned}
		\end{equation}
	\end{enumerate}
\end{lemma}

\begin{proof}
	\textbf{Step 1: Scalar compactness.}
	Set $\widetilde v^\e:=E^\e v^\e$.
	By Lemma~\ref{Emain}, the sequence
	$\{\widetilde v^\e\}$ is bounded in
	$L^2(S;H^1(\O))$. The fixed-domain unfolding compactness theorem
	therefore gives, after extraction,
	\[
	\widetilde v^\e
	\rightharpoonup v
	\quad\text{weakly in }L^2(S;H^1(\O)),
	\]
	and
	\[
	\begin{aligned}
		\Te(\widetilde v^\e)
		&\rightharpoonup v
		&&\text{weakly in }L^2(\O_T\times Y),\\
		\Te(\nabla\widetilde v^\e)
		&\rightharpoonup
		\nabla_xv+\nabla_y\widehat V
		&&\text{weakly in }
		L^2(\O_T\times Y)^3
	\end{aligned}
	\]
	for some $\widehat V
	\in L^2(\O_T;H^1_\#(Y))$.
	Here and below, the limit functions on $\O_T$ are identified with
	functions on $\O_T\times Y$ that are independent of $y$.
	Since
	\[
	\Te^\ast(v^\e)=\Te(\widetilde v^\e),
	\qquad
	\Te^\ast(\nabla v^\e)
	=
	\Te(\nabla\widetilde v^\e),
	\]
	this proves \eqref{PeriodicScalarCompactness}.
	
	Under the additional time-regularity assumptions,
	Lemma~\ref{lem:scalar-two-extensions} gives, after a further
	extraction,
	\[
	\widetilde v^\e
	\longrightarrow v
	\quad\text{strongly in }
	C([0,T];L^2(\O)).
	\]
	The continuity and consistency of the unfolding operator then give
	\[
	\Te^\ast(v^\e)
	=
	\Te(\widetilde v^\e)
	\longrightarrow v
	\quad\text{strongly in }L^2(\O_T\times Y).
	\]
	Moreover, $\|\Te^\ast(v^\e)\|_
	{L^\infty(S;L^6(\O\times Y))}
	\leq C$.
	Interpolation between the strong $L^2$ convergence and this
	uniform $L^\infty(S;L^6)$ bound yields
	\eqref{PeriodicScalarStrongCompactness} for every $2\leq s<6$.
	
	\textbf{Step 2: Weak velocity compactness.}
	Let $\widetilde\bv^\e:=\mathcal E^\e\bv^\e$
	be the primal extension and let $\overline\bv^\e$ be the zero
	extension. By Lemma~\ref{LemD02} and
	\eqref{VelocityCompactnessAssumptions},
	\begin{equation}\label{VelocityPrimalAndZeroBounds}
		\|\widetilde\bv^\e\|_
		{L^2(S;H^1_{0,\div}(\O)^3)}
		+
		\|\overline\bv^\e\|_
		{L^\infty(S;L^2(\O)^3)}
		\leq C,
	\end{equation}
	and
	\begin{equation}\label{VelocityPrimalZeroCloseness}
		\widetilde\bv^\e-\overline\bv^\e
		\longrightarrow0
		\quad\text{strongly in }L^2(\O_T)^3.
	\end{equation}
	Consequently, there exist $\bv\in L^2(S;H^1_{0,\div}(\O)^3)$ and $\widehat\bv
	\in L^2(\O_T;H^1_\#(Y)^3)$
	such that
	\begin{equation}\label{VelocityWeakFixedDomainConvergences}
		\widetilde\bv^\e
		\rightharpoonup\bv
		\quad\text{weakly in }
		L^2(S;H^1_0(\O)^3),\quad 
		\Te^\ast(\nabla\bv^\e)
		\rightharpoonup
		\nabla_x\bv+\nabla_y\widehat\bv
		\quad\text{weakly in }
		L^2(\O_T\times Y)^{3\times3}.
	\end{equation}
	
	\textbf{Step 3: Compactness of the velocity pairings.}
	Fix $\zeta\in C^\infty_{c,\div}(\O)^3$, and let $\Rc^\e\zeta\in\GH^1_{\div}(\O_p^\e)$
	be the field constructed in
	Lemma~\ref{lem:solenoidal-restriction}. Define
	\[
	A_\e^\zeta(t)
	:=
	(\bv^\e(t),\Rc^\e\zeta)_{L^2(\O_p^\e)}.
	\]
	Since $\Rc^\e\zeta$ is time-independent,
	\[
	\frac{\md}{\md t}A_\e^\zeta
	=
	\left\langle
	\partial_t\bv^\e,\Rc^\e\zeta
	\right\rangle_
	{\GH^{-1}_{\div}(\O_p^\e),
		\GH^1_{\div}(\O_p^\e)},\quad\text{in $\mathcal D'(S)$}.
	\]
	It follows from
	\eqref{VelocityCompactnessAssumptions} and
	\eqref{SolenoidalRestrictionBound} that
	\[
	\|A_\e^\zeta\|_{W^{1,4/3}(S)}
	\leq C_\zeta.
	\]
	Hence $\{A_\e^\zeta\}$ is relatively compact in
	$C([0,T])$. 
	
	Next, define
	\[
	B_\e^\zeta(t)
	:=
	(\overline\bv^\e(t),\zeta)_{L^2(\O)}
	=
	(\bv^\e(t),\zeta)_{L^2(\O_p^\e)}.
	\]
	Then
	\[
	\begin{aligned}
		\|B_\e^\zeta-A_\e^\zeta\|_{L^\infty(S)}
		&\leq
		\|\bv^\e\|_
		{L^\infty(S;L^2(\O_p^\e))}
		\|\zeta-\overline{\Rc^\e\zeta}\|_{L^2(\O)}
		\longrightarrow0
	\end{aligned}
	\]
	by \eqref{SolenoidalRestrictionConvergence}. Thus
	$\{B_\e^\zeta\}$ is relatively compact in $L^\infty(S)$, and
	every limit has a continuous representative.
	
	Furthermore,
	\[
	\begin{aligned}
		\bigl\|
		(\widetilde\bv^\e-\overline\bv^\e,\zeta)_{L^2(\O)}
		\bigr\|_{L^2(S)}
		&\leq
		\|\zeta\|_{L^2(\O)}
		\|\widetilde\bv^\e-\overline\bv^\e\|_{L^2(\O_T)}
		\longrightarrow0.
	\end{aligned}
	\]
	Therefore, for every
	$\zeta\in C^\infty_{c,\div}(\O)^3$, the sequence
	\[
	t\longmapsto
	(\widetilde\bv^\e(t),\zeta)_{L^2(\O)}
	\]
	is relatively compact in $L^2(S)$.
%
	
	\textbf{Step 4: Strong velocity compactness.}
	Let $\{\zeta_k\}_{k\geq1}
	\subset C^\infty_{c,\div}(\O)^3$
	be an orthonormal basis of the solenoidal pivot space $H^\sigma$
	 and let $P_N$ be the orthogonal projection onto $\operatorname{span}\{\zeta_1,\ldots,\zeta_N\}$.
	Using a diagonal extraction, the compactness of the scalar
	pairings obtained in Step~3 implies that, for every fixed $N$,
	\[
	\{P_N\widetilde\bv^\e\}
	\quad\text{is relatively compact in }L^2(\O_T)^3.
	\]
	
	The compact embedding $H^1_{0,\div}(\O)^3
	\Subset H^\sigma$
	and the strong convergence $P_Nz\to z$ in $H^\sigma$ imply
	\[
	\eta_N
	:=
	\sup_{\substack{
			z\in H^1_{0,\div}(\O)^3\\
			\|z\|_{H^1(\O)}\leq1}}
	\|(I-P_N)z\|_{L^2(\O)}
	\longrightarrow0.
	\]
	Consequently, by \eqref{VelocityPrimalAndZeroBounds},
	\[
	\sup_\e
	\|(I-P_N)\widetilde\bv^\e\|_{L^2(\O_T)}
	\leq C\eta_N
	\longrightarrow0.
	\]
	The finite-dimensional compactness and the uniform tail estimate
	show that $\{\widetilde\bv^\e\}$ is relatively compact in
	$L^2(\O_T)^3$. In view of
	\eqref{VelocityWeakFixedDomainConvergences},
	\begin{equation}\label{VelocityStrongFixedDomainConvergence}
		\widetilde\bv^\e
		\longrightarrow\bv
		\quad\text{strongly in }L^2(\O_T)^3.
	\end{equation}
	By \eqref{VelocityPrimalZeroCloseness},
	\[
	\overline\bv^\e
	\longrightarrow\bv
	\quad\text{strongly in }L^2(\O_T)^3.
	\]
	The uniform
	$L^\infty(S;L^2(\O)^3)$ bound for
	$\overline\bv^\e$, together with weak-star compactness, gives $\bv\in L^\infty(S;L^2(\O)^3)$.
	For every $\zeta\in C^\infty_{c,\div}(\O)^3$, the estimates of
	Step~3 and \eqref{VelocityStrongFixedDomainConvergence} give
	\[
	A_\e^\zeta
	\longrightarrow
	(\bv,\zeta)_{L^2(\O)}
	\quad\text{strongly in }L^2(S).
	\]
	Since $\{A_\e^\zeta\}$ is relatively compact in $C([0,T])$, we
	conclude that
	\begin{equation}\label{VelocityPairingConvergence}
		A_\e^\zeta
		\longrightarrow
		(\bv,\zeta)_{L^2(\O)}
		\quad\text{strongly in }C([0,T]).
	\end{equation}
	In particular, $\bv$ admits a representative in
	$C_w([0,T];H^\sigma)$.

	\textbf{Step 5: Unfolded velocity limits.}	
	Since $\Te^\ast(\bv^\e)=\Te(\widetilde\bv^\e)$,
	\eqref{VelocityStrongFixedDomainConvergence} and the continuity and
	consistency of $\Te$ give
	\[
	\Te^\ast(\bv^\e)
	\longrightarrow\bv
	\quad\text{strongly in }L^2(\O_T\times Y)^3.
	\]
	Together with
	\eqref{VelocityWeakFixedDomainConvergences}$_2$, this proves
	\eqref{PeriodicVelocityCompactness}.
	
	Finally, since
	$\div\,\widetilde\bv^\e=0$ in $\O$, the weak convergence in
	\eqref{VelocityWeakFixedDomainConvergences} gives $\nabla_x\cdot\bv=0$.
	Taking the trace in
	\eqref{VelocityWeakFixedDomainConvergences}$_2$ yields
	\[
	\nabla_x\cdot\bv+\nabla_y\cdot\widehat\bv=0\implies
	\nabla_y\cdot\widehat\bv=0
	\quad\text{in }\mathcal D'(\O_T\times Y).
	\]
	This completes the proof.
\end{proof}

\begin{lemma}[Compactness and unfolded convergence]
	\label{lem:limit-compactness}
	Let $(\Gu^\e,\phi^\e,\mu^\e)
	\in
	\mathfrak U^\e
	\times
	\mathfrak C^\e
	\times
	\mathfrak W^\e$
	be a weak solution of \eqref{MainP01} in the sense of
	Definition~\ref{Def01}. Then, up to a subsequence, there exist
	\[
	(\Gu,\phi,\mu)\in\fU\X\fC\times\fW,\quad\text{and}\quad
	\widehat\Gu
	\in
	L^2(\O_T;H^1_\#(Y)^3),
	\qquad
	\widehat\phi,\widehat\mu
	\in
	L^2(\O_T;H^1_\#(Y)),
	\]
	such that the following assertions hold.
	
	\smallskip
	\noindent\textnormal{(i) Periodic-scale convergences.}
	For every $2\leq s<6$,
	\begin{equation}
		\label{TSL0}
		\begin{aligned}
			\frac{1}{\sqrt{\lambda^\e}}
			\Te^\ast(\Gu^\e)
			&\longrightarrow
			\Gu
			&&\text{strongly in }
			L^2(S;L^s(\O\times Y)^3),
			\\
			\Te^\ast(\phi^\e)
			&\longrightarrow
			\phi
			&&\text{strongly in }
			C([0,T];L^s(\O\times Y)),
			\\
			\Te^\ast(\mu^\e)
			&\rightharpoonup
			\mu
			&&\text{weakly in }
			L^2(\O_T\times Y).
		\end{aligned}
	\end{equation}
	Moreover,
	\begin{equation}
		\label{TSL01}
		\begin{aligned}
			\frac{1}{\sqrt{\lambda^\e}}
			\Te^\ast(\nabla\Gu^\e)
			&\rightharpoonup
			\nabla_x\Gu+\nabla_y\widehat\Gu
			&&\text{weakly in }
			L^2(\O_T\times Y)^{3\times3},
			\\
			\Te^\ast(\nabla\phi^\e)
			&\rightharpoonup
			\nabla_x\phi+\nabla_y\widehat\phi
			&&\text{weakly in }
			L^2(\O_T\times Y)^3,
			\\
			\Te^\ast(\nabla\mu^\e)
			&\rightharpoonup
			\nabla_x\mu+\nabla_y\widehat\mu
			&&\text{weakly in }
			L^2(\O_T\times Y)^3.
		\end{aligned}
	\end{equation}
	Furthermore,
	\begin{equation}
		\label{eq:limit-divergence}
		\nabla_x\cdot\Gu=0
		\quad\text{in }\Dc'(\O_T),
		\qquad
		\nabla_y\cdot\widehat\Gu=0
		\quad\text{in }\Dc'(\O_T\times Y).
	\end{equation}
	
	\smallskip
	\noindent\textnormal{(ii) Convergence of the nonlinear terms.}
	For every $1\leq r<3/2$,
	\begin{equation}
		\label{TSL02}
		\begin{aligned}
			\chi_{p,\e}\,
			\Te^\ast(\phi^\e)
			\frac{\Te^\ast(\Gu^\e)}
			{\sqrt{\lambda^\e}}
			&\longrightarrow
			\phi\Gu
			&&\text{strongly in }
			L^2(\O_T\times Y)^3,
			\\
			\chi_{p,\e}\,
			F\bigl(\Te^\ast(\phi^\e)\bigr)
			&\longrightarrow
			F(\phi)
			&&\text{strongly in }
			L^1(\O_T\times Y),
			\\
			\chi_{p,\e}\,
			F'\bigl(\Te^\ast(\phi^\e)\bigr)
			&\rightharpoonup
			F'(\phi)
			&&\text{weakly in }
			L^2(\O_T\times Y),
			\\
			\chi_{p,\e}\,
			G\bigl(\Te^\ast(\phi^\e)\bigr)
			&\longrightarrow
			G(\phi)
			&&\text{strongly in }
			L^2(\O_T\times Y),
			\\
			\chi_{p,\e}\,
			\Te^\ast(\phi^\e)
			\Te^\ast(\nabla\mu^\e)
			&\rightharpoonup
			\phi
			\bigl(
			\nabla_x\mu+\nabla_y\widehat\mu
			\bigr)
			&&\text{weakly in }
			L^2(S;L^r(\O\times Y)^3),
			\\
			\chi_{p,\e}\,
			\frac{\Te^\ast(\Gu^\e)}
			{\sqrt{\lambda^\e}}
			\otimes
			\frac{\Te^\ast(\Gu^\e)}
			{\sqrt{\lambda^\e}}
			&\longrightarrow
			\Gu\otimes\Gu
			&&\text{strongly in }
			L^1(\O_T\times Y)^{3\times3}.
		\end{aligned}
	\end{equation}
\end{lemma}
\begin{proof}
	\textbf{Step 1: Periodic-scale convergences.}
	Set $\bv^\e:=\frac{\Gu^\e}{\sqrt{\lambda^\e}}$.
	The estimates in \eqref{MainE01} give
	\[
	\|\bv^\e\|_{L^\infty(S;L^2(\O_p^\e)^3)}
	+
	\|\bv^\e\|_{L^2(S;H^1(\O_p^\e)^3)}
	+
	\|\partial_t\bv^\e\|_
	{L^{4/3}(S;\GH^{-1}_{\div}(\O_p^\e))}
	\leq C.
	\]
	Lemma~\ref{lem:unfolding-compactness}(ii) therefore gives $\Gu\in\fU$
	and $\widehat\Gu
	\in L^2(\O_T;H^1_\#(Y)^3)$, $\nabla_y\cdot\widehat\Gu=0$,
	together with
	\[
	\frac{1}{\sqrt{\lambda^\e}}
	\Te^\ast(\Gu^\e)
	\longrightarrow
	\Gu
	\quad\text{strongly in }L^2(\O_T\times Y)^3
	\]
	and the velocity-gradient convergence in \eqref{TSL01}.
	
	The normalized energy estimate and the uniform extension theorem
	yield
	\[
	\left\|
	\frac{1}{\sqrt{\lambda^\e}}
	\Te^\ast(\Gu^\e)
	\right\|_{L^2(S;L^6(\O\times Y)^3)}
	\leq C.
	\]
	Interpolating this bound with the preceding strong $L^2$
	convergence gives \eqref{TSL0}$_1$
	for every $2\leq s<6$.
	
	Applying Lemma~\ref{lem:unfolding-compactness}(i) to
	$\phi^\e$ and $\mu^\e$ gives \eqref{TSL0}$_{3}$ and \eqref{TSL01}$_{2,3}$. 
	
	Set $\widetilde\phi^\e:=E^\e\phi^\e$ satisfy,
	after extraction,
	\[
	\widetilde\phi^\e
	\longrightarrow\phi
	\quad\text{strongly in }
	C([0,T];L^2(\O)),
	\]
	using 	Lemma~\ref{lem:scalar-two-extensions} and \eqref{MainE01}. Moreover, $\{\widetilde\phi^\e\}$ is bounded in
	$L^\infty(S;H^1(\O))$ and interpolation yields
	\[
	\widetilde\phi^\e
	\longrightarrow\phi
	\quad\text{strongly in }
	C([0,T];L^s(\O))
	\]
	for every $2\leq s<6$. The uniform consistency of the unfolding operator on compact subsets
	of $L^s(\O)$ gives \eqref{TSL0}$_2$; see
	\cite[Sections~1.4 and~2.2]{CDG}.

	\textbf{Step 2: Nonlinear terms.}
	For brevity, set
	\[
	W^\e
	:=
	\frac{\Te^\ast(\Gu^\e)}{\sqrt{\lambda^\e}},
	\qquad
	\Phi^\e:=\Te^\ast(\phi^\e).
	\]
	Since $\chi_{p,\e}\longrightarrow1$ strongly in $L^q(Y)$
	for every finite $q$, the convergence of $\Phi^\e$ gives
	\begin{equation}\label{CharacteristicStrongProducts}
		\chi_{p,\e}\Phi^\e
		\longrightarrow\phi
		\quad\text{strongly in }
		C([0,T];L^s(\O\times Y))
	\end{equation}
	for every $2\leq s<6$.
	
	Taking $s=4$ in
	\eqref{CharacteristicStrongProducts} and in the velocity
	convergence yields \eqref{TSL02}$_1$.
	
	The polynomial growth of $F$ and $F'$ and the strong convergence
	of $\Phi^\e$ in
	$C([0,T];L^s(\O\times Y))$, $s<6$, give \eqref{TSL02}$_2$.
	Moreover, $\{\chi_{p,\e}F'(\Phi^\e)\}$ is bounded in $L^2(\O_T\times Y)$,
	and
	\[
	\chi_{p,\e}F'(\Phi^\e)
	\longrightarrow F'(\phi)
	\quad\text{strongly in }L^q(\O_T\times Y)
	\]
	for every $1\leq q<2$. Hence, we obtain \eqref{TSL02}$_3$.
	The Lipschitz continuity of $G$ similarly gives \eqref{TSL02}$_4$.
	
	Let $1\leq r<3/2$ and choose $2\leq s<6$ such that
	\[
	\frac1r=\frac12+\frac1s.
	\]
	Combining \eqref{CharacteristicStrongProducts} with the weak
	convergence of $\Te^\ast(\nabla\mu^\e)$ gives \eqref{TSL02}$_5$.
	
	Finally, the strong convergence
	$W^\e\to\Gu$ in $L^2(\O_T\times Y)^3$ implies
	\[
	W^\e\otimes W^\e
	\longrightarrow
	\Gu\otimes\Gu
	\quad\text{strongly in }
	L^1(\O_T\times Y)^{3\times3}.
	\]
	Since $\chi_{p,\e}\to1$ strongly in every finite
	$L^q(Y)$, multiplication by $\chi_{p,\e}$ does not change the
	limit. This proves \eqref{TSL02}$_6$.
\end{proof}

\begin{lemma}[Identification of the initial conditions]
	\label{lem:limit-initial-data}
	Under assumption \eqref{ICA01}, the limit fields obtained in
	Lemma~\ref{lem:limit-compactness} satisfy
	\begin{equation}\label{eq:limit-initial-conditions}
		\Gu(0)=\Gu_0
		\quad\text{in }H^{-1}_{0,\div}(\O)^3,
		\qquad
		\phi(0)=\phi_0
		\quad\text{in }L^2(\O).
	\end{equation}
\end{lemma}

\begin{proof}
	Set $\bv^\e:=\frac{\Gu^\e}{\sqrt{\lambda^\e}}$
	and fix $\zeta\in C^\infty_{c,\div}(\O)^3$. Let
	$\Rc^\e\zeta\in\GH^1_{\div}(\O_p^\e)$ be the solenoidal restriction
	constructed in Lemma~\ref{lem:solenoidal-restriction}, and define
	\[
	A_\e^\zeta(t)
	:=
	\bigl(\bv^\e(t),\Rc^\e\zeta\bigr)_{L^2(\O_p^\e)}.
	\]
	By \eqref{VelocityPairingConvergence},
	\[
	A_\e^\zeta
	\longrightarrow
	(\Gu,\zeta)_{L^2(\O)}
	\quad\text{strongly in }C([0,T]).
	\]
	In particular, $\Gu$ admits a representative in
	$C_w([0,T];H^\sigma)$.
	
	Let $\overline{\Rc^\e\zeta}$ denote the zero extension of
	$\Rc^\e\zeta$. Since $\Gu^\e(0)=\Gu_0^\e$, we have
	\[
	A_\e^\zeta(0)
	=
	\left(
	\frac{\overline{\Gu_0^\e}}{\sqrt{\lambda^\e}},
	\overline{\Rc^\e\zeta}
	\right)_{L^2(\O)}.
	\]
	Therefore,
	\[
	\begin{aligned}
		\left|
		A_\e^\zeta(0)-(\Gu_0,\zeta)_{L^2(\O)}
		\right|
		&\leq
		\left\|
		\frac{\overline{\Gu_0^\e}}{\sqrt{\lambda^\e}}-\Gu_0
		\right\|_{L^2(\O)}
		\|\overline{\Rc^\e\zeta}\|_{L^2(\O)}
		+
		\|\Gu_0\|_{L^2(\O)}
		\|\overline{\Rc^\e\zeta}-\zeta\|_{L^2(\O)}.
	\end{aligned}
	\]
	The first term tends to zero by \eqref{ICA01}, while the second tends
	to zero by
	\eqref{SolenoidalRestrictionConvergence}. Hence
	\[
	A_\e^\zeta(0)
	\longrightarrow
	(\Gu_0,\zeta)_{L^2(\O)}.
	\]
	Evaluating the convergence of $A_\e^\zeta$ at $t=0$ gives
	\[
	(\Gu(0),\zeta)_{L^2(\O)}
	=
	(\Gu_0,\zeta)_{L^2(\O)}
	\qquad
	\forall\,\zeta\in C^\infty_{c,\div}(\O)^3.
	\]
	By density, $\Gu(0)=\Gu_0$ in
	$H^{-1}_{0,\div}(\O)^3$.
	
	For the phase field, let $\overline\phi^\e$ be the zero extension of
	$\phi^\e$. Lemmas~\ref{lem:scalar-two-extensions} and
	\ref{lem:limit-compactness} give
	\[
	\overline\phi^\e
	\longrightarrow
	\phi
	\quad\text{strongly in }C([0,T];L^2(\O)).
	\]
	Since $\phi^\e(0)=\phi_0|_{\O_p^\e}$,
	\[
	\|\overline\phi^\e(0)-\phi_0\|_{L^2(\O)}
	=
	\|\phi_0\|_{L^2(\O_s^\e)}
	\longrightarrow0,
	\]
	because $\phi_0\in L^2(\O)$ and $|\O_s^\e|\to0$. Evaluation at
	$t=0$ proves $\phi(0)=\phi_0$ in $L^2(\O)$.
\end{proof}

	\section{Homogenization of the NSCH system}
	\label{Sec06}
	
	In this section we derive the homogenized equations associated with the
	microscopic problem \eqref{MainP01}. The derivation relies on the
	uniform a priori estimates established in
	Section~\ref{SSec51}, the extension operators constructed in
	Section~\ref{Sec04}, and the periodic unfolding framework introduced in
	Section~\ref{SSec61}.
	
	We first identify the two-scale limits of the microscopic solutions and
	their gradients. These compactness results allow us to pass to the limit
	in each term of the weak formulation and identify the effective
	macroscopic equations.
	
	Throughout this section we distinguish the homogenized models according
	to the scaling parameter $\lambda$ introduced in
	\eqref{Lam01}. When $\lambda=0$, the inertial effects disappear and the
	homogenized system is of Stokes--Cahn--Hilliard type. When
	$\lambda>0$, the limiting equations retain the nonlinear convection term
	and lead to a Navier--Stokes--Cahn--Hilliard system.
	
	\subsection{Two-scale limit problem in the subcritical regime
		$\alpha>3$}
	\label{SSecLimitSpaces}
	
	We begin by introducing the periodic corrector spaces associated with
	the microscopic oscillations.
	
	The velocity correctors belong to
	$$
	\Vc_\#
	:=
	\left\{
	\bv\in H^1_{\mathrm{per}}(Y)^3:
	\nabla_y\cdot\bv=0,
	\quad
	\int_Y\bv(y)\,\md y=0
	\right\},
	$$
	whereas the scalar correctors belong to
	$$
	\mathcal H_\#
	:=
	\left\{
	v\in H^1_{\mathrm{per}}(Y):
	\int_Yv(y)\,\md y=0
	\right\}.
	$$
	For
	$
	a_\e=\e^\alpha$, $\alpha>3$,
	the rescaled obstacles satisfy
	$
	Y_{s,\e}
	=
	y_0+\frac{a_\e}{\e}T$, $|Y_{s,\e}|
	\longrightarrow0$.
	Consequently, the microscopic obstacles disappear at the periodic scale,
	and the corrector problems are posed on the whole reference cell $Y$.
	No boundary condition survives on the limiting point $y_0$.
	
	By Lemma~\ref{lem:limit-compactness}, the two-scale limit fields satisfy
	\begin{equation}
		\label{PeriodicLimitFieldSpace}
		\Gu
		\in
		\fU,
		\widehat\Gu
		\in
		L^2(\O_T;\Vc_\#),\quad
		\phi
		\in
		\fC,
		\widehat\phi
		\in
		L^2(\O_T;\mathcal H_\#),\quad
		\mu
		\in
		\fW,
		\widehat\mu
		\in
		L^2(\O_T;\mathcal H_\#).
	\end{equation}
	Moreover,
	$$
	\nabla_x\cdot\Gu=0
	\quad\text{in }\Dc'(\O_T),
	\qquad
	\nabla_y\cdot\widehat\Gu=0
	\quad\text{in }\Dc'(\O_T\times Y).
	$$
	
	Accordingly, the collection of admissible limit unknowns is
	$$
	\mathscr X_{>3}
	:=
	\fU
	\times
	L^2(\O_T;\Vc_\#)
	\times
	\fC
	\times
	L^2(\O_T;\mathcal H_\#)
	\times
	\fW
	\times
	L^2(\O_T;\mathcal H_\#).
	$$
	Hence,
	$$
	(\Gu,\widehat\Gu,\phi,\widehat\phi,\mu,\widehat\mu)
	\in
	\mathscr X_{>3}.
	$$
	To derive the effective equations we employ the test spaces
	\begin{equation}
		\label{LimitTestSpaces}
		\begin{aligned}
			\GL_1
			&:=
			L^4(S;H^1_{0,\div}(\O)^3)
			\times
			L^2(\O_T;\Vc_\#),
			\\
			\GL_2
			&:=
			L^2(S;H^1(\O))
			\times
			L^2(\O_T;\mathcal H_\#),\quad
			\GL_3
			:=
			L^2(S;H^1(\O))
			\times
			L^2(\O_T;\mathcal H_\#).
		\end{aligned}
	\end{equation}
	\begin{theorem}[Two-scale limit system for $\alpha>3$]
		\label{Th02}
		Let $\alpha>3$, assume \eqref{Lam01}, and let
		$
		(\Gu,\widehat\Gu,\phi,\widehat\phi,\mu,\widehat\mu)
		\in
		\mathscr X_{>3}
		$
		be the limit fields obtained in
		Lemma~\ref{lem:limit-compactness}. 
		Moreover,
		\[
		\partial_t\Gu
		\in
		L^{4/3}(S;H^{-1}_{0,\div}(\O)^3),
		\qquad
		\partial_t\phi
		\in
		L^2(S;H^{-1}(\O)).
		\]
		For
		$(\bv,\widehat\bv)$ and $(v,\widehat v)$, set
		$$
		\D(\bv,\widehat\bv)
		:=
		D_x\bv+D_y\widehat\bv,
		\qquad
		\G(v,\widehat v)
		:=
		\nabla_xv+\nabla_y\widehat v.
		$$
		Then, for every
		$
		(\varphi_i,\widehat\varphi_i)\in\GL_i$, $i\in\{1,2,3\}$,
		the following identities hold:
		\begin{equation}\label{TSLM01}
			\begin{aligned}
				&\int_S
				\langle\partial_t\Gu,\varphi_1\rangle\,\md t-
				\sqrt{\lambda}
				\int_{\O_T}
				(\Gu\otimes\Gu):\nabla_x\varphi_1
				\,\md(x,t)+
				\int_{\O_T\times Y}
				\A(t,y)
				\D(\Gu,\widehat\Gu):
				\D(\varphi_1,\widehat\varphi_1)
				\,\md(y,x,t)
				+
				\sqrt{\lambda}
				\int_{\O_T}
				\phi\nabla_x\mu\cdot\varphi_1
				\,\md(x,t)
				\\
				&\qquad \hspace{12cm}=
				\int_{\O_T}
				\Gg\cdot\varphi_1
				\,\md(x,t),
				\\[2mm]
				&\int_S
				\langle\partial_t\phi,\varphi_2\rangle\,\md t
				-
				\sqrt{\lambda}
				\int_{\O_T}
				\phi\Gu\cdot\nabla_x\varphi_2
				\,\md(x,t)
				+
				\int_{\O_T\times Y}
				\G(\mu,\widehat\mu)
				\cdot
				\G(\varphi_2,\widehat\varphi_2)
				\,\md(y,x,t)
				+
				\int_{\O_T}
				G(\phi)\varphi_2
				\,\md(x,t)
				=0,
				\\[2mm]
				&\int_{\O_T}
				\mu\varphi_3
				\,\md(x,t)
				-
				\int_{\O_T\times Y}
				\G(\phi,\widehat\phi)
				\cdot
				\G(\varphi_3,\widehat\varphi_3)
				\,\md(y,x,t)
				-
				\int_{\O_T}
				F'(\phi)\varphi_3
				\,\md(x,t)
				=0.
			\end{aligned}
		\end{equation}
	\end{theorem}
	
	\begin{proof}
		We first establish \eqref{TSLM01} for smooth test functions. Set
		$$
		\GT_1
		:=
		C_c^\infty
		\bigl(
		S;C^\infty_{c,\div}(\O)^3
		\bigr)
		\times
		C_c^\infty
		\bigl(
		\O_T;C^\infty_{\#,\div}(Y)^3
		\bigr),
		$$
		and, for $i\in\{2,3\}$,
		$$
		\GT_i
		:=
		C_c^\infty
		\bigl(
		S;C^\infty(\overline\O)
		\bigr)
		\times
		C_c^\infty
		\bigl(
		\O_T;C^\infty_\#(Y)
		\bigr).
		$$
		
		Let
		$(\varphi_1,\widehat\varphi_1)\in\GT_1$. By
		Lemma~\ref{lem:test-functions-subcritical}, there exists an admissible
		solenoidal test function $\varphi_1^\e$ such that, with
		$$
		\Psi_1^\e
		:=
		\frac{\varphi_1^\e}{\sqrt{\lambda^\e}},
		$$
		we have
		\begin{equation}\label{eq:velocity-test-convergence}
			\begin{aligned}
				\Te^\ast(\Psi_1^\e)
				&\longrightarrow
				\varphi_1
				&&\text{strongly in }
				L^2(\O_T\times Y)^3,
				\\
				\Te^\ast(\nabla\Psi_1^\e)
				&\longrightarrow
				\nabla_x\varphi_1+\nabla_y\widehat\varphi_1
				&&\text{strongly in }
				L^2(\O_T\times Y)^{3\times3},
				\\
				\Te^\ast(\partial_t\Psi_1^\e)
				&\longrightarrow
				\partial_t\varphi_1
				&&\text{strongly in }
				L^2(\O_T\times Y)^3.
			\end{aligned}
		\end{equation}
		Consequently,
		$$
		\Te^\ast(D\Psi_1^\e)
		\longrightarrow
		\D(\varphi_1,\widehat\varphi_1)
		\quad\text{strongly in }
		L^2(\O_T\times Y)^{3\times3}.
		$$
		
		For
		$(\varphi_i,\widehat\varphi_i)\in\GT_i$,
		$i\in\{2,3\}$, define
		\begin{equation}\label{Test01}
			\varphi_i^\e(t,x)
			:=
			\varphi_i(t,x)
			+
			\e
			\widehat\varphi_i
			\left(
			t,x,\frac{x}{\e}
			\right).
		\end{equation}
		The restrictions of ${\varphi_i^\e}_{|\O^\e_p}$ to $\O_p^\e$ are admissible in
		\eqref{weak-phi1} and \eqref{weak-mu1}. Moreover, the standard
		consistency properties of periodic unfolding yield
		\begin{equation}\label{eq:scalar-test-convergence}
			\begin{aligned}
				\Te(\varphi_i^\e)
				&\longrightarrow
				\varphi_i
				&&\text{strongly in }L^2(\O_T\times Y),
				\\
				\Te(\partial_t\varphi_i^\e)
				&\longrightarrow
				\partial_t\varphi_i
				&&\text{strongly in }L^2(\O_T\times Y),
				\\
				\Te(\nabla\varphi_i^\e)
				&\longrightarrow
				\G(\varphi_i,\widehat\varphi_i)
				&&\text{strongly in }L^2(\O_T\times Y)^3;
			\end{aligned}
		\end{equation}

		We first pass to the limit in the phase and chemical-potential
		equations. Since
		$$
		\nabla\cdot\Gu^\e=0
		\quad\text{in }\O_p^\e,
		\qquad
		\Gu^\e\cdot\Gn^\e=0
		\quad\text{on }\partial\O_p^\e,
		$$
		the transport term satisfies
		$$
		\int_{S\times\O_p^\e}
		(\Gu^\e\cdot\nabla\phi^\e)\varphi_2^\e
		\,\md(x,t)
		=
		-
		\int_{S\times\O_p^\e}
		\phi^\e\Gu^\e\cdot\nabla\varphi_2^\e
		\,\md(x,t).
		$$
		Using \eqref{TSL02}, \eqref{Lam01}, and
		\eqref{eq:scalar-test-convergence}, we obtain
		\begin{equation}\label{eq:limit-phase-convection}
			\int_{S\times\O_p^\e}
			(\Gu^\e\cdot\nabla\phi^\e)\varphi_2^\e
			\,\md(x,t)
			\longrightarrow
			-
			\sqrt{\lambda}
			\int_{\O_T}
			\phi\Gu\cdot\nabla_x\varphi_2
			\,\md(x,t).
		\end{equation}
		Here the contribution containing
		$\nabla_y\widehat\varphi_2$ vanishes because $\phi$ and $\Gu$ are
		independent of $y$ and
		$$
		\int_Y\nabla_y\widehat\varphi_2\,\md y=0.
		$$
		
		The remaining terms follow from
		\eqref{TSL0}--\eqref{TSL02},
		\eqref{eq:scalar-test-convergence}, and the integration identities
		of Proposition~\ref{PropPerforatedUnfolding}. More precisely,
		\begin{equation}\label{eq:limit-phase-terms}
			\begin{aligned}
				\int_{S\times\O_p^\e}
				\nabla\mu^\e\cdot\nabla\varphi_2^\e
				\,\md(x,t)
				&\longrightarrow
				\int_{\O_T\times Y}
				\G(\mu,\widehat\mu)
				\cdot
				\G(\varphi_2,\widehat\varphi_2)
				\,\md(y,x,t),
				\\
				\int_{S\times\O_p^\e}
				G(\phi^\e)\varphi_2^\e
				\,\md(x,t)
				&\longrightarrow
				\int_{\O_T}
				G(\phi)\varphi_2
				\,\md(x,t),\quad
				\int_{S\times\O_p^\e}
				\mu^\e\varphi_3^\e
				\,\md(x,t)
				\longrightarrow
				\int_{\O_T}
				\mu\varphi_3
				\,\md(x,t),
				\\
				\int_{S\times\O_p^\e}
				\nabla\phi^\e\cdot\nabla\varphi_3^\e
				\,\md(x,t)
				&\longrightarrow
				\int_{\O_T\times Y}
				\G(\phi,\widehat\phi)
				\cdot
				\G(\varphi_3,\widehat\varphi_3)
				\,\md(y,x,t),
				\\
				\int_{S\times\O_p^\e}
				F'(\phi^\e)\varphi_3^\e
				\,\md(x,t)
				&\longrightarrow
				\int_{\O_T}
				F'(\phi)\varphi_3
				\,\md(x,t).
			\end{aligned}
		\end{equation}
		Since $\varphi_2^\e$ is compactly supported in time, integration by
		parts gives
		\[
		\int_S
		\left\langle
		\partial_t\phi^\e,\varphi_2^\e
		\right\rangle\,\md t
		=
		-\int_{S\times\O_p^\e}
		\phi^\e\partial_t\varphi_2^\e\,\md(x,t).
		\]
		Using \eqref{TSL0}, \eqref{eq:scalar-test-convergence}, and the integration
		formula of Proposition~\ref{PropPerforatedUnfolding}, we obtain
		\begin{equation}\label{eq:limit-time-phase}
			\int_S
			\left\langle
			\partial_t\phi^\e,\varphi_2^\e
			\right\rangle\,\md t
			\longrightarrow
			-\int_{\O_T}
			\phi\,\partial_t\varphi_2\,\md(x,t).
		\end{equation}
		Combining \eqref{eq:limit-phase-convection}--\eqref{eq:limit-time-phase} gives the
		second identity in \eqref{TSLM01} in the sense of distributions in
		time. The third identity follows directly from the limits in
		\eqref{eq:limit-phase-terms}.
		
		We next use $\varphi_1^\e$ in \eqref{weak-momentum1} and divide the
		resulting identity by $\lambda^\e$. Set
		$$
		\mathbf U^\e
		:=
		\frac{\Gu^\e}{\sqrt{\lambda^\e}}.
		$$
		Since $\Psi_1^\e$ is compactly supported in time, integration by parts
		gives
		\[
		\begin{aligned}
			\frac{1}{\lambda^\e}
			\int_S
			\left\langle
			\partial_t\Gu^\e,\varphi_1^\e
			\right\rangle\,\md t
			&=
			-\int_{S\times\O_p^\e}
			\mathbf U^\e\cdot
			\partial_t\Psi_1^\e\,\md(x,t).
		\end{aligned}
		\]
		Using the strong convergence of the normalized velocity from
		Lemma~\ref{lem:limit-compactness} and the test-function convergence
		\eqref{eq:velocity-test-convergence}, we obtain
		\begin{equation}\label{VelocityTimeLimit}
			\frac{1}{\lambda^\e}
			\int_S
			\left\langle
			\partial_t\Gu^\e,\varphi_1^\e
			\right\rangle\,\md t
			\longrightarrow
			-\int_{\O_T}
			\Gu\cdot\partial_t\varphi_1\,\md(x,t).
		\end{equation}
		
		Using the divergence constraint, the convective term can be written
		as
		$$
		\frac{1}{\lambda^\e}
		\int_{S\times\O_p^\e}
		(\Gu^\e\cdot\nabla)\Gu^\e\cdot\varphi_1^\e
		\,\md(x,t)
		=
		-
		\sqrt{\lambda^\e}
		\int_{S\times\O_p^\e}
		(\mathbf U^\e\otimes\mathbf U^\e):
		\nabla\Psi_1^\e
		\,\md(x,t).
		$$
		By \eqref{TSL02}, \eqref{Lam01}, and
		\eqref{eq:velocity-test-convergence},
		\begin{equation}\label{eq:limit-convection-velocity}
			\frac{1}{\lambda^\e}
			\int_{S\times\O_p^\e}
			(\Gu^\e\cdot\nabla)\Gu^\e\cdot\varphi_1^\e
			\,\md(x,t)
			\longrightarrow
			-
			\sqrt{\lambda}
			\int_{\O_T}
			(\Gu\otimes\Gu):\nabla_x\varphi_1
			\,\md(x,t).
		\end{equation}
		The microscopic contribution vanishes since $\int_Y\nabla_y\widehat\varphi_1\,\md y=0$.
		
		For the capillary term, \eqref{TSL02} and
		\eqref{eq:velocity-test-convergence} imply
		\begin{equation}\label{eq:limit-capillary}
			\begin{aligned}
				\frac{1}{\lambda^\e}
				\int_{S\times\O_p^\e}
				\lambda^\e\phi^\e\nabla\mu^\e
				\cdot\varphi_1^\e
				\,\md(x,t)
				&=
				\sqrt{\lambda^\e}
				\int_{S\times\O_p^\e}
				\phi^\e\nabla\mu^\e
				\cdot\Psi_1^\e
				\,\md(x,t)\\
				&\longrightarrow
				\sqrt{\lambda}
				\int_{\O_T\times Y}
				\phi
				\G(\mu,\widehat\mu)
				\cdot\varphi_1
				\,\md(y,x,t)=
				\sqrt{\lambda}
				\int_{\O_T}
				\phi\nabla_x\mu\cdot\varphi_1
				\,\md(x,t).
			\end{aligned}
		\end{equation}
		The last identity follows from the periodicity and zero mean of
		$\widehat\mu$.
		
		For the viscous term, we have
		$$
		\frac{1}{\lambda^\e}
		\int_{S\times\O_p^\e}
		\A^\e D\Gu^\e:D\varphi_1^\e
		\,\md(x,t)
		=
		\int_{S\times\O_p^\e}
		\A^\e D\mathbf U^\e:D\Psi_1^\e
		\,\md(x,t).
		$$
		Using \eqref{PeriodicCoefficientConvergence},
		\eqref{TSL01}, and
		\eqref{eq:velocity-test-convergence}, we obtain
		\begin{equation}\label{eq:limit-viscosity-subcritical}
			\begin{aligned}
				&\frac{1}{\lambda^\e}
				\int_{S\times\O_p^\e}
				\A^\e D\Gu^\e:D\varphi_1^\e
				\,\md(x,t)
				\longrightarrow
				\int_{\O_T\times Y}
				\A(t,y)
				\D(\Gu,\widehat\Gu):
				\D(\varphi_1,\widehat\varphi_1)
				\,\md(y,x,t).
			\end{aligned}
		\end{equation}
		
		Finally, the normalized force convergence
		\eqref{Assgg01} and the strong convergence of $\Psi_1^\e$ give
		\begin{equation}\label{eq:limit-force}
			\frac{1}{\lambda^\e}
			\int_{S\times\O_p^\e}
			\Gg^\e\cdot\varphi_1^\e
			\,\md(x,t)
			\longrightarrow
			\int_{\O_T}
			\Gg\cdot\varphi_1
			\,\md(x,t).
		\end{equation}
		Combining
		\eqref{VelocityTimeLimit}--\eqref{eq:limit-force}
		gives the first identity in \eqref{TSLM01}.
		
		At this stage, the first two identities have been obtained for smooth
		test functions, with the time derivatives understood in the
		distributional sense:
		\[
		-\int_{\O_T}
		\Gu\cdot\partial_t\varphi_1\,\md(x,t),
		\qquad
		-\int_{\O_T}
		\phi\,\partial_t\varphi_2\,\md(x,t).
		\]
		
		Taking $\widehat\varphi_1=0$ in the first identity and using the
		estimates for the viscous, convective, capillary, and forcing terms as in Step 5 of Lemma \ref{Est01}
		gives $\partial_t\Gu
		\in
		L^{4/3}(S;H^{-1}_{0,\div}(\O)^3)$.
		Indeed, the convective term belongs to
		$L^{4/3}(S;H^{-1}_{0,\div}(\O)^3)$, while the remaining terms belong
		to $L^2(S;H^{-1}_{0,\div}(\O)^3)$.
		
		Similarly, taking $\widehat\varphi_2=0$ in the second identity and
		using the diffusion, transport, and source estimates gives $\partial_t\phi
		\in
		L^2(S;H^{-1}(\O))$.
		Consequently, the distributional time terms can be written as
		\[
		\int_S
		\langle\partial_t\Gu,\varphi_1\rangle\,\md t,
		\qquad
		\int_S
		\langle\partial_t\phi,\varphi_2\rangle\,\md t.
		\]
%
		We have proved the identities for test functions in
		$ \GT_i$, $i\in\{1,2,3\}$. The regularity obtained above,
		the uniform ellipticity and boundedness of $\mathbb A$, and the
		Sobolev embeddings in dimension three show that every term is
		continuous on the corresponding space $\GL_i$. The density of
		$\GT_i$ in $\GL_i$ therefore extends the identities to
		all admissible test functions.
	\end{proof}
	
\subsection{Homogenized system via cell problems}
\label{SSecHomogenizedSystem}
%

We first derive the cell problems associated with \eqref{TSLM01}. 

\textbf{Scalar correctors:}
Taking $\varphi_2=0$ in the phase equation and $\varphi_3=0$ in the
chemical-potential equation of \eqref{TSLM01}, and then localizing with
respect to $(t,x)\in\O_T$, gives
\begin{equation}\label{CellScalar}
	\left.\begin{aligned}
		\int_Y
		\left(
		\nabla_x\mu+\nabla_y\widehat\mu
		\right)
		\cdot\nabla_y\widehat\varphi_2
		\,\md y
		&=0
		&&
		\forall\,
		\widehat\varphi_2\in\mathcal H_\#,
		\\
		\int_Y
		\left(
		\nabla_x\phi+\nabla_y\widehat\phi
		\right)
		\cdot\nabla_y\widehat\varphi_3
		\,\md y
		&=0
		&&
		\forall\,
		\widehat\varphi_3\in\mathcal H_\#,
	\end{aligned}\right\}\quad\text{for a.e.\ $(t,x)\in\O_T$.}
\end{equation}
Since $\nabla_x\phi$ and $\nabla_x\mu$ are independent of $y$,
periodicity gives
\[
\int_Y
\nabla_x\phi\cdot\nabla_y\eta\,\md y
=
\int_Y
\nabla_x\mu\cdot\nabla_y\eta\,\md y
=0
\qquad
\forall\,\eta\in\mathcal H_\#.
\]
Consequently, \eqref{CellScalar} reduces to
\[
\int_Y
\nabla_y\widehat\mu\cdot\nabla_y\widehat\varphi_2
\,\md y
=0,
\qquad
\int_Y
\nabla_y\widehat\phi\cdot\nabla_y\widehat\varphi_3
\,\md y
=0.
\]
Taking
$\widehat\varphi_2=\widehat\mu$
and
$\widehat\varphi_3=\widehat\phi$
yields
\[
\nabla_y\widehat\mu=0,
\qquad
\nabla_y\widehat\phi=0.
\]
Since both correctors have zero mean over $Y$, we conclude that
\begin{equation}\label{VanishingScalarCorrectors}
	\widehat\phi=0,
	\qquad
	\widehat\mu=0
	\quad\text{a.e. in }\O_T\times Y.
\end{equation}
Thus the scalar effective tensors are the identity.
\smallskip

\textbf{Velocity corrector:}
Taking $\varphi_1=0$ in the momentum equation of \eqref{TSLM01} and
localizing with respect to $(t,x)\in\O_T$ gives
\begin{equation}\label{CellVelocity}
	\int_Y
	\A(t,y)
	\left(
	D_x\Gu(t,x)+D_y\widehat\Gu(t,x,y)
	\right)
	:
	D_y\widehat\varphi_1(y)
	\,\md y
	=0
	\qquad
	\forall\,
	\widehat\varphi_1\in\Vc_\#
\end{equation}
for a.e.\ $(t,x)\in\O_T$.

For every
$M\in\mathbb R_{\mathrm{sym}}^{3\times3}$
and a.e.\ $t\in S$, let $\chi_M(t,\cdot)\in\Vc_\#$
be the unique solution of
\begin{equation}\label{VelocityBasisCellProblem}
	\int_Y
	\A(t,y)
	\left(
	M+D_y\chi_M(t,y)
	\right)
	:
	D_y\psi(y)
	\,\md y
	=0
	\qquad
	\forall\,\psi\in\Vc_\#.
\end{equation}
The existence and uniqueness of $\chi_M(t,\cdot)$ follow from the
boundedness and uniform ellipticity of $\A$, the periodic Korn inequality,
and the Lax--Milgram theorem. The map $M\longmapsto\chi_M(t,\cdot)$
is linear.

Comparing \eqref{CellVelocity} and
\eqref{VelocityBasisCellProblem}, and using uniqueness, gives
\begin{equation}\label{VelocityCorrectorRepresentation}
	\widehat\Gu(t,x,\cdot)
	=
	\chi_{D_x\Gu(t,x)}(t,\cdot),\quad\text{for a.e.\ $(t,x)\in\O_T$}.
\end{equation}
For
$M,N\in\mathbb R_{\mathrm{sym}}^{3\times3}$,
define the effective viscosity tensor by
\begin{equation}\label{EffectiveViscosityTensor}
	\A^{\hom}(t)M:N
	:=
	\int_Y
	\A(t,y)
	\left(
	M+D_y\chi_M(t,y)
	\right)
	:
	\left(
	N+D_y\chi_N(t,y)
	\right)
	\,\md y.
\end{equation}
Testing \eqref{VelocityBasisCellProblem} for $\chi_M$ with
$\psi=\chi_N$ shows that
\begin{equation}\label{EffectiveViscosityFlux}
	\A^{\hom}(t)M:N
	=
	\int_Y
	\A(t,y)
	\left(
	M+D_y\chi_M(t,y)
	\right)
	:
	N
	\,\md y.
\end{equation}
Equivalently,
\[
\A^{\hom}(t)M
=
\int_Y
\A(t,y)
\left(
M+D_y\chi_M(t,y)
\right)
\,\md y.
\]
The tensor $\A^{\hom}(t)$ describes the effective viscous response of the
homogenized system.
\begin{lemma}
	\label{lem:effective-tensors}
	For a.e.\ $t\in S$, the effective viscosity tensor
	$\A^{\hom}(t)$ is symmetric, bounded, and uniformly elliptic.
	More precisely, there exist constants $c_1,c_2>0$, independent of
	$t$, such that
	\[
	c_1|M|^2
	\leq
	\A^{\hom}(t)M:M
	\leq
	c_2|M|^2
	\qquad
	\forall\,
	M\in\mathbb R_{\mathrm{sym}}^{3\times3}.
	\]
	Moreover,
	\[
	\A^{\hom}
	\in
	L^\infty
	\left(
	S;
	\mathcal L
	\left(
	\mathbb R_{\mathrm{sym}}^{3\times3}
	\right)
	\right).
	\]
\end{lemma}

\begin{proof}
	For a.e.\ $t\in S$ and every
	$M\in\mathbb R_{\mathrm{sym}}^{3\times3}$,
	the cell problem has the variational characterization
	\[
	\A^{\hom}(t)M:M
	=
	\min_{\psi\in\Vc_\#}
	\int_Y
	\A(t,y)
	\left(
	M+D_y\psi
	\right)
	:
	\left(
	M+D_y\psi
	\right)
	\,\md y.
	\]
	The symmetry of $\A^{\hom}(t)$ follows from the symmetry of
	$\A(t,y)$ and \eqref{EffectiveViscosityTensor}.
	
	Let $\kappa_0,\kappa_1>0$ be the uniform boundedness and ellipticity
	constants of $\A$. Since $\int_Y D_y\psi\,\md y=0$ forall $\psi\in\Vc_\#$,
	we have
	\[
	\int_Y
	|M+D_y\psi|^2\,\md y
	=
	|Y||M|^2
	+
	\int_Y|D_y\psi|^2\,\md y.
	\]
	As $|Y|=1$, the lower ellipticity bound gives
	\[
	\A^{\hom}(t)M:M
	\geq
	\kappa_1|M|^2.
	\]
	Choosing $\psi=0$ in the variational characterization gives
	\[
	\A^{\hom}(t)M:M
	\leq
	\kappa_0|M|^2.
	\]
	The asserted boundedness follows by polarization. Measurability in
	$t$ follows from the parameter-dependent variational problem. See
	\cite{JikovKozlovOleinik94,Hornung97,CDG} for the standard
	periodic homogenization argument.
\end{proof}

\begin{theorem}[Homogenized Navier--Stokes--Cahn--Hilliard system]
	\label{ThHomogenized}
	Let
	$(\Gu^\e,\phi^\e,\mu^\e)$
	be a weak solution of the microscopic problem
	\eqref{MainP01}, and let
	$(\Gu,\phi,\mu)$
	be the limit fields obtained in
	Theorem~\ref{Th02} (see also Lemma \ref{lem:limit-compactness}). Then
	$(\Gu,\phi,\mu)$
	is a weak solution of the homogenized system.
	
	Moreover,
	\[
	\Gu(0)=\Gu_0
	\quad\text{in }H^{-1}_{0,\div}(\O)^3,
	\qquad
	\phi(0)=\phi_0
	\quad\text{in }L^2(\O),\quad\text{and}\quad
	\nabla\cdot\Gu=0
	\quad\text{in }\mathcal D'(\O_T).
	\]
	More precisely,
	$(\Gu,\phi,\mu)$ satisfies
	\begin{equation}
		\label{HomogenizedWeakSystem}
		\begin{aligned}
			&
			\int_S
			\langle\partial_t\Gu,\varphi_1\rangle
			\,\md t
			-
			\sqrt{\lambda}
			\int_{\O_T}
			(\Gu\otimes\Gu):
			\nabla\varphi_1
			\,\md(x,t)
			+
			\int_{\O_T}
			\A^{\hom}(t)
			D\Gu:
			D\varphi_1
			\,\md(x,t)
			\\
			&\hspace{2cm}
			+
			\sqrt{\lambda}
			\int_{\O_T}
			\phi\nabla\mu\cdot\varphi_1
			\,\md(x,t)
			=
			\int_{\O_T}
			\Gg\cdot\varphi_1
			\,\md(x,t),
			\\[2mm]
			&
			\int_S
			\langle\partial_t\phi,\varphi_2\rangle
			\,\md t
			-
			\sqrt{\lambda}
			\int_{\O_T}
			\phi\Gu\cdot\nabla\varphi_2
			\,\md(x,t)
			+
			\int_{\O_T}
			\nabla\mu\cdot\nabla\varphi_2
			\,\md(x,t)
			+
			\int_{\O_T}
			G(\phi)\varphi_2
			\,\md(x,t)
			=
			0,
			\\[2mm]
			&
			\int_{\O_T}
			\mu\varphi_3
			\,\md(x,t)
			-
			\int_{\O_T}
			\nabla\phi\cdot\nabla\varphi_3
			\,\md(x,t)
			-
			\int_{\O_T}
			F'(\phi)\varphi_3
			\,\md(x,t)
			=
			0
		\end{aligned}
	\end{equation}
	for every $\varphi_1
	\in
	L^4(S;H^1_{0,\div}(\O)^3)$, $\varphi_2,\varphi_3
	\in
	L^2(S;H^1(\O))$.
\end{theorem}

\begin{proof}
	By \eqref{VanishingScalarCorrectors}, $\widehat\phi=0$ and $\widehat\mu=0$.
	Consequently, the scalar diffusion terms in the two-scale limit
	system \eqref{TSLM01} reduce to $\int_{\O_T}
	\nabla\mu\cdot\nabla\varphi_2
	\,\md(x,t)$, and $\int_{\O_T}
	\nabla\phi\cdot\nabla\varphi_3
	\,\md(x,t)$,
	respectively.
	
	For the velocity equation,
	\eqref{VelocityCorrectorRepresentation} gives
	\[
	\widehat\Gu(t,x,\cdot)
	=
	\chi_{D_x\Gu(t,x)}(t,\cdot).
	\]
	Using the velocity cell problem
	\eqref{VelocityBasisCellProblem}, the effective-flux identity
	\eqref{EffectiveViscosityFlux}, and $|Y|=1$, we obtain
	\[
	\begin{aligned}
		&
		\int_{\O_T\times Y}
		\A(t,y)
		\left(
		D_x\Gu+D_y\widehat\Gu
		\right)
		:
		\left(
		D_x\varphi_1+D_y\widehat\varphi_1
		\right)
		\,\md(y,x,t)
		=
		\int_{\O_T}
		\A^{\hom}(t)
		D\Gu:
		D\varphi_1
		\,\md(x,t).
	\end{aligned}
	\]
	Substituting these identities into the two-scale limit system
	\eqref{TSLM01} eliminates the microscopic variables and yields
	\eqref{HomogenizedWeakSystem}.
	The initial conditions follow from
	Lemma~\ref{lem:limit-initial-data}, while the divergence constraint
	follows from Theorem~\ref{Th02}.
\end{proof}

\begin{remark}
	The homogenized system \eqref{HomogenizedWeakSystem} distinguishes
	two regimes according to the limit parameter $\lambda$.
	
	If $\lambda=0$, the convection term in the momentum equation, the
	transport term in the phase equation, and the capillary coupling
	term vanish. The velocity then satisfies an effective unsteady
	Stokes system, while $(\phi,\mu)$ satisfies a Cahn--Hilliard system
	with source term $G(\phi)$. The two systems are decoupled.
	
	If $\lambda>0$, the convection, transport, and capillary terms
	remain, and the limit retains the full
	Navier--Stokes--Cahn--Hilliard coupling. The periodic viscosity
	oscillations are represented by $\A^{\hom}(t)$. The subcritical
	obstacles do not modify the scalar diffusion operators and leave no
	additional term in the homogenized equations.
\end{remark}

\begin{remark}
	The homogenized systems
	\eqref{HomStr01} and \eqref{LHomStr01}
	are stated in strong form for readability and include a pressure
	$p$. The analysis is carried out in the solenoidal formulation. We
	neither estimate the microscopic pressure $p^\e$ nor identify a
	limit of $p^\e$.
	
	The pressure in the homogenized systems is the Lagrange multiplier
	associated with the constraint
	$\nabla\cdot\Gu=0$. It can be recovered from
	\eqref{HomogenizedWeakSystem} by the standard de~Rham theorem on
	the fixed domain $\O$. This pressure is not asserted to be the
	limit of the microscopic pressures.
\end{remark}
	
The following convergence holds for the same subsequence as in
Lemma~\ref{lem:limit-compactness}.

\begin{theorem}[Time-integrated energy convergence]
	\label{Th04}
	We have
	\begin{equation}\label{EC01}
		\lim_{\e\to0}
		\int_S
		\frac{
			\GT^\e(\Gu^\e(t),\phi^\e(t))
		}{
			\lambda^\e
		}
		\,\md t
		=
		\int_S
		\GT(\Gu(t),\phi(t))
		\,\md t,
	\end{equation}
	where
	\[
	\GT(\Gu,\phi)
	:=
	\frac12\int_\O|\Gu|^2\,\md x
	+
	\frac12\int_\O|\nabla\phi|^2\,\md x
	+
	\int_\O F(\phi)\,\md x.
	\]
\end{theorem}

\begin{proof}
	We first prove the strong convergence of the phase gradients in the
	fluid part of the unfolded cell:
	\begin{equation}\label{SS}
		\chi_{p,\e}
		\Te^\ast(\nabla\phi^\e)
		\longrightarrow
		\nabla\phi
		\quad\text{strongly in }
		L^2(\O_T\times Y)^3.
	\end{equation}
	Here we have used
	$\widehat\phi=0$, as established in
	\eqref{VanishingScalarCorrectors}.
	Since $\chi_{p,\e}\longrightarrow1$ strongly in $L^q(Y)$
	for every finite $q$, the weak convergence in
	Lemma~\ref{lem:limit-compactness} gives
	\[
	\chi_{p,\e}
	\Te^\ast(\nabla\phi^\e)
	\rightharpoonup
	\nabla\phi
	\quad\text{weakly in }L^2(\O_T\times Y)^3.
	\]
	Moreover, the unfolding identity gives
	\[
	\begin{aligned}
		\int_{\O_T\times Y}
		\chi_{p,\e}
		\left|
		\Te^\ast(\nabla\phi^\e)
		\right|^2
		\,\md(y,x,t)
		&=
		\int_{S\times\widehat\O_p^\e}
		|\nabla\phi^\e|^2
		\,\md(x,t)
		&\leq
		\int_{S\times\O_p^\e}
		|\nabla\phi^\e|^2
		\,\md(x,t).
	\end{aligned}
	\]
	Hence,
	\begin{equation}\label{GradientLowerBound}
		\|\nabla\phi\|_{L^2(\O_T)}^2
		\leq
		\liminf_{\e\to0}
		\int_{\O_T\times Y}
		\chi_{p,\e}
		\left|
		\Te^\ast(\nabla\phi^\e)
		\right|^2
		\,\md(y,x,t)
		\leq
		\liminf_{\e\to0}
		\int_{S\times\O_p^\e}
		|\nabla\phi^\e|^2
		\,\md(x,t).
	\end{equation}
	
	Testing the microscopic chemical-potential identity with
	$\phi^\e$ gives
	\begin{equation}\label{MicroscopicGradientIdentity}
		\int_{S\times\O_p^\e}
		|\nabla\phi^\e|^2
		\,\md(x,t)
		=
		\int_{S\times\O_p^\e}
		\mu^\e\phi^\e
		\,\md(x,t)
		-
		\int_{S\times\O_p^\e}
		F'(\phi^\e)\phi^\e
		\,\md(x,t).
	\end{equation}
	The convergences established in
	Lemma~\ref{lem:limit-compactness}, together with the integration
	formulas of Proposition~\ref{PropPerforatedUnfolding}, give
	\[
	\begin{aligned}
		\int_{S\times\O_p^\e}
		\mu^\e\phi^\e
		\,\md(x,t)
		&\longrightarrow
		\int_{\O_T}
		\mu\phi
		\,\md(x,t),\quad
		\int_{S\times\O_p^\e}
		F'(\phi^\e)\phi^\e
		\,\md(x,t)
		&\longrightarrow
		\int_{\O_T}
		F'(\phi)\phi
		\,\md(x,t).
	\end{aligned}
	\]
	On the other hand, taking
	$\varphi_3=\phi$ in the limit chemical-potential equation gives
	\[
	\int_{\O_T}
	\mu\phi
	\,\md(x,t)
	-
	\int_{\O_T}
	|\nabla\phi|^2
	\,\md(x,t)
	-
	\int_{\O_T}
	F'(\phi)\phi
	\,\md(x,t)
	=0.
	\]
	Therefore,
	\begin{equation}\label{GradientNormConvergence}
		\lim_{\e\to0}
		\int_{S\times\O_p^\e}
		|\nabla\phi^\e|^2
		\,\md(x,t)
		=
		\int_{\O_T}
		|\nabla\phi|^2
		\,\md(x,t).
	\end{equation}
	Combining
	\eqref{GradientLowerBound} and
	\eqref{GradientNormConvergence} shows that
	\[
	\lim_{\e\to0}
	\int_{\O_T\times Y}
	\chi_{p,\e}
	\left|
	\Te^\ast(\nabla\phi^\e)
	\right|^2
	\,\md(y,x,t)
	=
	\|\nabla\phi\|_{L^2(\O_T)}^2.
	\]
	The weak convergence and convergence of the norms prove
	\eqref{SS}.
	
	We now pass to the limit in each component of the energy. By
	Lemmas~\ref{Lem08+} and~\ref{lem:limit-compactness}, we have
	\begin{equation}\label{KineticEnergyConvergence}
		\int_{S\times\O_p^\e}
		\frac{|\Gu^\e|^2}{\lambda^\e}
		\,\md(x,t)
		\longrightarrow
		\int_{\O_T}
		|\Gu|^2
		\,\md(x,t),\quad 		\int_{S\times\O_p^\e}
		F(\phi^\e)
		\,\md(x,t)
		\longrightarrow
		\int_{\O_T}
		F(\phi)
		\,\md(x,t).
	\end{equation}
	Finally,
	the above convergences combined with \eqref{GradientNormConvergence}
	yield \eqref{EC01}.
\end{proof}

	\section*{Acknowledgements}
	The authors would like to thank Markus Gahn (University of Augsburg) for his comments on the Poincaré inequality for velocity fields satisfying
	the free-slip boundary condition. This work was largely completed when Amartya Chakrabortty was a PhD student and research assistant in the department of Processes and Material Simulation, Fraunhofer ITWM, Kaiserslautern, Germany. Haradhan Dutta gratefully acknowledges the PMRF (ID-2402788) scheme and IIT Kharagpur for funding his PhD position.
	
	
	\section*{Data availability }
	No datasets were generated or analysed during the current study.
	
	\section*{Declarations and Competing interests}
	The authors declare no competing interests.
	
\appendix
\section{Auxiliary results}\label{AppSec}
\subsection{Proofs of some preliminary results}\label{App01}
\begin{proof}[\textbf{Proof of Lemma \ref{Bogov01+}}]
	Let $\mathcal B_T:L^2_0(T)\to H^1_0(T)^3$ be a fixed Bogovski\u{\i} operator
	satisfying
	\[
	\div_z\mathcal B_T(f)=f,
	\qquad
	\|\mathcal B_T(f)\|_{H^1(T)}\leq C\|f\|_{L^2(T)}.
	\]
	For $\kappa\in\Kc_\e$, set $g_{\e,\kappa}(z):=g(x_{\e,\kappa}+a_\e z)$,
	$z\in T$. The componentwise mean-zero condition gives $g_{\e,\kappa}\in L^2_0(T)$.
	Define, for $x\in T_{\e,\kappa}$,
	\[
	\mathcal B^\e_s(g)(x)
	:=a_\e\,\mathcal B_T(g_{\e,\kappa})\Bigl(\frac{x-x_{\e,\kappa}}{a_\e}\Bigr).
	\]
	This operator is linear, maps into $H^1_0(\O_s^\e)^3$, and satisfies
	$\div\mathcal B^\e_s(g)=g$ in $\O_s^\e$. Rescaling on each obstacle gives
	\[
	\|\nabla\mathcal B^\e_s(g)\|_{L^2(T_{\e,\kappa})}\leq C\|g\|_{L^2(T_{\e,\kappa})},
	\qquad
	\|\mathcal B^\e_s(g)\|_{L^2(T_{\e,\kappa})}\leq Ca_\e\|g\|_{L^2(T_{\e,\kappa})}.
	\]
	Summing over $\kappa\in\Kc_\e$ and using $a_\e=\e^\alpha\leq1$ yields
	\eqref{Bogov02+-gradient}, with $C_s$ independent of $\e$.
\end{proof}

\begin{proof}[\textbf{Proof of Lemma \ref{lem:scalar-two-extensions}}]
	The primal bound in \eqref{ScalarPrimalBound} is
	\eqref{eq:extension-primal-estimates} applied for a.e.\ $t\in S$.
	
	The zero extension inherits its time derivative from restriction. For
	$\eta\in H^1(\O)$, set
	\[
	\left\langle
	\partial_t\overline v^\e,\eta
	\right\rangle_
	{H^{-1}(\O),H^1(\O)}
	:=
	\left\langle
	\partial_tv^\e,
	\eta|_{\O_p^\e}
	\right\rangle_
	{H^{-1}(\O_p^\e),H^1(\O_p^\e)}.
	\]
	Since $\O_s^\e$ is time-independent, this functional is the distributional
	derivative of $\overline v^\e$. From
	$\|\eta|_{\O_p^\e}\|_{H^1(\O_p^\e)}\leq\|\eta\|_{H^1(\O)}$ we obtain
	\[
	\|\partial_t\overline v^\e\|_{H^{-1}(\O)}
	\leq
	\|\partial_tv^\e\|_{H^{-1}(\O_p^\e)},
	\]
	which gives the remaining bounds in \eqref{ScalarPrimalBound}. Moreover, the zero-extension operator is an isometry from
	$L^2(\O_p^\e)$ into $L^2(\O)$. Hence $\overline v^\e\in H^1(S;H^{-1}(\O))$.
	
	Since $\widetilde v^\e=v^\e$ in $\O_p^\e$, $\widetilde v^\e-\overline v^\e$ is equal to $0$ in $\O_p^\e$, and equal to $\widetilde v^\e$ in $\O_s^\e$.
	Hence, by H\"older's inequality and the Sobolev embedding
	$H^1(\O)\hookrightarrow L^6(\O)$,
	\[
	\|\widetilde v^\e-\overline v^\e\|_{L^\infty(S;L^2(\O))}
	\leq
	|\O_s^\e|^{1/3}
	\|\widetilde v^\e\|_{L^\infty(S;L^6(\O))}
	\leq
	C|\O_s^\e|^{1/3}.
	\]
	By \eqref{EqSolidVol}, $|\O_s^\e|=O(\e^{3(\alpha-1)})$, so
	$|\O_s^\e|^{1/3}\leq C\e^{\alpha-1}=Ca_\e/\e$, which is
	\eqref{ScalarExtensionsDifference}.
	
	By the standard Hilbert-triple result $v^\e\in C([0,T];L^2(\O_p^\e))$.
	Since both extension operators are time-independent and bounded on the
	corresponding $L^2$ spaces $\widetilde v^\e,\overline v^\e
	\in C([0,T];L^2(\O))$.
	Set
	\[
	d_\e
	:=
	\|\widetilde v^\e-\overline v^\e\|_
	{L^\infty(S;L^2(\O))},
	\]
	so that $d_\e\to0$. For $t,s\in[0,T]$, the interpolation
	inequality
	\[
	\|q\|_{L^2(\O)}^2
	\leq
	\|q\|_{H^1(\O)}
	\|q\|_{H^{-1}(\O)}
	\]
	applied to $q=\widetilde v^\e(t)-\widetilde v^\e(s)$ gives
	\begin{align*}
		\|\widetilde v^\e(t)-\widetilde v^\e(s)\|_
		{L^2(\O)}^2
		&\leq
		C
		\|\widetilde v^\e(t)-\widetilde v^\e(s)\|_
		{H^{-1}(\O)}\leq
		C
		\left(
		\|\overline v^\e(t)-\overline v^\e(s)\|_
		{H^{-1}(\O)}
		+
		2d_\e
		\right)\leq
		C\left(
		|t-s|^{1/2}
		+
		2d_\e
		\right),
	\end{align*}
	the last step using the time-derivative bound in
	\eqref{ScalarPrimalBound}. Thus $\{\widetilde v^\e\}$ is asymptotically
	equicontinuous in $L^2(\O)$. Since it is uniformly bounded in
	$L^\infty(S;H^1(\O))$ and $H^1(\O)\Subset L^2(\O)$, the Arzel\`a--Ascoli
	theorem gives, after extraction,
	\[
	\widetilde v^\e
	\longrightarrow v
	\quad\text{strongly in }
	C([0,T];L^2(\O)).
	\]
	The uniform bound in
	$L^\infty(S;H^1(\O))$ also gives, after a further extraction,
	\[
	\widetilde v^\e
	\overset{*}{\rightharpoonup}
	v
	\quad\text{in }L^\infty(S;H^1(\O)).
	\]
	Consequently, $v\in\frak C$.
	The convergence of $\overline v^\e$ follows from
	\eqref{ScalarExtensionsDifference}.
\end{proof}

\begin{proof}[\textbf{Proof of Lemma \ref{LemD02}}]
	Let $v\in V_\e^\sigma$ and set $w:=(E^\e v_1,E^\e v_2,E^\e v_3)$, with $E^\e$
	the scalar extension of Lemma~\ref{Emain}. Then $w|_{\O_p^\e}=v$, $w=0$ on
	$\partial\O$, and
	\begin{equation}\label{VectorScalarExtensionEstimates}
		\|w\|_{L^2(\O)}\leq C\|v\|_{L^2(\O_p^\e)},
		\qquad
		\|\nabla w\|_{L^2(\O)}\leq C\|\nabla v\|_{L^2(\O_p^\e)}.
	\end{equation}
	Set $g:=(\div w)|_{\O_s^\e}$. For each $\kappa\in\Kc_\e$, the traces of $w$ from
	the solid side and of $v$ from the fluid side agree on $\partial T_{\e,\kappa}$,
	so with $\Gn_s^\e=-\Gn^\e$ and the impermeability condition,
	\[
	\int_{T_{\e,\kappa}}g\,\md x
	=\int_{\partial T_{\e,\kappa}}w\cdot\Gn_s^\e\,\md S
	=-\int_{\partial T_{\e,\kappa}}v\cdot\Gn^\e\,\md S=0,
	\]
	hence $g\in L^2_{0,\mathrm{comp}}(\O_s^\e)$. Let $z_\e:=\mathcal B_s^\e(g)\in H_0^1(\O_s^\e)^3$
	and extend $z_\e$ by zero to $\O$. Since
	$\O_s^\e\Subset\O$ and $z_\e$ has zero trace on
	$\partial\O_s^\e$, its zero extension belongs to $H_0^1(\O)^3$. The
	bounds of Lemma~\ref{Bogov01+} and $\|g\|_{L^2(\O_s^\e)}\leq C\|\nabla w\|_{L^2(\O_s^\e)}$
	give
	\begin{equation}\label{SolidCorrectionH1Bounds}
		\|\nabla z_\e\|_{L^2(\O_s^\e)}\leq C\|\nabla w\|_{L^2(\O_s^\e)},
		\qquad
		\|z_\e\|_{L^2(\O_s^\e)}\leq Ca_\e\|\nabla w\|_{L^2(\O_s^\e)}.
	\end{equation}
	Define $\mathcal E^\e v:=w-z_\e$. Since $z_\e=0$ in $\O_p^\e$,
	$(\mathcal E^\e v)|_{\O_p^\e}=v$, and
	\[
	\div(\mathcal E^\e v)=
	\begin{cases}
		\div v=0, & \text{in }\O_p^\e,\\
		g-\div z_\e=0, & \text{in }\O_s^\e.
	\end{cases}
	\]
	Since $w,z_\e\in H^1(\O)^3$, we have
	$\mathcal E^\e v\in H^1(\O)^3$. The identities above on
	$\O_p^\e$ and $\O_s^\e$ therefore imply
	\[
	\div(\mathcal E^\e v)=0
	\quad\text{in }\mathcal D'(\O),\quad\text{and}\quad \mathcal E^\e v=0,\quad\text{on $\partial\O$}
	\]
	since $w=0$ there and $z_\e$ is supported in
	$\O_s^\e\Subset\O$. Thus $\mathcal E^\e v\in V^\sigma$. 
	The linearity of $\mathcal E^\e$ follows from the linearity of
	$E^\e$, the divergence operator, and $\mathcal B_s^\e$.
	Combining
	\eqref{VectorScalarExtensionEstimates} and \eqref{SolidCorrectionH1Bounds},
	\[
	\|\mathcal E^\e v\|_{L^2(\O)}\leq\|w\|_{L^2(\O)}+\|z_\e\|_{L^2(\O_s^\e)}
	\leq C\bigl(\|v\|_{L^2(\O_p^\e)}+a_\e\|\nabla v\|_{L^2(\O_p^\e)}\bigr),
	\qquad
	\|\nabla\mathcal E^\e v\|_{L^2(\O)}\leq C\|\nabla v\|_{L^2(\O_p^\e)},
	\]
	which are \eqref{SolenoidalPrimalExtensionBounds}; adding them and using
	$a_\e\leq1$ gives \eqref{SolenoidalPrimalExtensionH1Bound}.
	
	For the zero extension, the norm identity in \eqref{SolenoidalZeroExtension} is
	immediate. 
	For $\phi\in C_c^\infty(\O)$, integration by parts gives
	\[
	\begin{aligned}
		\int_\O \overline v^\e\cdot\nabla\phi\,\md x
		&=
		\int_{\O_p^\e}v\cdot\nabla\phi\,\md x=
		-\int_{\O_p^\e}(\div v)\phi\,\md x
		+
		\int_{\Gamma_s^\e}
		(v\cdot\Gn^\e)\phi\,\md S
		=0.
	\end{aligned}
	\]
	Thus
	\[
	\div\overline v^\e=0
	\quad\text{in }\mathcal D'(\O).
	\]	
	Finally, since $\mathcal E^\e v^\e=v^\e$ in $\O_p^\e$, the difference
	$\mathcal E^\e v^\e-\overline v^\e$ vanishes on $\O_p^\e$ and equals
	$\mathcal E^\e v^\e$ on $\O_s^\e$. By H\"older's inequality
	($\tfrac12=\tfrac13+\tfrac16$), the Sobolev embedding
	$H^1(\O)\hookrightarrow L^6(\O)$, and \eqref{SolenoidalPrimalExtensionH1Bound},
	\[
	\|\mathcal E^\e v^\e-\overline v^\e\|_{L^2(S\times\O)}
	=\|\mathcal E^\e v^\e\|_{L^2(S;L^2(\O_s^\e))}
	\leq|\O_s^\e|^{1/3}\|\mathcal E^\e v^\e\|_{L^2(S;L^6(\O))}
	\leq C|\O_s^\e|^{1/3}\|v^\e\|_{L^2(S;H^1(\O_p^\e))}.
	\]
	Since $|\O_s^\e|\leq Ca_\e^3/\e^3$, we have $|\O_s^\e|^{1/3}\leq Ca_\e/\e$, which
	gives \eqref{PrimalZeroExtensionDifference}. As $a_\e/\e=\e^{\alpha-1}\to0$, the
	stated convergence follows.
\end{proof}

\subsection{Solenoidal test functions}
\label{sec:test_functions}

{
	\begin{lemma}[Subcritical admissible solenoidal test functions]
		\label{lem:test-functions-subcritical}
		Let $\alpha>3$ and let
		$(\varphi_1,\widehat\varphi_1)\in\GT_1$.
		Then there exists $\varphi_1^\e
		\in
		L^2\bigl(S;H^1_{0,\div}(\O_p^\e)^3\bigr)
		\subset
		L^2\bigl(S;\mathbf H^1_{\div}(\O_p^\e)\bigr)$
		such that
		\begin{equation}\label{eq:test-subcritical-uniform-H1}
			\frac{1}{\sqrt{\lambda^\e}}
			\|\varphi_1^\e\|_
			{L^2(S;H^1(\O_p^\e)^3)}
			\leq C,
		\end{equation}
		and
		\begin{equation}\label{eq:test-subcritical-convergence}
			\begin{aligned}
				\frac{1}{\sqrt{\lambda^\e}}
				\Te^\ast(\varphi_1^\e)
				&\to
				\varphi_1
				&&\text{strongly in }
				L^2(\O_T\times Y)^3,
				\\
				\frac{1}{\sqrt{\lambda^\e}}
				\Te^\ast(\nabla\varphi_1^\e)
				&\to
				\nabla_x\varphi_1
				+
				\nabla_y\widehat\varphi_1
				&&\text{strongly in }
				L^2(\O_T\times Y)^{3\times3},
				\\
				\frac{1}{\sqrt{\lambda^\e}}
				\Te^\ast(\partial_t\varphi_1^\e)
				&\to
				\partial_t\varphi_1
				&&\text{strongly in }
				L^2(\O_T\times Y)^3.
			\end{aligned}
		\end{equation}
	\end{lemma}
	
	\begin{proof}
		We first construct a solenoidal test field on the full domain.
		Set $P_\e(t,x)
		:=
		\e\widehat\varphi_1
		\left(
		t,x,\frac{x}{\e}
		\right)$.
		Since
		$\nabla_y\cdot\widehat\varphi_1=0$, the chain rule gives
		\[
		\nabla_x\cdot P_\e
		=
		\e
		\nabla_x\cdot\widehat\varphi_1
		\left(
		t,x,\frac{x}{\e}
		\right)
		=:
		g_\e.
		\]
		The field $P_\e$ has compact support in $\O$. Therefore,
		\[
		\int_\O g_\e(t,x)\,\md x
		=
		\int_\O\nabla_x\cdot P_\e(t,x)\,\md x
		=0
		\qquad
		\text{for every }t\in S.
		\]
		Let $\mathcal B_\O:
		L^q_0(\O)
		\longrightarrow
		W^{1,q}_0(\O)^3$
		be a fixed Bogovski\u{\i} operator on $\O$, and define
		\[
		b_\e(t,\cdot)
		:=
		\mathcal B_\O[g_\e(t,\cdot)],
		\qquad
		Q_\e:=P_\e-b_\e,
		\qquad
		Z_\e:=\varphi_1+Q_\e.
		\]
		Then
		\[
		\nabla_x\cdot Q_\e=0,
		\qquad
		\nabla_x\cdot Z_\e=0
		\quad\text{in }\O.
		\]
		Moreover, since $\varphi_1$ and
		$\widehat\varphi_1$ have compact support with respect to $x$, $Z_\e
		\in
		L^2(S;H^1_{0,\div}(\O)^3)$.
		
		For every $1<q<\infty$, the boundedness of
		$\mathcal B_\O$ gives
		\begin{equation}\label{eq:test-subcritical-global-Bogovskii}
			\|b_\e\|_{L^2(S;W^{1,q}_0(\O)^3)}
			+
			\|\partial_tb_\e\|_
			{L^2(S;W^{1,q}_0(\O)^3)}
			\leq
			C_q\e.
		\end{equation}
		Consequently, for every fixed $q>3$,
		\begin{equation}\label{eq:subcritical-Ze-regularity}
			\begin{aligned}
				\|Z_\e\|_{L^2(S;W^{1,q}(\O)^3)}
				+
				\|\partial_tZ_\e\|_
				{L^2(S;W^{1,q}(\O)^3)}
				&\leq C_q,\quad
				\|Z_\e\|_{L^2(S;L^\infty(\O)^3)}
				+
				\|\partial_tZ_\e\|_
				{L^2(S;L^\infty(\O)^3)}
				&\leq C_q.
			\end{aligned}
		\end{equation}
		Here the second estimate follows from the Sobolev embedding
		$W^{1,q}(\O)\hookrightarrow L^\infty(\O)$ for $q>3$.
		
		The standard periodic unfolding properties and
		\eqref{eq:test-subcritical-global-Bogovskii} yield
		\begin{equation}\label{eq:subcritical-Ze-properties}
			\begin{aligned}
				\Te(Z_\e)
				&\to
				\varphi_1
				&&\text{strongly in }
				L^2(\O_T\times Y)^3,
				\\
				\Te(\nabla Z_\e)
				&\to
				\nabla_x\varphi_1
				+
				\nabla_y\widehat\varphi_1
				&&\text{strongly in }
				L^2(\O_T\times Y)^{3\times3},
				\\
				\Te(\partial_tZ_\e)
				&\to
				\partial_t\varphi_1
				&&\text{strongly in }
				L^2(\O_T\times Y)^3.
			\end{aligned}
		\end{equation}
		
		We next modify $Z_\e$ in neighborhoods of the holes. This is the
		local cutoff-Bogovski\u{\i} construction underlying the restriction
		argument for small holes; see
		\cite{AllaireTinyHolesI,AllaireTinyHolesII}.
		
		Choose $R>0$ such that $\overline T\subset B_R$.
		For every $\kappa\in\Kc_\e$, set $x_{\e,\kappa}
		:=
		\e(\kappa+y_0)$,
		and define
		\[
		A_{\e,\kappa}
		:=
		B_{2Ra_\e}(x_{\e,\kappa})
		\setminus
		\overline{B_{Ra_\e}(x_{\e,\kappa})}.
		\]
		Since $a_\e/\e\to0$, for all sufficiently small $\e$ the balls $B_{2Ra_\e}(x_{\e,\kappa})$, $\kappa\in\Kc_\e$,
		are mutually disjoint and contained in their corresponding
		$\e$-cells.
		
		Let $\eta\in
		C^\infty
		\bigl(
		\overline{B_{2R}\setminus B_R}
		\bigr)$
		satisfy
		\[
		\eta=0
		\quad\text{on }\partial B_R,
		\qquad
		\eta=1
		\quad\text{on }\partial B_{2R}.
		\]
		For every $\kappa\in\Kc_\e$, define
		\[
		\eta_{\e,\kappa}(x)
		:=
		\eta
		\left(
		\frac{x-x_{\e,\kappa}}{a_\e}
		\right).
		\]
		Then
		\begin{equation}\label{eq:subcritical-cutoff-gradient}
			\eta_{\e,\kappa}=0
			\quad\text{on }
			\partial B_{Ra_\e}(x_{\e,\kappa}),
			\qquad
			\eta_{\e,\kappa}=1
			\quad\text{on }
			\partial B_{2Ra_\e}(x_{\e,\kappa}),\quad
			|\nabla\eta_{\e,\kappa}|
			\leq
			\frac{C}{a_\e}.
		\end{equation}
		
		The cutoff field
		$\eta_{\e,\kappa}Z_\e$ has the required traces on the two
		boundaries of the annulus, but it is not divergence-free. Indeed,
		\[
		\nabla\cdot(\eta_{\e,\kappa}Z_\e)
		=
		\nabla\eta_{\e,\kappa}\cdot Z_\e,
		\]
		because $\nabla\cdot Z_\e=0$.
		Set
		\[
		f_{\e,\kappa}
		:=
		\nabla\eta_{\e,\kappa}\cdot Z_\e
		\quad\text{in }A_{\e,\kappa}.
		\]
		We verify that $f_{\e,\kappa}$ has zero mean. Using
		$\nabla\cdot Z_\e=0$, the divergence theorem, and the boundary
		values of $\eta_{\e,\kappa}$, we obtain
		\begin{align*}
			\int_{A_{\e,\kappa}}
			f_{\e,\kappa}\,\md x
			&=
			\int_{A_{\e,\kappa}}
			\nabla\cdot
			(\eta_{\e,\kappa}Z_\e)
			\,\md x
			=
			\int_{\partial B_{2Ra_\e}(x_{\e,\kappa})}
			Z_\e\cdot\bn\,\md S
			=
			\int_{B_{2Ra_\e}(x_{\e,\kappa})}
			\nabla\cdot Z_\e\,\md x
			=0.
		\end{align*}
		Thus, $f_{\e,\kappa}
		\in
		L^2_0(A_{\e,\kappa})$.
		Let
		\[
		\mathcal B_{A_{\e,\kappa}}:
		L^2_0(A_{\e,\kappa})
		\longrightarrow
		H^1_0(A_{\e,\kappa})^3
		\]
		be the Bogovski\u{\i} operator obtained by translation and
		rescaling from a fixed Bogovski\u{\i} operator on
		$B_{2R}\setminus\overline{B_R}$; see, for example,
		\cite[Lemma~III.3.1]{galdi2011}. Define
		\[
		c_{\e,\kappa}
		:=
		\mathcal B_{A_{\e,\kappa}}
		[f_{\e,\kappa}].
		\]
		Then
		\[
		c_{\e,\kappa}
		\in
		L^2(S;H^1_0(A_{\e,\kappa})^3),
		\qquad
		\nabla\cdot c_{\e,\kappa}
		=
		f_{\e,\kappa}
		\quad\text{in }A_{\e,\kappa}.
		\]
		The operator is chosen linearly and independently of $t$.
		
		By rescaling the Bogovski\u{\i} estimate from the fixed reference
		annulus, we have
		\[
		\|\nabla c_{\e,\kappa}\|_
		{L^2(A_{\e,\kappa})}
		\leq
		C
		\|f_{\e,\kappa}\|_
		{L^2(A_{\e,\kappa})},
		\]
		where $C$ is independent of $\e$ and $\kappa$. In view of
		\eqref{eq:subcritical-cutoff-gradient},
		\begin{equation}\label{eq:local-Bogovskii-estimate}
			\|\nabla c_{\e,\kappa}\|_
			{L^2(A_{\e,\kappa})}
			\leq
			\frac{C}{a_\e}
			\|Z_\e\|_
			{L^2(A_{\e,\kappa})}.
		\end{equation}
		The rescaled Poincar\'e inequality also gives
		\begin{equation}\label{eq:local-Bogovskii-L2-estimate}
			\|c_{\e,\kappa}\|_
			{L^2(A_{\e,\kappa})}
			\leq
			Ca_\e
			\|\nabla c_{\e,\kappa}\|_
			{L^2(A_{\e,\kappa})}
			\leq
			C
			\|Z_\e\|_
			{L^2(A_{\e,\kappa})}.
		\end{equation}
		
		Define $\psi_1^\e$ on $\O$ by
		\[
		\psi_1^\e(t,x)
		:=
		\begin{cases}
			0,
			&
			x\in B_{Ra_\e}(x_{\e,\kappa}),
			\\[2mm]
			\eta_{\e,\kappa}(x)Z_\e(t,x)
			-c_{\e,\kappa}(t,x),
			&
			x\in A_{\e,\kappa},
			\\[2mm]
			Z_\e(t,x),
			&
			x\notin
			\displaystyle
			\bigcup_{\kappa\in\Kc_\e}
			B_{2Ra_\e}(x_{\e,\kappa}).
		\end{cases}
		\]
		The definition is unambiguous because the balls are mutually
		disjoint.
		On the inner boundary of each annulus,
		\[
		\eta_{\e,\kappa}Z_\e-c_{\e,\kappa}=0,
		\]
		whereas on the outer boundary,
		\[
		\eta_{\e,\kappa}Z_\e-c_{\e,\kappa}=Z_\e.
		\]
		Thus the traces of the three pieces agree, and $\psi_1^\e
		\in
		L^2(S;H^1_0(\O)^3)$.
		
		In every annulus,
		\begin{align*}
			\nabla\cdot\psi_1^\e
			&=
			\nabla\cdot
			(\eta_{\e,\kappa}Z_\e-c_{\e,\kappa})
			=
			\nabla\eta_{\e,\kappa}\cdot Z_\e
			+
			\eta_{\e,\kappa}\nabla\cdot Z_\e
			-
			\nabla\cdot c_{\e,\kappa}
			=
			f_{\e,\kappa}-f_{\e,\kappa}
			=0.
		\end{align*}
		The field is also divergence-free inside the inner balls and
		outside the outer balls. Since the traces agree across the
		interfaces, no distributional surface term appears. Hence $\psi_1^\e
		\in
		L^2(S;H^1_{0,\div}(\O)^3)$.
		
		Moreover, $T_{\e,\kappa}
		\subset
		B_{Ra_\e}(x_{\e,\kappa})$,
		and $\psi_1^\e$ vanishes identically in this ball. Therefore,
		\[
		v_\e
		:=
		\psi_1^\e{}_{|\O_p^\e}
		\in
		L^2(S;H^1_{0,\div}(\O_p^\e)^3)
		\subset
		L^2(S;\mathbf H^1_{\div}(\O_p^\e)).
		\]
		
		We next show that the modification of $Z_\e$ tends strongly to
		zero. Set
		\[
		\mathscr N_\e
		:=
		\bigcup_{\kappa\in\Kc_\e}
		B_{2Ra_\e}(x_{\e,\kappa}).
		\]
		Since $\#\Kc_\e\leq C\e^{-3}$ and the balls are mutually
		disjoint,
		\begin{equation}\label{eq:small-modification-set}
			|\mathscr N_\e|
			\leq
			C\frac{a_\e^3}{\e^3}.
		\end{equation}
		
		Using
		\eqref{eq:subcritical-cutoff-gradient},
		\eqref{eq:local-Bogovskii-estimate},
		\eqref{eq:local-Bogovskii-L2-estimate}, and the disjointness of
		the balls, we obtain
		\begin{align}
			\|\psi_1^\e-Z_\e\|_
			{L^2(S;H^1(\O)^3)}
			\leq
			C\Bigg(
			&
			\|\nabla Z_\e\|_
			{L^2(S\times\mathscr N_\e)}
			+
			\frac{1}{a_\e}
			\|Z_\e\|_
			{L^2(S\times\mathscr N_\e)}
			\Bigg).
			\label{eq:local-restriction-error}
		\end{align}
		
		Fix $q>3$. By H\"older's inequality,
		\eqref{eq:small-modification-set}, and
		\eqref{eq:subcritical-Ze-regularity},
		\begin{align*}
			\|\nabla Z_\e\|_
			{L^2(S\times\mathscr N_\e)}
			&\leq
			|\mathscr N_\e|^{\frac12-\frac1q}
			\|\nabla Z_\e\|_
			{L^2(S;L^q(\O)^{3\times3})}
			\leq
			C_q
			\left(
			\frac{a_\e^3}{\e^3}
			\right)^{\frac12-\frac1q}
			\longrightarrow0.
		\end{align*}
		Furthermore,
		\begin{align*}
			\frac{1}{a_\e}
			\|Z_\e\|_
			{L^2(S\times\mathscr N_\e)}
			&\leq
			\frac{|\mathscr N_\e|^{1/2}}{a_\e}
			\|Z_\e\|_
			{L^2(S;L^\infty(\O)^3)}
			\leq
			C
			\left(
			\frac{a_\e}{\e^3}
			\right)^{1/2}.
		\end{align*}
		Since $a_\e=\e^\alpha$ and $\alpha>3$,
		\[
		\left(
		\frac{a_\e}{\e^3}
		\right)^{1/2}
		=
		\e^{(\alpha-3)/2}
		\longrightarrow0.
		\]
		Therefore,
		\begin{equation}\label{eq:subcritical-strong-restriction}
			\|\psi_1^\e-Z_\e\|_
			{L^2(S;H^1(\O)^3)}
			\longrightarrow0.
		\end{equation}
		
		Since the cutoff functions and the Bogovski\u{\i} operators are
		independent of $t$, and since the Bogovski\u{\i} operators are
		linear,
		\[
		\partial_t c_{\e,\kappa}
		=
		\mathcal B_{A_{\e,\kappa}}
		\left[
		\nabla\eta_{\e,\kappa}
		\cdot\partial_tZ_\e
		\right].
		\]
		The same argument, applied to $\partial_tZ_\e$, gives
		\begin{equation}\label{eq:subcritical-time-restriction}
			\|\partial_t\psi_1^\e-\partial_tZ_\e\|_
			{L^2(\O_T)^3}
			\longrightarrow0.
		\end{equation}
		
		Let $\widetilde v_\e$ denote the zero extension of $v_\e$ through
		the holes. Since $\psi_1^\e$ vanishes in a neighborhood of every
		hole,
		\[
		\widetilde v_\e=\psi_1^\e
		\quad\text{in }\O,\quad\text{and}\quad
		\nabla\widetilde v_\e
		=
		\nabla\psi_1^\e,
		\qquad
		\partial_t\widetilde v_\e
		=
		\partial_t\psi_1^\e.
		\]
		For this no-slip sequence, we use the zero extension in the
		definition of $\Te^\ast$. Equivalently, we use the consistency
		property of the chosen extension operator on fields having zero
		trace on the holes. Thus,
		\[
		\Te^\ast(v_\e)=\Te(\psi_1^\e),\quad
		\Te^\ast(\nabla v_\e)
		=
		\Te(\nabla\psi_1^\e),
		\qquad
		\Te^\ast(\partial_tv_\e)
		=
		\Te(\partial_t\psi_1^\e).
		\]
		
		By the continuity of the periodic unfolding operator,
		\eqref{eq:subcritical-strong-restriction} gives
		\[
		\begin{aligned}
			\|\Te^\ast(v_\e)-\Te(Z_\e)\|_
			{L^2(\O_T\times Y)^3}
			&\to0,\quad
			\|\Te^\ast(\nabla v_\e)-\Te(\nabla Z_\e)\|_
			{L^2(\O_T\times Y)^{3\times3}}
			\to0.
		\end{aligned}
		\]
		Likewise, \eqref{eq:subcritical-time-restriction} gives
		\[
		\|\Te^\ast(\partial_tv_\e)
		-\Te(\partial_tZ_\e)\|_
		{L^2(\O_T\times Y)^3}
		\to0.
		\]
		Combining these convergences with
		\eqref{eq:subcritical-Ze-properties}, we obtain
		\begin{equation}\label{eq:subcritical-ve-convergence}
			\begin{aligned}
				\Te^\ast(v_\e)
				&\to
				\varphi_1
				&&\text{strongly in }
				L^2(\O_T\times Y)^3,
				\\
				\Te^\ast(\nabla v_\e)
				&\to
				\nabla_x\varphi_1
				+
				\nabla_y\widehat\varphi_1
				&&\text{strongly in }
				L^2(\O_T\times Y)^{3\times3},
				\\
				\Te^\ast(\partial_tv_\e)
				&\to
				\partial_t\varphi_1
				&&\text{strongly in }
				L^2(\O_T\times Y)^3.
			\end{aligned}
		\end{equation}
		
		Finally, define $\varphi_1^\e
		:=
		\sqrt{\lambda^\e}\,v_\e$.
		By \eqref{eq:subcritical-strong-restriction} and the uniform
		boundedness of $Z_\e$ in $L^2(S;H^1_0(\O)^3)$,
		\begin{align*}
			\frac{1}{\sqrt{\lambda^\e}}
			\|\varphi_1^\e\|_
			{L^2(S;H^1(\O_p^\e)^3)}
			&=
			\|v_\e\|_
			{L^2(S;H^1(\O_p^\e)^3)}
			\leq
			\|\psi_1^\e\|_
			{L^2(S;H^1(\O)^3)}\\
			&\leq
			\|Z_\e\|_
			{L^2(S;H^1(\O)^3)}
			+
			\|\psi_1^\e-Z_\e\|_
			{L^2(S;H^1(\O)^3)}
			\leq C.
		\end{align*}
		This proves \eqref{eq:test-subcritical-uniform-H1}.
		Since $\lambda^\e$ is independent of $t$ and $x$,
		\[
		\nabla\varphi_1^\e
		=
		\sqrt{\lambda^\e}\,\nabla v_\e,
		\qquad
		\partial_t\varphi_1^\e
		=
		\sqrt{\lambda^\e}\,\partial_tv_\e.
		\]
		Therefore, by the linearity of $\Te^\ast$,
		\[
		\frac{1}{\sqrt{\lambda^\e}}
		\Te^\ast(\varphi_1^\e)
		=
		\Te^\ast(v_\e),\quad
		\frac{1}{\sqrt{\lambda^\e}}
		\Te^\ast(\nabla\varphi_1^\e)
		=
		\Te^\ast(\nabla v_\e),
		\qquad
		\frac{1}{\sqrt{\lambda^\e}}
		\Te^\ast(\partial_t\varphi_1^\e)
		=
		\Te^\ast(\partial_tv_\e).
		\]
		The convergences in
		\eqref{eq:subcritical-ve-convergence} now give
		\eqref{eq:test-subcritical-convergence}.
	\end{proof}
}
	\begin{lemma}[Subcritical solenoidal restriction]\label{lem:solenoidal-restriction}
		Let $\alpha>3$. For every $\zeta\in C_{c,\div}^\infty(\O)^3$ there exists
		$\Rc^\e\zeta\in\GH^1_{\div}(\O_p^\e)$, vanishing in a neighbourhood of
		$\Gamma_s^\e$, such that its zero extension $\overline{\Rc^\e\zeta}$ belongs to
		$H^1_{0,\div}(\O)^3$ and
		\begin{equation}\label{SolenoidalRestrictionConvergence}
			\overline{\Rc^\e\zeta}\longrightarrow\zeta\quad\text{strongly in }H^1(\O)^3,
		\end{equation}
		\begin{equation}\label{SolenoidalRestrictionBound}
			\|\Rc^\e\zeta\|_{H^1(\O_p^\e)}\leq C_\zeta,
		\end{equation}
		with $C_\zeta>0$ independent of $\e$. The map $\zeta\mapsto \Rc^\e\zeta$ is linear
		and independent of time.
	\end{lemma}
	
	\begin{proof}
		Apply the cutoff--Bogovski\u{\i} construction of
		Lemma~\ref{lem:test-functions-subcritical} to the time-independent
		field $\zeta$, i.e.\ with $\phi_1=\zeta$ and $\widehat\phi_1=0$, and retain the
		resulting solenoidal field before the final scaling by $\sqrt{\lambda^\e}$; call
		it $\Rc^\e\zeta$. By construction $\Rc^\e\zeta\in\GH^1_{\div}(\O_p^\e)$ vanishes near
		$\Gamma_s^\e$, so its zero extension lies in $H^1_{0,\div}(\O)^3$, and the
		uniform bound proved there gives \eqref{SolenoidalRestrictionBound} with $C_\zeta$
		depending only on $\zeta$. The same estimates show that the local correction near
		the obstacles tends to zero in $H^1(\O)^3$, since $a_\e/\e^3=\e^{\alpha-3}\to0$,
		while the underlying macroscopic field converges to $\zeta$ by standard
		consistency; together these give \eqref{SolenoidalRestrictionConvergence}.
	\end{proof}

	
	\addcontentsline{toc}{section}{References}
	{\small
		\bibliographystyle{plain}
		\bibliography{ref}
	}

\end{document}